\documentclass[9pt]{article}
\usepackage{amsmath,amssymb,amsthm,amsfonts,mathrsfs,latexsym,
times,color,epsfig,bm}
\usepackage{hyperref}

\newtheorem{thm}{Theorem}
\newtheorem{defn}{Definition}

\newtheorem{lem}{Lemma}

\newtheorem{rmk}{Remark}

\newenvironment{pf}{{\noindent \it \bf Proof.}}{{\hfill$\Box$}\\}

\newcommand{\vertiii}[1]{{\left\vert\kern-0.25ex\left\vert\kern-0.25ex\left\vert #1 
    \right\vert\kern-0.25ex\right\vert\kern-0.25ex\right\vert}}

\title{An iterative projection method for unsteady Navier-Stokes equations with high Reynolds numbers}
\date{}
\author{Xiaoming Zheng, Kun Zhao, Jiahong Wu, Weiwei Hu, Dapeng Du}

\author{Xiaoming Zheng\footnote{Department of Mathematics, Central Michigan University, Mount Pleasant, MI 48858. Corresponding Author. Email: zheng1x@cmich.edu}, 
Kun Zhao\footnote{School of Mathematical Sciences, Harbin Engineering University, Harbin 150001, Heilongjiang, China \& Department of Mathematics, Tulane University, New Orleans, LA 70118, USA. Email: kzhao@tulane.edu}, 
Jiahong Wu\footnote{Department of Mathematics,  University of Notre Dame, 
Notre Dame, IN 46556. Email: jwu29@nd.edu},
Weiwei Hu\footnote{Department of Mathematics, University of Georgia, Athens, GA 30602. Email: Weiwei.Hu@uga.edu},
Dapeng Du\footnote{School of Mathematics and Statistics, Northeast Normal University, Changchun, Jilin Province, 130024, China P.R..  Email: dudp954@nenu.edu.cn}
}

\begin{document}
%\showkeys

\maketitle
%\tableofcontents

\abstract{A new iterative projection method is proposed to solve the unsteady Navier-Stokes equations with high Reynolds numbers. The convectional projection method attempts to project the intermediate velocity to the divergence free space only once per time step. However, such a velocity is not genuinely divergence free in general practice, which can yield large errors when the Reynolds number is high. The new method has several important features: the BDF2 time discretization, the skew-symmetric convection in a semi-implicit form, two modulating parameters, and the iterative projections in each time step. A major difficulty in the proof of iteration convergence is the nonlinear convection. We solve this problem by first analyzing the non-convective scheme with a focus on the spectral properties of the iterative matrix, and then employing a delicate perturbation analysis for the convective scheme. The work achieves the weakly divergence free velocity (strongly divergence free for divergence free finite element spaces), and the rigorous stability and error analysis when the iterations converge. The three dimensional numerical tests confirm that this new method can effectively treat high Reynolds numbers with only a few iterations per time step, where the convectional projection method and the iterative projection method with the explicit convection would fail. 
 } 

\vspace{-0.5cm}

\section{Introduction}
The incompressible Navier-Stokes equations are predominant in viscous fluid dynamics, including problems with high Reynolds numbers such as turbulence flows. However, it is extremely challenging to design efficient and robust numerical schemes for these problems. A major difficulty is the treatment of the incompressibility constraint,
or the divergence free condition, which couples the velocity and pressure.
Consider a bounded {three-dimensional} domain $\Omega$ and the Navier-Stokes equations
\begin{align}
{u}_t - \nu \Delta u + \text{NL}(u,u)  + \nabla p &= f,
\label{ns001}
\\
\nabla\cdot{u} &= 0,
\label{ns002}
\end{align}
with the nonlinear convection $\text{NL}(u,v)=(u\cdot \nabla)v$,  the Dirichlet boundary condition 
$u|_{\partial\Omega} = 0$,
where $u=(u_1, u_2, u_3)$ is the velocity,  $p$ is  the kinematic pressure,  $\nu$ is the kinematic viscosity, and $f$ is the force. There are various equivalent forms of the convection term (e.g., see \cite{CHARNYI2017289}), and this work focuses on the skew-symmetric form
\begin{align}
\text{NL}(u,v) = (u\cdot\nabla) v + \frac12 (\nabla\cdot u) v.\label{convection_skew}
\end{align}

The projection method is a classic and popular method to solve the Navier-Stokes equations first developed in \cite{Chorin1968, Temam1968} in 1968, where a systematic review can be found in \cite{GuermondShenReview2006}. Its crucial operation is projecting the intermediate velocity, which is not divergence free, to the divergence free space in a certain sense. 
One of the most recommended projection schemes in \cite{GuermondShenReview2006}
can be put as follows, when the BDF2 (Backward Differentiation Formula 2) time integration is employed,
\begin{align}
(1.5-k\nu\Delta) u^{n+1,*} + k \text{NL}(w^{n+1}, u^{n+1,*})  &=-k\nabla \tilde{p}^{n+1}+ F^{n+1}, 
\label{pp001}\\
-\Delta \phi^{n+1} &= -\frac{1}{k} \nabla\cdot u^{n+1,*},
\label{pp002}\\
p^{n+1} &= \tilde{p}^{n+1}  +1.5\phi^{n+1} - \nu  \nabla\cdot { u}^{n+1,*}, 
\label{pp003}\\
u^{n+1} &= u^{n+1,*} - k \nabla \phi^{n+1},
\label{pp004}
\end{align}
where 
$w^{n+1}=2u^n-u^{n-1}$, 
$\tilde{p}^{n+1}=2p^n-p^{n-1}$,
$F^{n+1}=2{ u}^n-0.5{ u}^{n-1} + k f^{n+1}$, 
$u^{n+1,*}|_{\partial\Omega}=0$, 
$\frac{\partial\phi^{n+1}}{\partial{ n}}\Big|_{\partial\Omega} = 0$.
The variable $\phi^{n+1}$ is from the Hodge decomposition $u^{n+1,*}= u^{n+1} + k \nabla \phi^{n+1}$, where $u^{n+1}$ is supposed to be divergence free. 
This method was first proposed in \cite{Timmermans1996} in 1996 and a similar one was  introduced in \cite{Brown2001}. 
We call this method ``{\em rotational projection method}'' in this paper. 
Another popular projection scheme replaces the pressure update  \eqref{pp003} with 
\begin{equation}
p^{n+1} = \tilde{p}^{n+1}  +1.5\phi^{n+1}.  \label{standard_proj}
\end{equation}
We  call it ``{\em standard projection method}'' in this work.  This one is also widely used in literature, such as \cite{Guermond2002} and \cite{ALBENSOEDER2005536}.

The main strength of the projection method lies in its simplicity and efficiency. It decouples the velocity and pressure fields in the Navier-Stokes equations and splits the original system into two smaller standard problems: one advection-diffusion equation for $u^{n+1,*}$ and one Poisson equation for $\phi^{n+1}$, followed by two updates to get the end-of-step velocity and pressure. Therefore, its computational cost is far less than that of any iterative methods including the Uzawa and Newton methods. 
However,  it has some serious defects. First, the velocity field obtained is not divergence free even in the weak sense in the mixed finite element implementations (see \cite[Remark 3.5]{Guermond1998} and Appendix\,\ref{sec_non_divfree}), and even when the finite element method admits divergence free velocity field (see Definition\,\ref{def_divfreeFEM} and  \cite{SchroederLube2017}). This would induce mass loss when the 
velocity is used to transport a density function. 
Second, the stability proof for the BDF2 schemes is not available. But we note that a recent work \cite{HuangShen2023} obtains stability with a general second order BDF scheme and a generalized scalar auxiliary variable approach. 
Third, this scheme is unable to treat high Reynolds numbers due to the violation of the divergence free condition  (see Section\,\ref{sec_Re2300}).  

We propose the following iterative projection method, with the iteration index $s=0,1,\cdots$,
\begin{align}
(1.5 -k\nu \Delta) u^{n+1,s} 
+ k \text{NL}(u^{n+1,s-1}, u^{n+1,s})
& = - k\nabla p^{n+1,s} + F^{n+1},
 \label{mol001}\\
-\Delta \phi^{n+1,s} &= -\frac{1}{k} \nabla\cdot u^{n+1,s},
\label{mol002}\\
p^{n+1,s+1} &= {p}^{n+1,s} + \alpha \phi^{n+1,s} - \rho  \nabla\cdot u^{n+1,s}, 
\label{mol003}
\end{align}
with $u^{n+1,-1}=2u^n-u^{n-1}$, $p^{n+1,0}=2p^n-p^{n-1}$, 
and boundary conditions $u^{n+1,s}|_{\partial\Omega}=0$ and  
$\frac{\partial\phi^{n+1,s}}{\partial{ n}}\Big|_{\partial\Omega} = 0$. 
The control parameters $\alpha$ and $\rho$ are  used to optimize the convergence speed.  Note when $\alpha=1.5$ and $\rho=\nu$, it is the iteration of the rotational projection scheme. When $\alpha=1.5$ and $\rho=0$, it is the iteration of the standard projection scheme. When $\alpha=0$, it is the Uzawa method  (see \cite[Chapter 2]{Fortin1983}).
When the iterations converge, the limit $(u^{n+1}, p^{n+1})$ is the solution of the following nonlinear coupled scheme,
\begin{align}
(1.5-k\nu \Delta ) u^{n+1}
+ k \text{NL}(u^{n+1}, u^{n+1}) 
 + k\nabla p^{n+1} &= F^{n+1}, 
\label{eqn01}\\
\nabla\cdot u^{n+1} &= 0.\label{eqn02}
\end{align}

The iterative projection approach with the fully explicit convection has been utilized in some previous work including  \cite{Zheng2005,Doering2018,Aoussou2018}. 
In \cite{Zheng2005}, it was found that the repeated projections can reduce the spurious errors generated from the singular surface tension forces in the free boundary problems. 
In \cite{Doering2018},  the iterative projections provide accurate velocity fields in the stratification of temperature in the Boussinesq flows.  
In \cite{Aoussou2018}, it is found that iterations could improve splitting errors of projection method. 
However, the numerical test in  Section\,\ref{sec_IMEX}  shows the iterative projections with explicit convection is not so stable as the implicit treatment. 
It is worth noting that \cite{Aoussou2018} briefly mentioned an iterative projection method, (64b) in \cite{Aoussou2018}, without providing analysis or simulation. This method employs the original form of convection, $(u^{n+1,s-1}\cdot \nabla) u^{n+1,s}$, rather than the skew-symmetric form, $\text{NL}(u^{n+1,s-1}, u^{n+1,s})$,  in \eqref{mol001}. The skew-symmetric form significantly simplifies the stability and error estimation proofs for the limit scheme, as the 
$L^2$ inner product $(\text{NL}(w,v),v)$ vanishes when 
$w,v\in H^1_0(\Omega)$. In contrast, the original convection form lacks this advantage.

In Section\,\ref{sec_analysis_iterations}, the convergence of the proposed iterative projection method to the limit scheme \eqref{eqn01} and \eqref{eqn02} at each individual time step is presented through the normal mode and finite element analysis. 
Note the convection term  in the iterative method \eqref{mol001}, $\text{NL}(u^{n+1,s-1}, u^{n+1,s})$, is explicit on the first component and implicit on the second, while the convection in the limit scheme \eqref{eqn01} is fully implicit. 
These implicit convection terms causes the main challenge in analysis. We first analyze the case without convection through the eigenvalue study of the iterative matrix for the pressure.  The convective case is then regarded as a perturbed system of the non-convective system, where the time step size $k$ is the perturbation parameter. A delicate induction process is implemented to establish the iteration convergence.

In Section\,\ref{sec_stabilityanderror}, we provide the stability and  error analysis results of the limit scheme in the context of the mixed finite element method, along with a brief literature review of the stability and error analysis of the fully implicit schemes with Galerkin finite element method for the Navier-Stokes equations.
In Section\,\ref{sec_numerical_tests}, we test the proposed scheme with three dimensional finite element method and high Reynolds numbers.  The conclusions and discussions are given in Section\,\ref{sec_conclusions}.

\vspace{-0.4cm}
 
\section{Convergence of projection iterations at a single time step}
\label{sec_analysis_iterations}
This section is devoted to the study of the convergence  of the iterative projection method \eqref{mol001}, \eqref{mol002}, \eqref{mol003} to \eqref{eqn01} and \eqref{eqn02} when $s$ goes to infinity. To simplify notations, we delete the time step superscript $(n+1)$  in this section.

\vspace{-0.4cm}

\subsection{Normal mode analysis without convection} \label{sec_normalmode}
To perform the normal mode analysis, we neglect the convection term  and assume the solution is smooth and periodic in the region $\Omega=[0,2\pi]^3$. 
Note that the intermediate value $\phi^s$ in \eqref{mol002} can be written as
$
\phi^s = -\frac{1}{k} (-\Delta)^{-1} (\nabla\cdot { u}^s), 
$
where $(-\Delta)^{-1}$ refers to the inverse operator of the  Neumann problem \eqref{mol002}.
So the iterative scheme \eqref{mol001}, \eqref{mol002}, \eqref{mol003} when the convection is removed can be simplified to, for $s=0, 1, \cdots$, 
\begin{align}
1.5 { u}^{s} 
- k\nu \Delta { u}^s + k\nabla p^{s} &= F,
\hspace{.3 in} { u}^s|_{\partial\Omega}=0, 
 \label{weiwei001}\\
p^{s+1} &= {p}^{s} -\big(\rho +  \frac{\alpha}{k} (-\Delta)^{-1} \big)  (\nabla\cdot { u}^{s}). 
\label{weiwei002}
\end{align}

Consider the solution $({u},p)$ satisfying \eqref{eqn01} and \eqref{eqn02} without the convection, and the sequence $({u}^s,p^s)$ of \eqref{weiwei001} and \eqref{weiwei002}.
Let $\bar{u}^s={u}^s-{u}$ and $\bar{p}^s=p^s-p$. Then
\begin{align}
1.5 \bar{{u}}^s - k\nu \Delta \bar{u}^s + k \nabla\bar{p}^s &=0,
\label{pup01}\\
\bar{p}^{s+1} &=\bar{p}^{s} - \big(\rho + \frac{\alpha}{k} (-\Delta)^{-1}\big) (\nabla\cdot \bar{u}^{s}). \label{pup02}
\end{align}
Denote $\bar{u}^s=(\bar{u}^s_1, \bar{u}^s_2, \bar{u}^s_3)$, the multi-wavenumber $\xi=(\xi_1, \xi_2, \xi_3)\in\mathbb{Z}^3$, 
$\xi\cdot x=\xi_1 x_1 + \xi_2 x_2 + \xi_3 x_3$.  
Denote the Fourier series of $\bar{u}^s_j$ and $\bar{p}^s$ as 
$\bar{u}_j^s(x) = \sum_{\xi\in \mathbb{Z}^3} \widehat{\bar{u}_j^s}(\xi) e^{i\xi\cdot x}$, $j=1,2,3$, 
$\bar{p}^s(x) = \sum_{\xi\in\mathbb{Z}^3} \widehat{\bar{p}^s}(\xi) e^{i\xi\cdot x}$, where $i=\sqrt{-1}$.
Then \eqref{pup01} and \eqref{pup02} become
\begin{align}
1.5 \widehat{\bar{u}_j^s}(\xi) + k\nu |\xi|^2 \widehat{\bar{u}_j^s}(\xi) 
+ ki\xi_j \widehat{\bar{p}^s}(\xi)&=0,\hspace{.3 in} j=1,2,3; 
\label{mop001}
\\
\widehat{\bar{p}^{s+1}}(\xi) &= 
\widehat{\bar{p}^{s}}(\xi) - 
\Big(
 \rho + \frac{\alpha}{k|\xi|^2} \Big) 
 \sum_{j=1}^3 i\xi_j \widehat{\bar{u}^s_j}(\xi).
 \label{mop002}
\end{align}
It can be calculated from \eqref{mop001} that 
$
\widehat{\bar{u}^s_j}(\xi) = \frac{-ki\xi_j \widehat{\bar{p}^s}}{1.5+k\nu |\xi|^2}
$, which is substituted into \eqref{mop002} to achieve
\begin{equation}
\widehat{\bar{p}^{s+1}}(\xi) =  
\frac{ (1.5-\alpha) + k(\nu-\rho)|\xi|^2}{1.5+k\nu|\xi|^2}
\cdot \widehat{\bar{p}^{s}}(\xi)
\triangleq C(\alpha,\rho, \xi)  
\cdot \widehat{\bar{p}^{s}}(\xi).
\end{equation}
For the sequence $\widehat{\bar{p}^{s}}$ to converge to zero, it is equivalent to set $|C(\alpha,\rho, \xi)|<1$.
Note the optimal convergence occurs when $\alpha=1.5$ and $\rho=\nu$, that is, the constant $C(\alpha,\rho, \xi)=0$. This corresponds to the iterative rotational projection method.  Therefore, this method obtains the exact solution in just one iteration. However, this analysis is based on the smooth solutions and without convection, which is not the case for the full Navier-Stokes equations with the convection term and  the finite element solutions of limited regularity.

If $\alpha=1.5$, then $|C(\alpha,\rho, \xi)|<1$ leads to 
$
-\frac{1.5}{k |\xi|^2} < \rho < 2\nu + \frac{1.5}{k|\xi|^2}.
$
For these inequalities to hold for all $\xi\in\mathbb{Z}^3$, it is equivalent to have $0\le \rho\le 2\nu$.
If $\alpha=1.5$ and $\rho=0$, that is, when the iterative standard projection method is used, then 
$C(1.5, 0, \xi) =\frac{k\nu|\xi|^2}{1.5+k\nu|\xi|^2}.
$
In this case, the iteration always converges but the convergent constant $C$ is close to 0 when $k\nu|\xi|^2$ is small, and approaching 1 when $k\nu|\xi|^2$ is large. Thus, this iterative scheme is preferred only when $\nu$ is sufficiently small. 
If $\alpha=0$, then this scheme reduces to the Uzawa scheme and the convergence constant is 
$
C_u(\rho, \xi)=
\frac{1.5 + k(\nu-\rho) |\xi|^2}{1.5+k\nu|\xi|^2}.
$
The convergence requires
$ 0 < \rho \le 2\nu$. 
It is easy to see that $|C_u|\to 1$ when $\nu\to 0$.  
Thus, the Uzawa iterations would be applicable only when $\nu$ is sufficiently large. This is consistent with the fact that the Uzawa method is mainly used for the Stokes equations whose viscosity is considerably large.
The above analysis is summarized in Table\,\ref{normal_mode}. 
\begin{table}[htbp]
\scriptsize
\begin{center} 
\caption{\scriptsize Normal mode analysis when the convection is neglected. The  pressure update is $p^{s+1}={p}^{s} -(\rho +  \frac{\alpha}{k} (-\Delta)^{-1} )  (\nabla\cdot { u}^{s})$.
\label{normal_mode}
}
\begin{tabular}{|p{1.7cm} | p{2cm} | p{2cm}|p{2cm}|p{4cm}|}
\hline
method & parameters  & Convergence constant & parameter range for convergence &notes\\
\hline
iterative rotational projection & $\alpha=1.5$, $\rho=\nu$
& 0  & always convergent & super-convergence\\
\hline
iterative standard projection & $\alpha=1.5$, $\rho=0$ &
$\frac{k\nu|\xi|^2}{1.5+k\nu|\xi|^2}$
& always convergent & convergence is fast when $\nu$ is small \\
\hline 
Uzawa & $\alpha=0$ &
$\frac{1.5 + k(\nu-\rho) |\xi|^2}{1.5+k\nu|\xi|^2}
$
& $0<\rho\le 2\nu$ & convergence is fast when $\nu$ is large.\\
\hline
\end{tabular}
\end{center}
\end{table}

%%%%%%%%%%%%%%%%%%%%%%%%%%%%%%%5

\vspace{-0.8cm}
%%%%%%%%%%%%%%%%%%%%%%%%%%%
\subsection{Iteration analysis with mixed finite element method}
\label{sec_iterative_analysis}
The paper \cite{Aoussou2018} has proved the convergence of iterative projections when there is no convection or the convection term is treated utterly explicitly, and there are no freely chosen parameters. In contrast, our scheme handles the convection implicitly and has two parameters $\alpha$ and $\rho$ that can be modulated to obtain optimal convergence.

We denote $L^2(\Omega)$ as the Lebesgue space of {square} integrable functions with inner product $(u,v) = \int_{\Omega} u(x) v(x) \mathrm{d}x$,
and the $L^2$ norm as $\|u\|_{L^2(\Omega)} = \sqrt{(u,u)}$. Denote 
$H^1_0(\Omega) = \{f: f, \frac{\partial f}{\partial x_i}\in L^2(\Omega), i=1,2, 3, f|_{\partial\Omega}=0 \}$. We introduce 
\begin{align}
V&= \{u=(u_1, u_2, u_3): u_i \in H^1_0(\Omega), i=1,2,3 \} ,
\label{def_V}\\
Q&=L^2(\Omega)/\mathbb{R} =\{q\in L^2(\Omega): \frac{1}{|\Omega|} \int_\Omega q(x) \mathrm{d} x=0 \},
\label{def_Q}
\end{align}
with the norm $|u|_1 = (\int_{\Omega} |\nabla u(x)|^2 \mathrm{d}x)^{1/2}$,
where $|\nabla u|^2=\sum_{i,j=1}^3 |\frac{\partial u_i}{\partial x_j}|^2 
$ and $\frac{1}{|\Omega|} \int_\Omega q(x) \mathrm{d}x $ is the average of $q$ over $\Omega$. Let $V_h\subset V$ and $Q_h\subset Q$ be a pair of conforming mixed finite element spaces, 
which satisfies the inf-sup condition (a.k.a. the Ladyzhenskaya-Babuka-Brezzi (LBB) condition, see e.g. \cite{Fortin1983}), i.e.,
\begin{equation}
\inf_{q\in Q_h} \sup_{v\in V_h}
\frac{(\nabla\cdot v, q)}{\|q\|_{L^2(\Omega)/\mathbb{R}} |v|_1} > 0.
\label{inf_sup}
\end{equation}
Denote the basis of $V_h$ as $w_i, i=1, \cdots, l$, and the basis of $Q_h$ as $q_j, j=1,\cdots, m$.

\begin{defn}
\label{def_divfreeFEM}
A velocity field $u_h\in V_h$ is weakly divergence free if $(\nabla\cdot u_h, q_h)=0$ for all $q_h\in Q_h$ and strongly divergence free if $\nabla\cdot u_h$ is almost everywhere zero in $\Omega$. 
A pair of mixed finite element spaces $V_h$ and $Q_h$ is called  divergence free if $\nabla\cdot V_h\subset Q_h$, where $\nabla\cdot V_h$ is the range of the divergence of $V_h$. Furthermore, a finite element method with a pair of divergence free finite element spaces is called a div-free FEM, and otherwise a non-div-free FEM.
\end{defn}

In this work, we use a weak incompressibility measure \eqref{weak_divfreemeasure} (a strong incompressiblity measure \eqref{eqn_strong_measure})  to check whether a velocity is weakly (strongly) divergence free. A weakly divergence free velocity may be not pointwise divergence free, as seen in Section\,\ref{sec_test_iterations}.
However, under a pair of divergence free FEM spaces, a weakly divergence free velocity is also strongly divergence free. Indeed, if $u_h$ satisfies $(\nabla\cdot u_h, q_h)=0$ for all $q_h\in Q_h$ and $\nabla\cdot V_h\subset Q_h$, then $\|\nabla\cdot u_h\|_{L^2(\Omega)}=0$ by choosing $q_h=\nabla\cdot u_h$, which suggests  
$\nabla\cdot u_h$ is zero almost everywhere in $\Omega$. 
One popular example of divergence free elements is the Scott–Vogelius pair \cite{ScottVogelius1984}. However, many $H^1$-conforming, inf-sup stable mixed finite elements pairs are not  divergence free, including the very popular Tayor-Hood elements.

%%%%%%%%%%%%%%%%%%%%%%%%%%%%%%%%%%%%%%%
%%%%%%%%%%%%%%%%%%%%%%%%%%%%%%%%%%%%%%%
\vspace{-0.4cm}
%%%%%%%%%%%%%%%%%%%%%%%%%%%%%%%%%%%%%%%
\subsubsection{Iteration analysis without convection}
\label{sec_iter_no_convection}
Without  convection, the finite element solution
$(u^s, \phi^s, p^s)\in V_h\times Q_h\times Q_h$
of  \eqref{mol001}, \eqref{mol002} and 
\eqref{mol003} satisfies, after dropping time step index $(n+1)$, 
\begin{align}
1.5({u^s}, {w_i}) + k\nu (\nabla {u^s}, \nabla{w_i})
%+ k(({w}\cdot\nabla) {u^s}, {w_i})
- k(\nabla\cdot {w_i}, p^s) = (F,{w_i}),
&\qquad\forall\, i=1,\cdots, l, 
\label{ip001}\\
(\nabla \phi^s, \nabla q_i) =
-\frac{1}{k}(\nabla \cdot u^s, q_i),  &\qquad\forall\, i=1,\cdots, m, 
\label{ip002} \\
(p^{s+1}, q_i) = (p^s, q_i) + \alpha (\phi^s, q_i) - \rho (\nabla\cdot{ u}^s, q_i), &\qquad\forall\, i=1, \cdots, m.
\label{ip003}
\end{align}
The matrices $A_{l\times l}$, $B_{m\times l}$, 
$G_{m\times m}$, and $M_{m\times m}$ are introduced as follows, 
$A_{ij} = 1.5 (w_i, w_j) + k\nu (\nabla w_i, \nabla w_j), i,j=1,\cdots, l$;
$B_{ij} = -(q_i, \nabla\cdot w_j), i=1,\cdots, m, j=1,\cdots, l$;
${G}_{ij} = (\nabla q_i, \nabla q_j)$, 
$M_{ij} = (q_j, q_j), i,j=1, \cdots, m$.
The {transpose} of $B$ is written as $B^T$. 
Let the vector 
$\vec{f}=\left( 
\begin{array}{ccc}
(F,w_1) & \hdots  & (F,w_l)
\end{array}\right)^T$.
We express $u^s = \sum_{i=1}^l u^s_i w_i$, 
$\phi^s = \sum_{i=1}^m \phi^s_i q_i$, 
$p^s = \sum_{i=1}^m p^s_i q_i$ and 
denote the corresponding vector forms as
$
\vec{u}^{s} = \left(
\begin{array}{ccc}
u^s_1 & \hdots & u^s_l
\end{array} \right)^T$,
$\vec{\phi}^{s} = \left(
\begin{array}{ccc}
\phi^s_1 & \hdots &  \phi^s_m
\end{array} \right)^T$, 
$\vec{p}^{\,s} = \left(
\begin{array}{ccc}
p^s_1 & \hdots & p^s_m
\end{array}\right)^T.
$
Thereafter, the conversion between the function and  vector forms of a quantity will be used this way. 
Thus, the numerical scheme, 
\eqref{ip001}, \eqref{ip002}, and \eqref{ip003}, 
turns into the following matrix-vector form, 
\begin{align}
A\vec u^s + k B^t \vec{p}^{\,s} = \vec f, \quad 
G\vec{\phi}^s =\frac{1}{k} B\vec{u^s}, \quad 
M\vec{p}^{\,s+1} = M\vec{p}^{\,s} + \alpha M\vec{\phi}^s + \rho B\vec{u}^s.
\label{wei041}
\end{align}

Without the convection,  the weak solution of \eqref{eqn01} and \eqref{eqn02}  $(u,p)\in V_h\times Q_h$ satisfies
\begin{align}
1.5(u, {w_i}) + k\nu (\nabla u, \nabla{w_i})
%+ k(({w}\cdot\nabla) u, {w_i})
- k(\nabla\cdot {w_i}, p) = (F,{w_i}),
&\qquad\forall\, i=1,\cdots, l, 
\label{ip011}\\
(\nabla\cdot u, q_i) = 0, &\qquad\forall\, i=1, \cdots, m.
\label{ip013}
\end{align}
The corresponding vector form can be written as 
\begin{align}
A\vec u + k B^t \vec p = \vec f,\quad 
G\vec{\phi} =\frac{1}{k} B\vec{u}, \quad 
M\vec{p} = M\vec{p} + \alpha M\vec{\phi} + \rho B\vec{u}.
\label{wei044}
\end{align}
The last two equations in \eqref{wei044}  can be reduced to $\left(G + \frac{\alpha}{k \rho} M \right)\vec{\phi} = 0$, whose solution is $\vec{\phi}=0$ since both $G$ and $M$ are symmetric and positive definite. Thus, the last two equations in \eqref{wei044}  are equivalent to 
$B\vec{u}=\vec{0}$.

Denote the errors in the function form as 
$
{\bar{u}}^{s} = {u}^s - {u},
{\bar{\phi}}^{s} = {\phi}^s - {\phi},
{\bar{p}}^{s} = {p}^s - {p},
$
and the respective vector forms as 
$\vec{\bar{u}}^{s}$,
$\vec{\bar{\phi}}^{s}$,
$\vec{\bar{p}}^{\,s}$.
The subtractions of the equations,  
\eqref{wei041}-\eqref{wei044},  give rise to
\begin{align}
A\vec{\bar{u}}^s + k B^t \vec{\bar{p}}^{\,s} = 0,
\quad 
G \vec{\bar{\phi}}^s = \frac{1}{k} B \vec{\bar{u}}^s, 
\quad 
M\vec{\bar{p}}^{\,s+1} = M\vec{\bar{p}}^{\,s} + \alpha M\vec{\bar{\phi}}^s + \rho B \vec{\bar{u}}^s.
\label{kun001}
\end{align}
From \eqref{kun001} we get
$\vec{\bar{u}}^s = -k A^{-1}B^t \vec{\bar{p}}^{\,s}$ and 
$\vec{\bar{\phi}}^s = - G^{-1}BA^{-1} B^t \vec{\bar{p}}^{\,s}$.
Let 
$
D=BA^{-1}B^t,
$
which is the negative Schur complement of $A$ (e.g., see \cite{Benzi2005}), 
then 
\begin{equation}
\vec{\bar{p}}^{\,s+1} = ( I-(\alpha G^{-1} + \rho k M^{-1} )D ) 
\vec{\bar{p}}^{\,s}
\triangleq (I-K) \vec{\bar{p}}^{\,s},
\label{keyrelation01}
\end{equation}
where
\begin{equation}
K_{m\times m}=(\alpha G^{-1}  + \rho k M^{-1} )D.
\label{def_K}
\end{equation}

To study the convergence, we analyze the eigenvalues of the matrix $K$.
First, we introduce the following lemmas. 
\begin{lem}\label{lem1}
Let $S$, $T$ be symmetric matrices. If $S$ is positive definite, then $S^{-1}T$ is diagonalizable.
\end{lem}
\begin{pf}\\
Since $S$ is symmetric and positive definite (SPD), its inverse $S^{-1}$ is also SPD. Denote its square root as $S^{-1/2}$. Then $S^{1/2} S^{-1}TS^{-1/2}=S^{-1/2}TS^{-1/2}$,
where the latter is symmetric because $T$ and $S^{-1/2}$ are both symmetric. Therefore, $S^{-1}T$ is similar to a symmetric matrix and hence diagonalizable.\end{pf}
\begin{lem}
\label{lem_eigenvaluepositive}
All the eigenvalues of $K$ are positive and $K$ has a set of linearly independent eigenvectors that span $\mathbb{R}^m$, which is called an eigenbasis of $K$. 
\end{lem}
\begin{pf}\\
First, note that without the convection term the matrix $A$ is SPD. The inf-sup condition guarantees $B$ has the full row rank and then the matrix $D$ is also SPD. In addition, the matrices $G$ and $M$ are also SPD.

Suppose $\lambda$ is an eigenvalue of $K$ and $\vec{\xi}\in \mathbb{R}^m$ is an associated eigenvector, i.e., 
$K\vec{\xi}=\lambda\vec{\xi}$.   Left-multiplying both sides by $\vec{\xi}^{\,t}D$ leads to 
\begin{equation}
\alpha \vec{\xi}^{\,t} DG^{-1}D \vec{\xi}
+ \rho k \vec{\xi}^{\,t} DM^{-1}D \vec{\xi} = \lambda \vec{\xi}^{\,t} D \vec{\xi}.
\label{kun006}
\end{equation}
Since the quantities $\alpha$, $\rho$, $k$, $\vec{\xi}^{\,t} DG^{-1}D \vec{\xi}$, $\vec{\xi}^{\,t} DM^{-1}D \vec{\xi}$, $\vec{\xi}^{\,t} D \vec{\xi}$ are all positive, the value of $\lambda$ must be positive.

Furthermore, a matrix  is diagonalizable if and only if there is a set of eigenbasis of this matrix (e.g., \cite{Horn2013} Theorem 1.3.7, p.\,59).
In the matrix $K$, the left factor $(\alpha G^{-1} + \rho k M^{-1})$ is SPD and the right factor $D$ is symmetric.
According to Lemma\,\ref{lem1},  $K$ has an eigenbasis. 
\end{pf}

Denote an eigenbasis of $K$ as 
$\{ \vec{\xi}_1, \cdots, \vec{\xi}_m \}$ and the corresponding eigenvalues as $\lambda_j, j=1,\cdots, m$. 
Let $\rho(I-K)$ be the spectral radius of $I-K$, that is, 
$
\rho(I-K)=  \max_{j=1,\cdots, m} |1-\lambda_j|.
$
Because all the eigenvalues are positive according to Lemma\,\ref{lem_eigenvaluepositive},  it is easy to deduce  that 
\begin{lem}
\label{lem_rhoI-K}
$\rho(I-K)<1$ if and only if $\lambda_{\max}<2$.
\end{lem}
The next lemma gives an upper bound of the largest eigenvalue $K$ dependent on the parameters in the numerical scheme.
\begin{lem}
\label{lem_lambdamax}
Let $\lambda_{\max}$ be the largest eigenvalue of $K$, then 
\begin{equation}
\lambda_{\max}
\le  \max\left( \frac{\alpha}{1.5}, \frac{\rho}{\nu}\right).
\label{cat104}
\end{equation}
\end{lem}
\begin{pf}\\
Let $\lambda$ be an eigenvalue of $K$ and $\vec{\xi}$ be an associated eigenvector. According to \eqref{kun006}, 
\begin{equation}
\lambda = \frac{\alpha \langle G^{-1}D\vec{\xi}, D\vec{\xi}\rangle 
+ \rho k \langle M^{-1}D\vec{\xi}, D\vec{\xi}\rangle}{\langle D\vec{\xi}, \vec{\xi}\rangle},
\label{kun007}
\end{equation}
where the notation $\langle \vec{\xi}, \vec{q}\rangle\triangleq \vec{\xi}^{\,t} \vec{q}$ for any two vectors in the same Euclidean space $\mathbb{R}^m$ or $\mathbb{R}^l$.
Denote $\vec{v}=A^{-1}B^t \vec{\xi}$, i.e., $A\vec{v}=B^t\vec{\xi}$ or
$D\vec{\xi}=B\vec{v}$.
 Then \eqref{kun007} becomes
\begin{equation}
\lambda 
=\frac{\alpha \langle G^{-1}B\vec{v}, B\vec{v}\rangle 
+ \rho k \langle M^{-1} B\vec{v}, B\vec{v}\rangle}{\langle \vec{v}, A\vec{v}\rangle}
=\frac{\alpha \langle G^{-1}B\vec{v}, B\vec{v}\rangle
+ \rho k \langle M^{-1} B\vec{v}, B\vec{v}\rangle}
{1.5\|v\|^2_{L^2(\Omega)} + \nu k\|\nabla v\|^2_{L^2(\Omega)}}
\triangleq R(\vec{v}).
\label{kun008}
\end{equation}
This equates the eigenvalues to the stationary values of the Rayleigh quotient $R(\vec{v})$ defined in \eqref{kun008}.
To estimate the two numerator terms in \eqref{kun008}, 
let $\vec{\phi}=G^{-1}B\vec{v}$, 
then $G\vec{\phi}=B\vec{v}$. Using the connection between the vector and function forms: 
$\vec{\phi}=(\phi_1 \cdots  \phi_m )^t$ and $\phi=\sum_{i=1}^m \phi_i q_i $, 
$\vec{v}=(v_1 \cdots v_l )^t$ and $v=\sum_{i=1}^l v_i w_i$,
we obtain
\begin{equation}
\langle G^{-1}B\vec{v}, B\vec{v}\rangle = \langle \vec{\phi}, B\vec{v}\rangle 
 = -\int_\Omega \phi (\nabla\cdot v) \mathrm{d}x.
 \label{kun009b}
\end{equation}
Similarly, we observe that $G\vec{\phi}=B\vec{v}$  is equivalent to 
\begin{equation}
\int_\Omega\nabla \phi \cdot \nabla \eta \mathrm{d}x = -\int_\Omega (\nabla\cdot v) \eta \mathrm{d}x, \quad
\forall \eta \in Q_h.
\label{kun009}
\end{equation}
Using $\eta=\phi$ in \eqref{kun009}, integration by parts, and Young's inequality, we get
$
\|\nabla \phi\|_{L^2(\Omega)}^2 = 
-\int_\Omega \phi (\nabla\cdot v)\mathrm{d}x
=\int_\Omega v\cdot  \nabla \phi \mathrm{d}x
\le \frac{1}{2}\big( \|v\|^2_{L^2(\Omega)} + \|\nabla \phi\|^2_{L^2(\Omega)}
\big),
$
which further gives $\|\nabla \phi\|_{L^2(\Omega)}^2\le \|v\|^2_{L^2(\Omega)}$.
Therefore, \eqref{kun009b}  becomes
\begin{equation}
\langle G^{-1}B\vec{v}, B\vec{v}\rangle  = \|\nabla \phi\|^2_{L^2(\Omega)} \le  \|v\|^2_{L^2(\Omega)}.
\label{kun009a}
\end{equation}

To study $\langle M^{-1} B\vec{v}, B\vec{v}\rangle $, we let $\vec{y}=M^{-1}B\vec{v}$, i.e., $M\vec{y}=B\vec{v}$.
Suppose  $\vec{y}=(y_1 \cdots y_m )^t$ and  the corresponding function  $y=\sum_{i=1}^m y_i q_i$. Then $M\vec{y}=B\vec{v}$ can be rewritten as
$
\int_\Omega y \eta \mathrm{d}x = -\int_\Omega (\nabla\cdot v) \eta \mathrm{d}x, \forall \eta \in Q_h.
$
Letting $\eta=y$, then 
$
\|y\|^2_{L^2(\Omega)} = -\int_\Omega (\nabla\cdot v)y \mathrm{d}x
\le \frac{1}{2}\big( 
\|\nabla\cdot v\|^2_{L^2(\Omega)} + \|y\|^2_{L^2(\Omega)}
\big).
$
This results in $\|y\|_{L^2(\Omega)} \le \|\nabla\cdot v\|_{L^2(\Omega)}$.
It is also known that for any $v\in [H_0^1(\Omega)]^3$, for example,  \cite{Temam1984} p.\,93, 
\begin{equation}
\|\nabla\cdot v\|_{L^2(\Omega)} \le \|\nabla v\|_{L^2(\Omega)}.
\label{kun204}
\end{equation}
All these operations sum to
\begin{equation}
\langle M^{-1} B\vec{v}, B\vec{v}\rangle  = \langle \vec{y}, B\vec{v}\rangle 
 =  -\int_\Omega (\nabla\cdot v) y \mathrm{d}x
 = \|y\|^2_{L^2(\Omega)} \le \|\nabla\cdot v\|_{L^2(\Omega)}^2
 \le \|\nabla v\|_{L^2(\Omega)}^2.
 \label{kun011}
\end{equation}
Inserting \eqref{kun009a} and \eqref{kun011} to \eqref{kun008}, we obtain
\begin{equation}
\lambda_{\max}
\le \frac{\alpha \|v\|^2_{L^2(\Omega)} 
+ \rho k \|\nabla v\|^2_{L^2(\Omega)} }
{1.5\|v\|^2_{L^2(\Omega)} + \nu k\|\nabla v\|^2_{L^2(\Omega)}} 
\le \max\left( \frac{\alpha}{1.5}, \frac{\rho}{\nu}\right).
\label{cat004}
\end{equation}

\end{pf}

We introduce the maximum vector norm with respect to the eigenbasis  $\{\vec{\xi}_i\}_{i=1}^m$ of $K$, that is,  
\begin{align}
\|\vec{q}\|_{K,\infty} = \max_{1\le j\le m} |q_{K,j}| \text{ when } 
\vec{q} = \sum_{j=1}^m q_{K,j} \vec{\xi}_j. 
\end{align}

\begin{thm}[Iteration convergence without convection]
The pressure solution error $\vec{\bar{p}}^{\,s}$ between the scheme \eqref{wei041} and the system \eqref{wei044} satisfies 
\label{thm_iterative_conv_no_convection}
\begin{equation}
\| \vec{\bar{p}}^{\,s+1}  \|_{K,\infty}
\le  \rho(I-K)  \|\vec{\bar{p}}^{\,s} \|_{K,\infty}.
\label{kun033}
\end{equation}
Furthermore,  if $\max\left( \frac{\alpha}{1.5}, \frac{\rho}{\nu}\right)<2$, then the iterative solution (${u}^s$, ${p}^{\,s}$) of  \eqref{ip001}, \eqref{ip002}, \eqref{ip003} converges linearly to the solution (${u}$, ${p}$) of \eqref{ip011}, \eqref{ip013}.
\end{thm}
\begin{pf}\\
Denote  
$\vec{\bar{p}}^{\,s}= \sum_{j=1}^m \bar{p}^s_{K,j} \vec{\xi}_j$. Then 
\eqref{keyrelation01} becomes
$
\sum_{j=1}^m \bar{p}^{s+1}_{K,j} \vec{\xi}_j
= \sum_{j=1}^m (1-\lambda_j) \bar{p}^s_{K,j} \vec{\xi}_j.
$
Due to the linear independency of $\{\vec{\xi}_j\}_{j=1}^m$, we arrive at
$
 \bar{p}^{s+1}_{K,j} = (1-\lambda_j) \bar{p}^s_{K,j},  j=1,\cdots, m.
$
Taking the maximum magnitudes on both sides yields \eqref{kun033}. 
By Lemma\,\ref{lem_rhoI-K} and Lemma\,\ref{lem_lambdamax}, if $\max\left( \frac{\alpha}{1.5}, \frac{\rho}{\nu}\right)<2$, then $\rho(I-K)<1$ and the convergence follows.  

\end{pf}
%%%%%%%%%%%%%%%%%%%%%%
%%%%%%%%%%%%%%%%%%%%%%
\vspace{-0.8cm}
%%%%%%%%%%%%%%%%%%%%%%%%%%%%%
\subsubsection{Iteration analysis with implicit convection}
\label{sec_FIconvection}
The finite element solution $(u^s,p^s)$ of the iterative projection method \eqref{mol001}, \eqref{mol002}, \eqref{mol003} satisfies, after dropping time step  index $(n+1)$,  
\begin{align}
1.5 (u^s, w_i) + k\nu (\nabla u^s, \nabla w_i)
+ k (\text{NL}(u^{s-1}, u^s), w_i) - k(p^s, \nabla\cdot w_i) &= (F,w_i),
\quad i=1,\cdots, l,
\label{aa003} 
\end{align}
along with \eqref{ip002} and \eqref{ip003}.
Here, $(u^{-1}, p^0)$ is the initial guess.
We introduce the following the matrices,
\begin{align}
N(w)_{ij} &= (\text{NL}(w,w_j), w_i),
\quad \forall\, i,j=1,\cdots, l, 
\label{matrix_Ac}\\
A_{N}(w) &= A+kN(w),
\label{matrix_Atilde}\\
K_N(w) &= (\alpha G^{-1} + \rho k M^{-1}) B A^{-1}_N(w) B^T.
\label{kun124}
\end{align}
Then the matrix form of \eqref{aa003}, \eqref{ip002}, \eqref{ip003} is 
\begin{align}
A_N(u^{s-1}) \vec{u}^{s} + k B^T {\vec p}^{\,s} = \vec{f}, 
\quad
G\vec{\phi}^s =\frac{1}{k} B\vec{u^s}, \quad 
M\vec{p}^{\,s+1} = M\vec{p}^{\,s} + \alpha M\vec{\phi}^s + \rho B\vec{u}^s.
\label{kunaa010}
\end{align}
The finite element solution of the limit scheme \eqref{eqn01} and \eqref{eqn02} $(u,p)\in V_h\times Q_h$ satisfies
\begin{align}
1.5(u, {w_i}) + k\nu (\nabla u, \nabla{w_i})
+ k (\text{NL}(u, u), w_i)
- k(\nabla\cdot {w_i}, p) = (F,{w_i}),
&\qquad\forall\, i=1,\cdots, l, 
\label{aa001}\\
(\nabla\cdot u, q_i) = 0, &\qquad\forall\, i=1, \cdots, m.
\label{aa002}
\end{align}
With the same argument for \eqref{wei044}, the corresponding matrix form can be written as 
\begin{align}
A_N(u)  \vec{u} + k B^T \vec{p} = \vec{f}, \quad
G\vec{\phi} =\frac{1}{k} B\vec{u}, \quad 
M\vec{p} = M\vec{p} + \alpha M\vec{\phi} + \rho B\vec{u}.
\label{kunaa011}
\end{align}

The proof of iteration convergence with the implicit convection is  by writing the system as a perturbation with the perturbation parameter $k$ of the non-convective system \eqref{keyrelation01}. This is reflected in 
\eqref{vi008} along with \eqref{kun602}. The perturbation terms with the coefficient $k$ depend on the matrices $N(w)$, $A^{-1}_N(w)$, and $I-K_N(w)$, whose estimates are given in Lemma\,\ref{lem_Ac} and \ref{lem_I-Kc}. 
In this analysis, we adopt the spectral matrix norm $\vertiii{\cdot}_2$ for square matrices, which is induced by the Euclidean metric in $\mathbb{R}^l$. For more details of this norm, see, e.g., \cite[Section 5.6]{Horn2013}. 
\begin{lem}
\label{lem_Ac}
There exists a constant $C$ only dependent on $h$ and the basis of $V_h$ such that 
\begin{align}
|(\text{NL}(u,v),w)| &\le  C \|u\|_{L^2(\Omega)} 
\|v\|_{L^2(\Omega)} \|w\|_{L^2(\Omega)},
\quad \forall u,v,w\in V_h,
\label{kun300}\\
\vertiii{ N(w) }_2 &\le C  \|w\|_{L^2(\Omega)}, 
\quad \forall w\in V_h.
\label{kun301}
\end{align}
\end{lem}
\begin{pf}\\
Using \cite[Theorem 5.4.17 and Theorem 5.6.2]{Horn2013},
$\vertiii{N(w)}_2 = \max_{\|\vec{\xi}\|_2=\|\vec{\eta}\|_2=1} {\vec{\xi}}^T N(w) \vec{\eta}$, where $\vec{\xi}, \vec{\eta}\in \mathbb{R}^l$. 
Using the connection between the vector and function forms: 
$\vec{\xi}=(\xi_1 \cdots  \xi_l )^T$,
 $\xi=\sum_{i=1}^l \xi_i w_i $, 
$\vec{\eta}=(\eta_1 \cdots  \eta_l )^t$,
 $\eta=\sum_{i=1}^l \eta_i w_i$,  
we obtain 
$\vec{\xi}^T N(w) \vec{\eta} 
=((w\cdot\nabla)\eta, \xi) + \frac{1}{2} ((\nabla\cdot w)\eta,\xi)$.
Hence, 
\begin{equation}
\vertiii{N(w)}_2  =\max_{\|\vec{\xi}\|_2=\|\vec{\eta}\|_2=1}  ((w\cdot\nabla)\xi, \eta) + \frac{1}{2} ((\nabla\cdot w)\eta,\xi).
\label{kun207}
\end{equation}
Note that
\begin{align}
|((w\cdot\nabla)\eta, \xi)| 
&\le  \|w\|_{L^4(\Omega)} \|\xi\|_{L^4(\Omega)} 
\|\nabla \eta\|_{L^2(\Omega)} \label{kun201}\\
&\le 2 \|w\|^{\frac14}_{L^2(\Omega)} \|\nabla w\|^{\frac34}_{L^2(\Omega)} 
\|\xi\|^{\frac14}_{L^2(\Omega)}
\|\nabla \xi\|^{\frac34}_{L^2(\Omega)} 
 \|\nabla \eta\|_{L^2(\Omega)}, 
\label{kun202}
\end{align}
where \eqref{kun201} is obtained by H\"older's inequality and 
\eqref{kun202} by Ladyzhenskaya's inequality in three-dimensional case, 
$\|v\|_{L^4(\Omega)} \le \sqrt{2} \|v\|^{1/4}_{L^2(\Omega)}
 \|\nabla v\|^{3/4}_{L^2(\Omega)}$, $\forall v\in V$ (see, e.g., \cite[Lemma 3.5, p.\,200]{Temam1984}). 
Next, by applying the inverse inequality $\|\nabla v\|_{L^2(\Omega)}
\le Ch^{-1} \|v\|_{L^2(\Omega)}$ for any $v\in V_h$ (e.g., \cite[Theorem 4.5.11, p.\,112]{BrennerScott}), where $C$ only depends on $\Omega$, the above estimates continue as
\begin{align}
|((w\cdot \nabla)\eta, \xi)| &\le 2 \|w\|^{\frac14}_{L^2(\Omega)}
(Ch^{-1})^{\frac34} \|w\|^{\frac34}_{L^2(\Omega)} 
\|\xi\|^{\frac14}_{L^2(\Omega)} 
(Ch^{-1})^{\frac34}
\|\xi\|^{\frac34}_{L^2(\Omega)}
Ch^{-1}\|\eta\|_{L^2(\Omega)}
\nonumber 
\\
&= C \|w\|_{L^2(\Omega)} 
\|\xi\|_{L^2(\Omega)}
\|\eta\|_{L^2(\Omega)}, 
\label{kun203}
\end{align}
where the last $C$ absorbs $h^{-2.5}$ in the previous step. 
Similarly, the term $\frac12((\nabla\cdot w)\eta,\xi)$ can be shown to have the same upper bound.
Note the inequality \eqref{kun203} also  holds for any $\xi,\eta \in V_h$, which leads to \eqref{kun300}.

Note $\|\xi\|^2_{L^2(\Omega)} = (\sum_{i=1}^l \xi_i w_i,\sum_{j=1}^l \xi_j w_j) = \sum_{i,j=1}^l \xi_i \xi_j (w_i, w_j) = \vec{\xi}^T A_0 \vec{\xi}$, where $A_0$ is the mass matrix of $V_h$ with $(A_0)_{ij}=(w_i,w_j)$.  Clearly, $A_0$ is symmetric and positive definite. Since $\|\vec{\xi}\|_2=1$, $\vec{\xi}^T A_0 \vec{\xi}\le \rho(A_0)$, the spectral radius of $A_0$. Thus, $\|\xi\|_{L^2(\Omega)} \le \sqrt{\rho(A_0)}$ and similarly 
$\|\eta\|_{L^2(\Omega)} \le \sqrt{\rho(A_0)}$. Therefore, 
\begin{align}
\bigg|((w\cdot \nabla)\eta, \xi) + \frac12((\nabla\cdot w)\eta,\xi) \bigg| \le C \rho(A_0) 
\|w\|_{L^2(\Omega)}
\le C \|w\|_{L^2(\Omega)}, 
\label{kun205}
\end{align}
where $C\rho(A_0)$ is denoted as a new $C$ in the final step,  which  only depends on $h$ and the basis of $V_h$. Finally,  \eqref{kun301} is achieved from \eqref{kun207} and \eqref{kun205}.
\end{pf}

The following lemma reveals that $I-K_N(w)$ used in \eqref{vi007} is  a first order perturbation of $I-K$.
\begin{lem}
\label{lem_I-Kc}
For each $w\in V_h$, there exists a constant $k_w>0$, such that when  $0<k<\min(k_w,1)$, the matrix $A_N(w)$ is invertible, 
\begin{equation}
\vertiii{A^{-1}_N(w)}_2 \le C(1+k \|w\|_{L^2(\Omega)}),
\label{kun500}
\end{equation}
and 
\begin{equation}
I-K_N(w) = I-K - kA_{N2}(w), 
\label{kun501}
\end{equation}
where $A_{N2}(w)$ is a matrix satisfying 
\begin{equation}
\vertiii{A_{N2}(w)}_2 \le  C \|w\|_{L^2(\Omega)}.
\label{kun502}
\end{equation}
Here, the constant $C$ is independent of $w$, but depends on $h$ and bases of $V_h$ and $Q_h$.
\end{lem}

\begin{pf}\\
Note $A_N(w) = A+kN(w) = A (I+kA^{-1} N(w))$. When $w$ and the basis of $V_h$ are fixed, the matrices $A$ and $N(w)$ are also fixed. Note $A_N$ is invertible if and only if $(I+kA^{-1} N(w))$ is invertible. 
According to \cite[Corollary 5.6.16]{Horn2013}, $(I+kA^{-1} N(w))$ being invertible is equivalent to $k\vertiii{A^{-1} N(w)}_2 <1$. 

Let 
\begin{equation}
k_w = \frac{0.5}{C \vertiii{A^{-1}}_2  \|w\|_{L^2(\Omega)}}, 
\label{def_kw}
\end{equation}
where the value of $C$ is the one used in \eqref{kun301}.
If $0<k<k_w$,  then 
$k \vertiii{A^{-1} N(w)}_2  \le 
k \vertiii{A^{-1} }_2  \vertiii{N(w)}_2
\le k C \vertiii{A^{-1} }_2  \|w\|_{L^2(\Omega)} <0.5$. 
Thus, the matrix $(I+kA^{-1} N(w))$ is invertible and  $(I+kA^{-1} N(w))^{-1}=I+kA_{N1}(w)$,  where 
$
A_{N1}(w) \triangleq  \sum_{j=1}^{\infty} (-1)^j k^{j-1} (A^{-1} N(w))^j
$
and 
\begin{equation}
\vertiii{A_{N1}(w)}_2\le \frac{\vertiii{A^{-1} N(w)}_2}{1-k\vertiii{A^{-1} N(w)}_2} \le 
2 \vertiii{A^{-1}}_2 \vertiii{ N(w)}_2.
\label{kun111}
\end{equation}
Thus, when $0<k<k_w$, $A_N(w)$ is invertible and  
\begin{align}
A^{-1}_N(w) = (I+k  A_{N1}(w)) A^{-1} 
= A^{-1} + kA_{N1}(w) A^{-1}. 
\label{formula_Atilde}
\end{align}
Furthermore,
\begin{align}
\vertiii{A^{-1}_N(w)}_2  
&= \vertiii{ (I+k  A_{N1}(w)) A^{-1}}_2
\le \vertiii{A^{-1}}_2 (1 + k \vertiii{A_{N1}(w)}_2) \label{kun499}\\
&\le \vertiii{A^{-1}}_2 (1 + 2k \vertiii{A^{-1}}_2 
\vertiii{ N(w)}_2) \label{kun498}\\
&\le \vertiii{A^{-1}}_2 (1 + 2C k \vertiii{A^{-1}}_2 
\|w\|_{L^2(\Omega)})
\label{kun497}\\
&\le  C' (1+k \|w\|_{L^2(\Omega)}),
\label{kun496}
\end{align}
where \eqref{kun111} is used in \eqref{kun499} to get \eqref{kun498}, and \eqref{kun301} is applied in \eqref{kun498} to get \eqref{kun497}. The $C'$ in \eqref{kun496} is $C'=\max(\vertiii{A^{-1}}_2,2C\vertiii{A^{-1}}^2_2)$ from \eqref{kun497}. Changing $C'$ to a new $C$ results in \eqref{kun500}.

Then \eqref{kun124} and \eqref{formula_Atilde}  lead to
\begin{align}
I-K_N(w)
&=
I-(\alpha G^{-1} + \rho k M^{-1}) B A^{-1} B^t 
- k (\alpha G^{-1} + \rho k M^{-1}) BA_{N1}(w) A^{-1} B^t  \\
&\triangleq I-K - k A_{N2}(w), 
\label{kun122}
\end{align} 
where $K$ is defined in \eqref{def_K} and 
$
A_{N2}(w) \triangleq (\alpha G^{-1} + \rho k M^{-1}) BA_{N1}(w) A^{-1} B^t.
$
Then
\begin{align}
\vertiii{A_{N2}(w)}_2 &\le
( \alpha\vertiii{ G^{-1}}_2 +   \rho k \vertiii{M^{-1}}_2 )
\vertiii{B}_2  
\vertiii{B^t}_2 \vertiii{A^{-1}}_2 \vertiii{A_{N1}(w)}_2 \\
&\le ( \alpha\vertiii{ G^{-1}}_2 +   \rho  \vertiii{M^{-1}}_2 )
\vertiii{B}_2  
\vertiii{B^t}_2 
\cdot 2 \vertiii{A^{-1}}^2_2  \vertiii{N(w)}_2 
\label{kun120}\\
&=  C_1 \vertiii{N(w)}_2 \le C  \|w\|_{L^2(\Omega)},
\label{kun121}
\end{align}
where we have used $0<k<1$ and \eqref{kun111} to get  \eqref{kun120} and Lemma\,\ref{lem_Ac} to obtain \eqref{kun121}. 
The constant $C_1$ in \eqref{kun121} denotes the product of all terms in front of $\vertiii{N(w)}_2$ in \eqref{kun120}, which only depends on the bases of $V_h$ and $Q_h$. 
The identity \eqref{kun122} and estimate \eqref{kun121}  conclude this lemma.

\end{pf}

%%%%%%%%%%%%%%%%%%%%%%%%%%%%%%%%%%%%%%%%%%%%%%%%%%%%%%
\begin{thm}[Iteration Convergence of Iterative Projection Method at a single time step]
\label{thm_with_FIconvection}
If $\rho(I-K)<1$ and the time step size $k$ is sufficiently  small, then there exists a constant $0<\gamma<1$ such that for any $s\ge 0$, 
the pressure and velocity errors between \eqref{kunaa010} 
and \eqref{kunaa011} satisfy
\begin{align}
\|\vec{\bar{p}}^{\,s} \|_{K,\infty} \le  
(\|\vec{\bar{p}}^{\,0}\|_{K,\infty}
+ \|\vec{\bar{u}}^{-1}\|_2) \gamma^s,
\quad 
\|\vec{\bar{u}}^{\,s} \|_{2} \le  \sqrt{k} 
(\|\vec{\bar{p}}^{\,0}\|_{K,\infty}
+ \|\vec{\bar{u}}^{-1}\|_2) \gamma^{s+1}. 
\label{kun604}
\end{align}
Thus, the solution $(u^s, p^s)\in V_h\times Q_h$ of \eqref{aa003}, \eqref{ip002} and \eqref{ip003} converges linearly to the solution $(u,p)\in V_h\times Q_h$ of \eqref{aa001} and \eqref{aa002}.
\end{thm}
\begin{pf}\\
{\bf Step 1. Derive the relations \eqref{kun605} and \eqref{kun606}.}
Subtracting \eqref{aa001} from \eqref{aa003} gives, for each  $i=1,\cdots, l$, 
\begin{eqnarray}
1.5(\bar{u}^s,w_i)  + k\nu (\nabla\bar{u}^s, \nabla w_i) + 
k(\text{NL}(u^{s-1}, \bar{u}^s), w_i) 
+ k(\text{NL}(\bar{u}^{s-1}, u), w_i)
-k(\bar{p}^s, \nabla\cdot w_i)=0.
\end{eqnarray}
This results in, by using the notations of vectors 
$\vec{u}^s$, $\vec{\phi}^s$, $\vec{p}^{\,s}$ and errors  
$\vec{\bar{u}}^s$, $\vec{\bar{\phi}}^s$, $\vec{\bar{p}}^{\,s}$ from Section\,\ref{sec_iter_no_convection}, 
\begin{eqnarray}
(A + kN(u^{s-1}) )\vec{\bar{u}}^s
+ kB^T\vec{\bar{p}}^{\,s} = -kA_d(u) \vec{\bar{u}}^{s-1},
\label{uh001}
\end{eqnarray}
where the matrix $A_d$ is defined as 
\begin{eqnarray}
(A_d)_{ij} &\triangleq & 
(\text{NL}(w_j, u), w_i), \quad i,j=1,\cdots, l.
\end{eqnarray}
According to Lemma\,\ref{lem_Ac}, $\vertiii{N(u^{s-1})}_2  \le C \|u^{s-1}\|_{L^2(\Omega)}$.
According to Lemma\,\ref{lem_I-Kc}, for $u^{s-1}$, there exists $k_{u^{s-1}}>0$ that depends on the norm $\|u^{s-1}\|_{L^2(\Omega)}$ (see \eqref{def_kw}) such that 
when $0<k<\min(k_{u^{s-1}},1)$, 
$A_N(u^{s-1})=A+kN(u^{s-1})$ is invertible.
Thus, we  solve \eqref{uh001} to get  
\begin{eqnarray}
\vec{\bar{u}}^s= - k A^{-1}_N(u^{s-1}) 
(B^t \vec{\bar{p}}^{\,s} + A_d \vec{\bar{u}}^{s-1}).
\label{kun600}
\end{eqnarray}
According to Lemma\,\ref{lem_Ac} and \eqref{kun300},  
$\vertiii{A_d}_2 \le C \|u\|_{L^2(\Omega)}$. 
We take $C \|u\|_{L^2(\Omega)}$ as a new $C$, which yields 
\begin{equation}
\vertiii{A_d}_2 \le C. 
\label{kun601}
\end{equation}
Applying  \eqref{kun500} and \eqref{kun601} on \eqref{kun600} gives  
\begin{align}
\|\vec{\bar{u}}^s\|_2 &\le k 
\vertiii{ A^{-1}_N(u^{s-1})}_2 
(\vertiii{B^t}_2 \|\vec{\bar{p}}^{\,s}\|_2 
+ \vertiii{A_d}_2 \|\vec{\bar{u}}^{s-1}\|_2 )\nonumber\\
&\le k C (1+k\|u^{s-1}\|_{L^2(\Omega)} )
( \|\vec{\bar{p}}^{\,s}\|_{K,\infty} 
+  \|\vec{\bar{u}}^{s-1}\|_2).
\label{kun602}
\end{align}

From the difference between \eqref{kunaa010} 
and  \eqref{kunaa011}, we obtain
\begin{align}
\vec{\bar{p}}^{\,s+1} 
&= \vec{\bar{p}}^{\,s} + 
\bigg(\frac{\alpha}{k} G^{-1} +\rho M^{-1}\bigg) B\vec{\bar{u}}^{\,s} \nonumber \\
&= \vec{\bar{p}}^{\,s} -
\big(\alpha G^{-1} +\rho  k M^{-1}\big) 
B  A^{-1}_N(u^{s-1}) 
(B^t \vec{\bar{p}}^{\,s} + A_d \vec{\bar{u}}^{s-1}) 
\label{vi006}
 \\
 &= (I-K_N(u^{s-1}))  \vec{\bar{p}}^{\,s}
- (\alpha G^{-1} +\rho  k M^{-1} ) 
 A^{-1}_N(u^{s-1}) A_d \vec{\bar{u}}^{s-1}
\label{vi007}
 \\
&= ( I - K - k A_{N2}(u^{s-1})) \vec{\bar{p}}^{\,s}
 - B_1(u^{s-1})  \vec{\bar{u}}^{s-1}, 
\label{vi008}
\end{align}
where the definition $K_N(u^{s-1})$ in \eqref{kun124} is used in \eqref{vi006} to obtain \eqref{vi007}, 
the formula \eqref{kun501} is employed to go from \eqref{vi007} to \eqref{vi008}, 
and $B_1(u^{s-1})$ is defined as 
\begin{equation}
B_1(u^{s-1}) \triangleq (\alpha G^{-1} +\rho  k M^{-1} ) 
 A^{-1}_N(u^{s-1}) A_d.
\end{equation}
Applying \eqref{kun500} and \eqref{kun601} on $B_1(u^{s-1})$ yields 
\begin{align}
\vertiii{B_1(u^{s-1})}_2 \le (\alpha \vertiii{G^{-1}}_2
+ \rho k \vertiii{M^{-1}}_2) 
\vertiii{(A^{-1}_N(u^{s-1}) }_2
\vertiii{A_d}_2 
\le C (1+k\|u^{s-1}\|_{L^2(\Omega)}).
\end{align}

We express $A_{N2}(u^{s-1}) \vec{\bar{p}}^{\,s}$
and $B_1(u^{s-1})  \vec{\bar{u}}^{s-1}$ in \eqref{vi008}
with the eigenbasis $\{\vec{\xi}_j\}_{j=1}^m$ of the matrix $K$:
\begin{align}
A_{N2}(u^{s-1}) \vec{\bar{p}}^{\,s} = 
\sum_{j=1}^m \eta_{K,j} \vec{\xi}_{j} = \vec{\eta}, 
\quad 
B_1(u^{s-1})  \vec{\bar{u}}^{s-1} 
= \sum_{j=1}^m \theta_{K,j} \vec{\xi}_{j} = \vec{\theta}.
\label{vi009}
\end{align}
It is easy to tell that
\begin{align}
\|\vec{\eta}\|_{K,\infty} \le 
C \|A_{N2}(u^{s-1}) \vec{\bar{p}}^{\,s} \|_2 
\le C \vertiii{A_{N2}(u^{s-1})}_2 
\|\vec{\bar{p}}^{\,s} \|_{K,\infty} 
\le C \|u^{s-1} \|_{L^2(\Omega)} 
\|\vec{\bar{p}}^{\,s} \|_{K,\infty}, 
\end{align}
where \eqref{kun502} is used in the last step. 
Similarly, 
\begin{align}
\|\vec{\theta}\|_{K,\infty} \le 
C \|B_1(u^{s-1})  \vec{\bar{u}}^{s-1} \|_2
\le C \vertiii{B_1(u^{s-1})}_2 
\|\vec{\bar{u}}^{s-1} \|_2  
\le C  (1+k\|u^{s-1}\|_{L^2(\Omega)})
\|\vec{\bar{u}}^{s-1} \|_2.
\end{align}

From \eqref{vi008} and \eqref{vi009}, we obtain
\begin{align}
 \vec{\bar{p}}^{\,s+1} 
 = (I - K) \vec{\bar{p}}^{\,s} -  k\vec{\eta}
 - \vec{\theta}.
\end{align}
This implies
\begin{align}
\|\vec{\bar{p}}^{\,s+1}\|_{K,\infty}
&\le \rho(I-K)\|\vec{\bar{p}}^{\,s}\|_{K,\infty} 
+ k \|\vec{\eta}\|_{K,\infty} 
+  \|\vec{\theta}\|_{K,\infty} \nonumber\\
&\le 
(\rho(I-K) + k C \|u^{s-1} \|_{L^2(\Omega)} )
\|\vec{\bar{p}}^{\,s}\|_{K,\infty}
+  C  (1+k\|u^{s-1}\|_{L^2(\Omega)})
\|\vec{\bar{u}}^{s-1} \|_2 \\
&= a_{s-1} \|\vec{\bar{p}}^{\,s}\|_{K,\infty}
 + b_{s-1} \|\vec{\bar{u}}^{s-1} \|_2, 
\label{kun603}
\end{align}
where for $s=-1, 0, \cdots$, 
\begin{align}
a_{s} &= \rho(I-K) + k C \|u^{s} \|_{L^2(\Omega)}
\le \rho(I-K) + kC ( \|\vec{\bar{u}}^s \|_2 + \|u\|_{L^2(\Omega)}), \label{def_as}
\\
b_{s} &= C  (1+k\|u^{s}\|_{L^2(\Omega)})
\le C( 1+ k ( \|\vec{\bar{u}}^s \|_2 + \|u\|_{L^2(\Omega)})).
\label{def_bs}
\end{align}
Therefore, using these two quantities, \eqref{kun602} and \eqref{kun603} become, for $s=0, 1, \cdots$,
\begin{align}
\|\vec{\bar{p}}^{\,s+1}\|_{K,\infty}
&\le  a_{s-1} \|\vec{\bar{p}}^{\,s}\|_{K,\infty}
 + b_{s-1} \|\vec{\bar{u}}^{s-1} \|_2, 
 \label{kun605}\\
\|\vec{\bar{u}}^s\|_2 &\le  k b_{s-1} 
( \|\vec{\bar{p}}^{\,s}\|_{K,\infty} 
+  \|\vec{\bar{u}}^{s-1}\|_2). 
\label{kun606}
\end{align}

{\bf Step 2. prove \eqref{kun604} by induction.}
Fix $C$ as the maximum of one and all $C$'s that appear in Step 1 (this enforces $C\ge 1$), which only depends on $u$, $h$ and bases of $V_h$ and $Q_h$. 
Let 
\begin{align}
C_1&=\|\vec{\bar{p}}^{\,0}\|_{K,\infty}
+ \|\vec{\bar{u}}^{-1}\|_2, \label{def_C1} \\
C_2&=\max(\|u^{-1}\|_{L^2(\Omega)} , CC_1+\|u\|_{L^2(\Omega} ), \label{def_C2} \\
P_1&=\rho(I-K) + C\sqrt{k} + kC(1+\sqrt{k})(\sqrt{k} C_1+\|u\|_{L^2(\Omega)}), \label{def_P1} \\
P_2&=\sqrt{k}(1+\sqrt{k}) C (1+k\sqrt{k})(\sqrt{k} C_1+\|u\|_{L^2(\Omega)}), \label{def_P2}  \\
P_3&=\sqrt{k} C (1+ k\| u^{-1}\|_{L^2(\Omega)}). 
\label{def_P3}
\end{align}
If $\rho(I-K)<1$, then there exists $k_0>0$ such that when 
$0<k<k_0$,   $\max(P_1, P_2, P_3)<1$. 
We fix a value $k$ such that 
\begin{align}
0< k<\min \left(1, k_0,  \frac{0.5}{CC_2 \vertiii{A^{-1}}_2} \right).
\label{def_k}
\end{align}
The value of $\gamma$ is chosen satisfying 
\begin{equation}
\max(P_1, P_2, P_3)<\gamma<1.
\label{def_gamma}
\end{equation}

In the induction, we need to prove \eqref{kun604} and show $A_N(u^{s})$ is invertible at each step $s\ge -1$.
In the initial step, we show $A_N(u^{-1})$ is invertible, $\|\vec{\bar{p}}^{\,0}\|_{K,\infty}\le C_1$, 
$\|\vec{\bar{u}}^{\,0}\|_{2}\le \sqrt{k} C_1 \gamma$,
and $A_N(u^0)$ is invertible. 
First, \eqref{def_k} implies 
$k<\frac{0.5}{CC_2 \vertiii{A^{-1}}_2}$, which in turn suggests 
$k< \frac{0.5}{C \vertiii{A^{-1}}_2 \|u^{-1}\|_{L^2(\Omega)} }$ due to \eqref{def_C2}. According to Lemma\,\ref{lem_Ac} and \eqref{def_kw}, $A_N(u^{-1})$ is invertible. 
When $s=0$, the inequality for $\|\vec{\bar{p}}^{\,0}\|_{K,\infty}$ in \eqref{kun604} is trivially true due to the definition of $C_1$ in \eqref{def_C1}. From \eqref{kun606}, 
\begin{align}
\|\vec{\bar{u}}^{\,0}\|_{2} \le 
k b_{-1} 
( \|\vec{\bar{p}}^{\,0}\|_{K,\infty} 
+  \|\vec{\bar{u}}^{-1}\|_2) 
= k C(1+k \|u^{-1}\|_{L^2(\Omega)}) C_1
= \sqrt{k} P_3 C_1  \le \sqrt{k} C_1 \gamma. 
\end{align}
Next, we show $A_N(u^0)$ is invertible. 
Note 
\begin{align}
\|u^0\|_{L^2(\Omega)} &\le 
\|\bar{u}^0\|_{L^2(\Omega)} + \|u\|_{L^2(\Omega)}  
\le 
 C \|\vec{\bar{u}}^0\|_2 + \|u\|_{L^2(\Omega)}   \nonumber\\
 &\le C \sqrt{k} C_1 \gamma + \|u\|_{L^2(\Omega)} 
 \le CC_1 + \|u\|_{L^2(\Omega)},
 \label{formula_u0}
\end{align}
where $0<k<1$ and $0<\gamma<1$ are used in the last step.
According to \eqref{def_k}, \eqref{def_C2}, \eqref{formula_u0}, 
\begin{align}
k < \frac{0.5}{C C_2 \vertiii{A^{-1}}_2} 
\le \frac{0.5 }{C \vertiii{A^{-1}}_2 (CC_1 + \|u\|_{L^2(\Omega)})}
\le \frac{0.5 }{C \vertiii{A^{-1}}_2  \|u^0\|_{L^2(\Omega)}}
= k_{u^0}. 
\end{align}
Due to Lemma\,\ref{lem_I-Kc} and \eqref{def_kw}, $A_N(u^0)$ is invertible.

Second, assume the inequalities in \eqref{kun604} hold and $A_N(u^{s})$ is invertible for $s\le s_0$ where $s_0\ge 0$.  
Then, \eqref{kun605} implies
\begin{align}
\|\vec{\bar{p}}^{\,s_0+1}\|_{K,\infty}
&\le  a_{s_0-1} \|\vec{\bar{p}}^{\,s_0}\|_{K,\infty}
 + b_{s_0-1} \|\vec{\bar{u}}^{s_0-1} \|_2 \nonumber\\
&\le 
[ \rho(I-K) + kC (\|\vec{\bar{u}}^{s_0-1}\|_2
+  \|u\|_{L^2(\Omega)}) ]
\|\vec{\bar{p}}^{\,s_0}\|_{K,\infty} \nonumber \\
&\quad
+ C [ 1+k  (\|\vec{\bar{u}}^{s_0-1}\|_2
+  \|u\|_{L^2(\Omega)}) ]
\|\vec{\bar{u}}^{s_0-1} \|_2  \label{kun542} \\
%%%
&\le [ \rho(I-K) + kC(\sqrt{k} C_1\gamma^{s_0} + 
\|u\|_{L^2(\Omega)} ) ] C_1 \gamma^{s_0} \nonumber\\
&\quad 
+ 
C[ 1 + k( \sqrt{k} C_1 \gamma^{s_0}+ \|u\|_{L^2(\Omega)} )]
\sqrt{k}C_1\gamma^{s_0} \label{kun544}\\
&=C_1\gamma^{s_0} \bigg[ 
\rho(I-k) + \sqrt{k} C + kC (1+\sqrt{k}) (\sqrt{k}C_1\gamma^{s_0}  + \|u\|_{L^2(\Omega)} )
\bigg]  \label{kun545} \\
&< C_1\gamma^{s_0} \bigg[ 
\rho(I-k) + \sqrt{k} C + kC (1+\sqrt{k}) (\sqrt{k}C_1 + \|u\|_{L^2(\Omega)} )
\bigg] \label{kun546}\\
&= C_1 \gamma^{s_0} P_1 \le  C_1 \gamma^{s_0+1},  \label{kun547}
\end{align}
where the inequalities in \eqref{def_as} and \eqref{def_bs} are used in attaining \eqref{kun542}, 
the inductions at $s_0-1$ and $s_0$ and the  $C_1$ notation in \eqref{def_C1} is applied to obtain \eqref{kun544}, 
$0<\gamma<1$ are applied from\eqref{kun545} to \eqref{kun546}, the definition $P_1$ in \eqref{def_P1} is employed to go from \eqref{kun546} to \eqref{kun547}, and the relation \eqref{def_gamma} is utilized in \eqref{kun547}.

Furthermore, \eqref{kun606} implies 
\begin{align}
\|\vec{\bar{u}}^{s_0+1}\|_2 &\le kb_{s_0}
(\|\vec{\bar{p}}^{s_0+1}\|_{K,\infty}
+ \|\vec{\bar{u}}^{s_0}\|_2) \nonumber \\
&\le kC [1+k(\|\vec{\bar{u}}^{s_0}\|_2
+ \|u\|_{L^2(\Omega)}  )] \cdot 
(\|\vec{\bar{p}}^{s_0+1}\|_{K,\infty}
+ \|\vec{\bar{u}}^{s_0}\|_2) \label{kun531}\\
&\le kC[ 1 + k(\sqrt{k}C_1 \gamma^{s_0+1} + \|u\|_{L^2(\Omega)} )] \cdot (C_1 \gamma^{s_0+1} + \sqrt{k} C_1 \gamma^{s_0+1} ) \label{kun532} \\
&= \sqrt{k} C_1 \gamma^{s_0+1} 
\bigg[ 
\sqrt{k} (1+\sqrt{k})  C 
( 1+ k(\sqrt{k}C_1 \gamma^{s_0+1} + \|u\|_{L^2(\Omega)} ) 
\bigg] \nonumber \\
&< \sqrt{k} C_1 \gamma^{s_0+1} 
\bigg[ 
\sqrt{k} (1+\sqrt{k})  C 
( 1+ k(\sqrt{k}C_1 + \|u\|_{L^2(\Omega)} ) 
\bigg] \label{kun534} \\
&= \sqrt{k} C_1 \gamma^{s_0+1} P_2 
\le 
 \sqrt{k} C_1 \gamma^{s_0+2}, \label{kun535}
\end{align}
where the inequality \eqref{def_bs} is used in achieving \eqref{kun531}, 
the induction of $\vec{\bar{u}}^{s_0}$ and \eqref{kun547} are  applied to get \eqref{kun532}, 
$0<\gamma<1$ is used in obtaining \eqref{kun534}, the definition $P_2$ in \eqref{def_P2} is employed from \eqref{kun534} to \eqref{kun535}, and the relation \eqref{def_gamma} is utilized in \eqref{kun535}.

Finally, following the same process as we prove that $A_N(u^0)$ is invertible, we can prove $A_N(u^{s_0+1})$ is invertible.  The induction proof is completed.

\end{pf}

\vspace{-0.8cm}
%%%%%%%%%%%%%%%%%%%%%%
\section{Stability and error analysis of the limit scheme}
\label{sec_stabilityanderror}

From here, we denote the solution at time step $n$ with  iteration number $s$ in the finite element space $V_h\times Q_h$  as $(u_h^{n,s}, p_h^{n,s})$.
Thus, the finite element formulation of the iterative scheme \eqref{mol001}--\eqref{mol003} is, at the time  $t_n$,  seeking ${u}^{n+1,s}_h\in V_h$, $p^{n+1,s}_h \in Q_h$, $\phi^{n+1,s}_h\in Q_h$, $s=0, 1, \cdots$, such that for any ${ v_h}\in V_h$, $q_h\in Q_h$,
\begin{align}
1.5 ( u_h^{n+1,s}, v_h) 
&+k\nu (\nabla u_h^{n+1,s}, \nabla  v_h)
+ k (\text{NL}(u_h^{n+1,s-1},u_h^{n+1,s}), v_h)
- k (\nabla\cdot v_h, p^{n+1,s}_h)
\nonumber \\
%% equation 1, 2nd half:
&= 
(2u^n_h -0.5u^{n-1}_h, v_h)
+ k ( f^{n+1}, v_h),
\label{eqn_u}
 \\
%% equation 2:
(\nabla \phi^{n+1,s}_h,  \nabla q_h ) 
&=
 -\frac{1}{k}
  (\nabla\cdot { u}^{n+1,s}_h, q_h),
\label{eqn_phi_FEM}
\\
%% equation 3:
(p^{n+1,s+1}_h, q_h) 
&=
( p^{n+1,s}_h + \alpha \phi_h^{n+1,s} - 
\rho (\nabla\cdot u^{n+1,s}_h),  q_h ),
\label{p_update}
\end{align}
where $u^{n+1,-1}_h = 2u^n_h -u^{n-1}_h$, 
$p_h^{n+1,0}=2p_h^n - p_h^{n-1}$. 
When the conditions in Theorem\,\ref{thm_with_FIconvection} are satisfied,   the iterative solution ($u_h^{n+1,s}$, $p_h^{n+1,s}$)  converges linearly to ($u_h^{n+1}$, $p_h^{n+1}$) of the following system
\begin{align}
1.5 ( u_h^{n+1}, v_h) 
&+
k\nu (\nabla u_h^{n+1}, \nabla  v_h)
+ k(\text{NL}(u_h^{n+1},u_h^{n+1}), v_h)
- k (\nabla\cdot v_h, p^{n+1}_h) 
\nonumber \\
%% equation 1, 2nd half:
&= 
(2u^n_h -0.5u^{n-1}_h, v_h)
+ k ( f^{n+1}, v_h),
\quad \forall v_h\in V_h,
\label{eqn_u_skew_conv}\\
(\nabla \cdot u^{n+1}_h, q_h) &= 0,
\quad \forall q_h\in Q_h.
\label{eqn_p_skew_conv}
\end{align}

\begin{defn}
A scheme is  robust if there is no negative power of viscosity in the $L^2$ error bound of velocity, and not robust otherwise.
\end{defn}

The stability and error analysis of fully implicit schemes for the Navier-Stokes equations with the Galerkin finite element method has been extensively studied, notably in a series of seminal works by Heywood and Rannacher \cite{HeywoodRannacher1982, HeywoodRannacher1986, HeywoodRannacher1988, HeywoodRannacher1990}. In particular, rigorous error estimates for the fully implicit Crank-Nicolson scheme  are provided in \cite{HeywoodRannacher1990}. Note the viscosity $\nu=1$ in these works. According to \cite[Section 4.1.1]{John2021}, the analysis developed in \cite{HeywoodRannacher1990} does not lead to a robust upper bound for small viscosity values.  
The fully implicit BDF2 scheme with mixed finite element method is studied with error analysis in \cite{deFrutos2009} but the viscosity is also set as $\nu=1$.
The existing error analysis results that directly account for viscosity in the upper bound can be found in \cite{John2021} (for general FEMs, discussed loosely between Theorem 4.3 and Example 4.4), \cite[Theorem 4.7]{SchroederLube2017} (for divergence-free FEMs), and \cite[Theorem 3.1]{OLSHANSKII2020113369} (for general FEMs). However, none of the error estimates
in \cite{John2021, SchroederLube2017, OLSHANSKII2020113369} incorporate time discretization.

On the other hand, various stability results are available in some of the papers mentioned above such as \cite[Prop.\,3.3]{HeywoodRannacher1986}, \cite[Prop.\,3.1-3.3]{HeywoodRannacher1990}, \cite[Theorem\,4.3]{John2021}, \cite[Lemma\,3.1]{SchroederLube2017}, but these stability estimates are either for a backward Euler scheme or for the semi-discrete, continuous-in-time approximation. 
As for the BDF2-type time discretizations,  one stability result is given in \cite{Akbas2017}, which analyzes a widely used linearly extrapolated BDF2 scheme with a convection term that is not fully implicit. Therefore, to the best of our knowledge, there are no available stability or error analysis results for the fully implicit BDF2 scheme \eqref{eqn_u_skew_conv} and \eqref{eqn_p_skew_conv}.
The primary challenge of the BDF2 scheme lies in establishing the recursive relationships for the quantities of interest across multiple time steps, which is solved by using the identity \eqref{Shen_identity}.
We provide the stability and error analysis for \eqref{eqn_u_skew_conv} and \eqref{eqn_p_skew_conv} below and include their proofs in the Appendix.

To simplify notations, we denote $H^m_x=H^m(\Omega)$, $m\in\mathbb{N}$, 
$L^2_x=L^2(\Omega)$,
$L^\infty_x=L^\infty(\Omega)$, $L^2_t=L^2(0,T)$, $L^\infty_t=L^\infty(0,T)$, 
$L^\infty_t L^2_x =L^\infty(0,T; L^2(\Omega))$,
$L^\infty_t H^m_x=L^\infty(0,T; H^m(\Omega))$,  and so on.

%\vspace{-0.4cm}

%%%%%%%%%%%%%%%%%%%%%%%%%%%%%%%%%%%%%%%%
%\subsection{Stability analysis of the iterative projection method}
\begin{thm}
\label{thm_stability}
Suppose either the non-div-free or div-free FEM is used. 
If $0<k<1$, then the velocity solution of the scheme
\eqref{eqn_u_skew_conv} and  \eqref{eqn_p_skew_conv}  satisfies for a given final time $T>0$, at any time step $2\le n\le \lfloor \frac{T}{k} \rfloor$, 
\begin{equation}
\|u_h^{n}\|_{L^2_x}^2 + \|2u_h^{n}-u_h^{n-1}\|_{L^2_x}^2
\le  \exp\Big(\frac{T}{1-k}\Big) \cdot  \Big( \|u_h^{1}\|_{L^2_x}^2 + \|2u_h^{1}-u_h^0\|_{L^2_x}^2  
+ 4k \sum_{j=2}^{n} \|f^{j}\|_{L^2_x}^2 \Big).
\label{cat8910}
\end{equation}
\end{thm}
%

%%%%%%%%%%%%%%%%%%%%%%%%
%\subsection{Error analysis of the iterative projection method}
%\label{sec_error_analysis}

\begin{thm}[Error estimate for iterative projection method]
\label{thm_errFI}
Suppose $u, u_t\in L^\infty_tH^{r+1}_x$ with $r\ge 2$, $u_{ttt}\in L^\infty_tL^2_x$ and $p\in L_t^\infty H^r_x$.  Let the time step size $k<1$ and the spatial mesh size $h<1$.  Then for any $2\le n \le \lfloor \frac{T}{k} \rfloor$, there exists a constant $C$ independent of the solution $(u,p)$ and the discretization parameters $k$ and $h$, such that the velocity solution of the scheme \eqref{eqn_u_skew_conv} and  \eqref{eqn_p_skew_conv}  satisfies, 
\begin{itemize}
\item[(i)] with non-div-free FEMs,
\begin{align}
&\,\|u^n - u^n_h\|^2_{L^2_x} + 
\frac{k\nu}{2} \|\nabla (u^n-u^n_h)\|^2_{L^2_x} 
\le Ch^{2r} (h^2+ \nu k) \|u\|^2_{L^\infty_t H^{r+1}_x}\nonumber\\
&\,+\exp{\Big( \frac{T M_1}{1-k  M_1}\Big)} \cdot 
\Big( \|\phi^1_h\|^2_{L^2_x} 
+ \|2\phi^1_h - \phi^0_h\|^2_{L^2_x}
+ \frac{k\nu}{2} \|\nabla \phi^1_h\|^2_{L^2_x}
+ \frac{CM^2_2}{M_1}
  \Big),
\label{error_anal_FI}
\end{align}
where $\phi^0_h$ and $\phi^1_h$ are defined in \eqref{err_defs}.
\item[(ii)] 
With div-free FEMs, the estimate  \eqref{error_anal_FI}  with $M_1$ and $M_2$ replaced with $M^{div0}_1$ and $M^{div0}_2$, respectively. 

\end{itemize}
The constants $M_1$, $M^{div0}_1$, $M_2$, and $M^{div0}_2$ are defined below, 
\begin{align}
M_1 =& M^{div0}_1 +  \frac{C}{\nu} \|u\|^2_{L^\infty_t H^{r+1}_x},
\quad  M^{div0}_1 =\, 1+ C \|u\|_{L^\infty_t H^{r+1}_x},
\label{def_M1} \\
M_2 =&\,
h^{r+1} \|u_t\|_{L^\infty_t H^{r+1}_x}
+ \frac{h^r}{\sqrt{\nu}} \|p\|_{L^\infty_t H^r_x}
+ h^r \|u\|^2_{L^\infty_t H^{r+1}_x}
+  k^2 \|u_{ttt}\|_{L^\infty_t L^2_x }
\label{def_M2}\\
M^{div0}_2 =&\,
h^{r+1} \|u_t\|_{L^\infty_t H^{r+1}_x}
+ h^{r+1} \|u\|^2_{L^\infty_t H^{r+1}_x}
+  k^2 \|u_{ttt}\|_{L^\infty_t L^2_x }.
\label{def_M2div0}
\end{align}
\end{thm}

\begin{rmk}
\label{rmk_notrobust}
In Case (i) of Theorem\,\ref{thm_errFI}, if all the norms of $u$ and $p$ are bounded, $kM_1<0.5$, 
and $\phi^1_h=\phi^0_h=0$, 
then 
when $\nu\to 0+$, 
\begin{align}
\|u^n-u^n_h\|_{L_x^2} =&\, 
O(h^{r+1}+h^r\sqrt{k\nu} )
+ \exp\Big(O\Big(\frac{1}{\nu} \Big)\Big)
\cdot O\big( h^{r+1}\sqrt{\nu} + 
 h^r + k^2\sqrt{\nu} ) \big),\\
\|\nabla(u^n-u^n_h)\|_{L_x^2} =&\,
O\left(\frac{h^{r+1}} {\sqrt{k\nu}} + h^r\right)  
+ \exp\Big(O\Big(\frac{1}{\nu} \Big)\Big)
\cdot O\bigg( \frac{h^{r}} {\sqrt{k\nu}} + \frac{h^r}{\sqrt{k}}+ k^{1.5}  \bigg).
\end{align}
The spatial accuracy order of both the $L^2$ and $H^1$ errors of velocity is $r$. On the other hand, 
when $h=0$, the  accuracy order in time step size for the velocity in $H^1$ norm is 1.5, the same as that for the rotational projection method \cite{GuermondShen2004} where the error is achieved without spatial discretizations.
\end{rmk}

\begin{rmk}
\label{rmk_robust}
In Case (ii) of Theorem\,\ref{thm_errFI}, if all the norms of $u$ and $p$ are bounded, $kM_1<0.5$, 
and $\phi^1_h=\phi^0_h=0$,  then 
when $\nu\to 0+$, 
\begin{align}
\|u^n-u^n_h\|_{L_x^2} = 
O(h^{r+1} + \sqrt{k\nu} h^r + k^2), \quad 
\|\nabla(u^n-u^n_h)\|_{L_x^2} =
\frac{1}{\sqrt{k\nu}} O(h^{r+1} + \sqrt{k\nu} h^r + k^2).
\end{align}
Since the $L^2$ error of velocity does not contain the negative power of $\nu$, this method with the div-free FEMs is robust. 
\end{rmk}

%%%%%%%%%%%%%%%%%%%%%%%%%%%%%%%
%%%%%%%%%%%%%%%%%%%%%%%%%%%%%%%
%\clearpage
\section{Numerical tests}
\label{sec_numerical_tests}
Two problems are used to test the performance of the proposed method,  where Problem 1 has the exact and smooth solution and Problem 2 is a classic lid driven cavity flow problem. The finite element space is chosen as the $P2/P1$ Taylor-Hood pair. 
The numerical integration uses a 24-point quadrature rule in \cite{KEAST1986339}, which is exact for 6th degree polynomials in a tetrahedron. 
The non-symmetric matrix equations \eqref{eqn_u} is solved by the ILUT preconditioned GMRES method \cite{Saad1986,Saad1994}.
The stopping criterion of the iteration \eqref{eqn_u} and \eqref{eqn_phi_FEM}, \eqref{p_update} is set as 
\begin{align}
( \|p^{n+1,s+1}_h - p^{n+1,s}_h\|_{\max} <\epsilon 
\text{ and } 
\|p^{n+1,s+1}_h - p^{n+1,s}_h\|_{L^2_x} <\epsilon )
\text{ or  } 
s\ge \text{IterMax}.
\label{stoppingcriterion}
\end{align}
In practice we use $\epsilon=10^{-2}$ and $\text{IterMax}=50$.

\subsection{Problem 1: exact solution is given and smooth}
In Problem 1, the spatial domain is $[0,1]^3$.
The Reynolds number is $Re=\frac{L U}{\nu}$,
where the length scale $L=1$, the velocity scale $U$ is the maximum velocity magnitude. The kinematic viscosity $\nu$ is modulated to give various values of $Re$.
The exact solution of the test problem is 
\begin{align}
u_1&=\cos(t) g(x)  g'(y) g'(z), \hspace{.38 in} u_2=\cos(t) g'(x)  g(y) g'(z),\\
u_3&=-2 \cos(t) g'(x)  g'(y) g(z), \hspace{.25 in} p=\cos(3t) \sin(2\pi x) \sin(2\pi y) z^3,
\end{align}
where $g(x)=10x^2(1-x)^2$. In this problem, $U=4.6$ and $Re=\frac{4.6}{\nu}$.
The regular tetrahedral meshes are used to discretize $\Omega$ used with $N$ subdivisions in each direction.

%%%%%%%%%%%%%%%%%%%%%%%%%%%
\subsubsection{Numerical tests of iteration convergence at a single time step}
\label{sec_test_iterations}
This group of tests examines the convergence of the projection iterations in one single time step of the scheme \eqref{eqn_u}, \eqref{eqn_phi_FEM}, and \eqref{p_update} of the test problem.   This is done by finding the velocity and pressure at $t_2$, where their values at $t_0$ and $t_1$ are chosen as the exact solutions.
To track the divergence of the numerical velocity through iterations, we adopt two measures. The first is a weak  measure, 
\begin{equation}
|\text{weak div}|_{\text{max}} =
|\pi_h(\nabla\cdot u_h^s)|_{\text{max}}
= \max_{q \text{ is a basis function of }Q_h} 
|(\nabla\cdot {u}^s_h, q)|,
\label{weak_divfreemeasure}
\end{equation}
where $\pi_h$ is the $L^2$ projection from $V_h$ to $Q_h$. 
The second measure is the strong measure,   
\begin{equation}
|\text{strong div}|_{\text{max}} 
= \max_{{\bf  x}\in \Omega} 
|\nabla\cdot {u}^s_h({x})|.
\label{eqn_strong_measure}
\end{equation}

\begin{figure}[!htbp]
\centering
\includegraphics[width=1.5in]{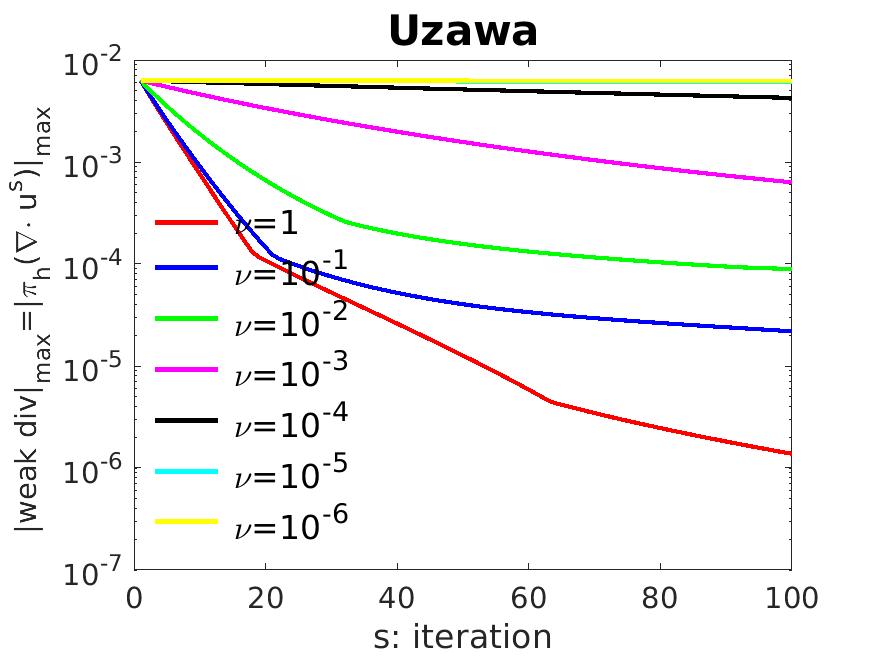}
\includegraphics[width=1.5in]{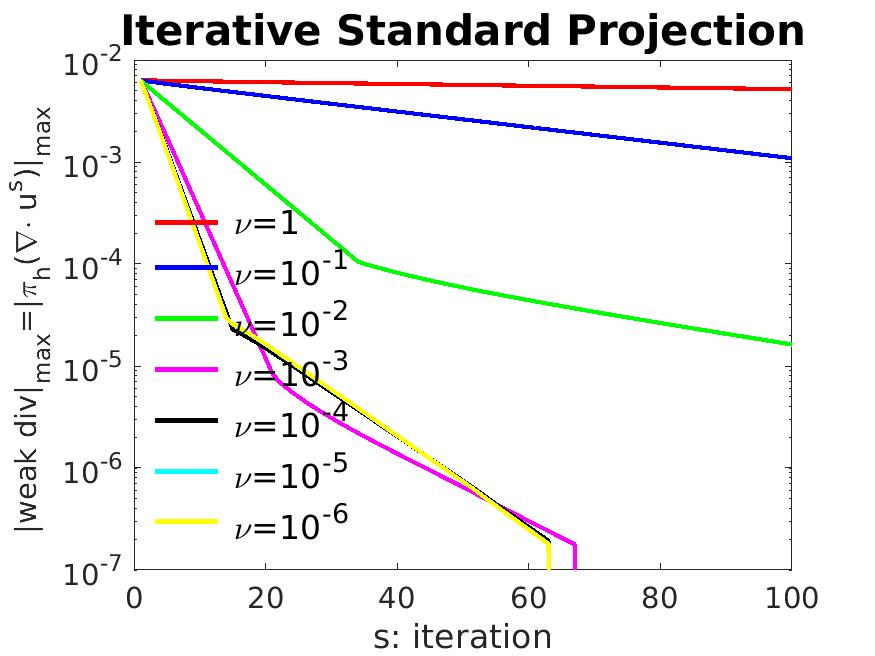}
\includegraphics[width=1.5in]{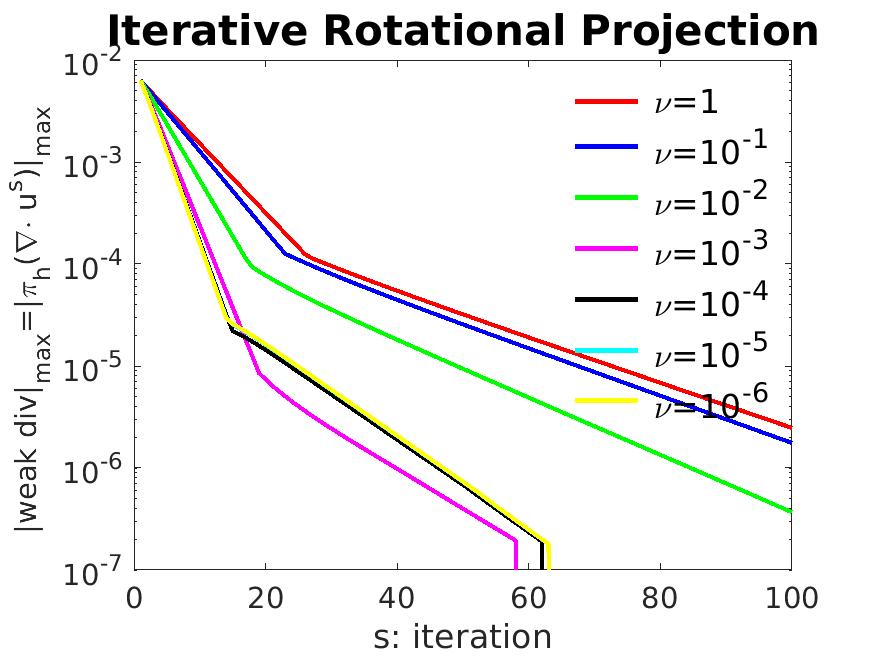}\\
\includegraphics[width=1.5in]{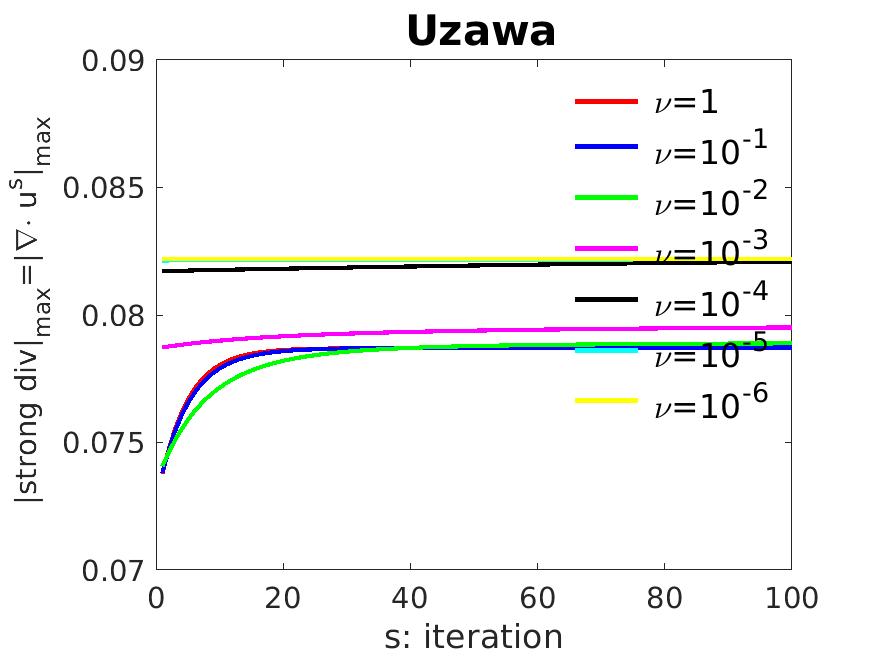}
\includegraphics[width=1.5in]{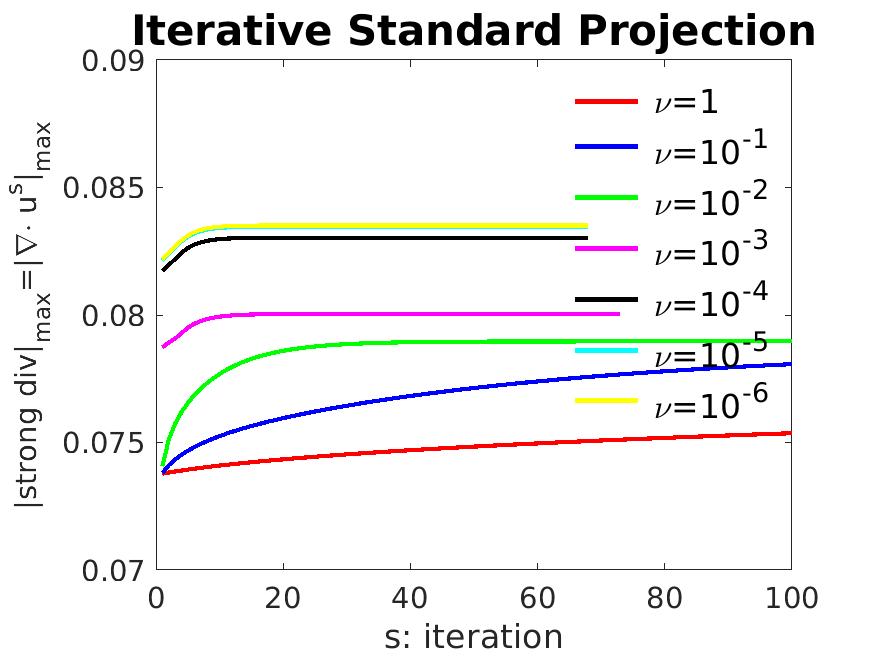}
\includegraphics[width=1.5in]{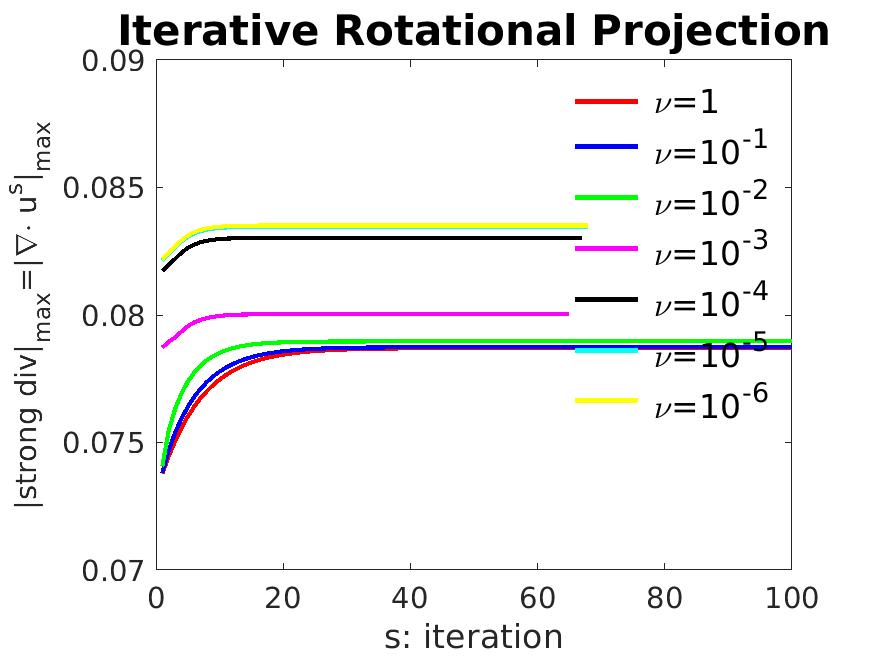}
\caption{
\scriptsize
Weak and strong measures of velocity divergence over iterations of  Problem 1. First row: weak measure. Second row: strong measure. 
%{\color{red}. The figures were generated in HPCC: DU/PaperUzawa/latex/ResultsIteration/uzawa/comp1.m.}
}
\label{Jan182019}
\end{figure}

The mesh resolution is chosen as $N=80$ and the time step size is $k=10^{-3}$. 
The computational results shown in Figure\,\ref{Jan182019} agree with the normal mode analysis. That is, the Uzawa iterations give the fast decay of the weak measure  only for large viscosity values ($\nu\ge 10^{-2}$)). The standard projection iterations converge fast only for small viscosity is small values ($\nu\le 10^{-2}$). The rotational projection iterations converge fast for all the viscosity values. The strong  measure tends to a nonzero steady state in all these simulations,  which implies the limit velocity  is not pointwise divergence free. 
This is because the Taylor-Hood FEMs are not divergence free.

Based on extensive numerical tests, we find that the optimal parameter values for achieving the fastest iteration convergence are problem dependent, which are typically found not at $\alpha=1.5$ and $\rho=\nu$, but rather within the range $\alpha\in (1.5,3)$ and $\rho\in (\nu,2\nu)$. It is recommended to test the optimal values at the initial time step with a coarse mesh before running a long time evolution. 

%\clearpage
\subsubsection{Numerical tests on accuracy of the iterative projection method}
\label{sec_test_Problem1}
This group of tests investigates the stability and accuracy of  the full scheme \eqref{eqn_u}, \eqref{eqn_phi_FEM}, and \eqref{p_update}. The simulations run from $t=0$ to $t=0.5$  with various Reynolds numbers starting from $Re=920$, where the error is evaluated at time $t=0.5$.  All the simulations in this section use $\alpha=1.5$, $\rho=\nu$, $k=0.001$.

\begin{enumerate}
\item[(1)] {\bf $Re=920$:  one iteration attains  accurate solutions}

\label{sec_Re920}
When $Re=920$ and only one iteration is used per time step, the numerical solutions show the optimal convergence (see Table\,\ref{table_pig1} (Left)): second order accuracy of  both the velocity (in $H^1$ norm) and  the pressure (in $L^2$ norm).
\begin{table}[!htbp]
\begin{center}
\scriptsize
\caption{
\scriptsize
Errors at $t=0.5$ when iteration=1. Left: $Re=920$.
Right: $Re=2300$. 
%{\color{red}HPCC: DU/Example3/YesConvection920/F1, F2, F3, F4, F5}
\label{table_pig1}}

\begin{tabular}{|c|c|c|c|c|}
\hline
 N &  $\|{ u}_h - { u}\|_{H^1}$ & 	rate & 	$|p_h-p|_{L^2}$	 &rate \\
\hline
 60	&  $0.20\times 10^{-1}$	&	1.86	 & $0.15\times 10^{-4}$	&	2.01 \\
 70	&  $0.15\times 10^{-1}$	&	2.32	 & $0.11\times 10^{-4}$	&	2.11 \\
 80	&  $0.11\times 10^{-1}$  &	1.89	 & $0.83\times 10^{-5}$	&	2.08 \\
 90	&  $0.88\times 10^{-2}$	&	2.04	 & $0.65\times 10^{-5}$	&	2.12 \\
100 &  $0.71\times 10^{-2}$	&		 & $0.52\times 10^{-5}$ &  \\
\hline
\end{tabular}
\quad 
\begin{tabular}{|c|c|c|c|c|}
\hline
 N &  $\|{ u}_h - { u}\|_{H^1}$ & 	rate & 	$|p_h-p|_{L^2}$	 &rate \\
\hline
 60 &     $0.12\times 10^{+3}$   & 11    & $0.31\times 10^{+1}$  &	 6 \\
 70 &     $0.42\times 10^{+2}$   & 33    & $0.12\times 10^{+1}$  &	35\\ 
 80 &     $0.51\times 10^{0}$   & 33    & $0.10\times 10^{-1}$ 	 &	62 \\
 90 &     $0.94\times 10^{-2}$   & 2.14  & $0.66\times 10^{-5}$  &  1.90\\
 100 &    $0.75\times 10^{-2}$   &      & $0.54\times 10^{-5}$  & \\
\hline
\end{tabular}

\end{center}
\end{table}

%%%%%%%%%%%%%%%%%%%%
\item[(2)] {\bf $Re=2300$: iterations reduce/remove spurious solutions}

\label{sec_Re2300}
The simulations with only one projection per time step leads to large errors when $N=60, 70, 80$, but  accurate solutions when $N=90, 100$ (see Table\,\ref{table_pig1} (Right)).  These errors coincide with sizes of the weak and strong incompressibility measures as shown in Figure\,\ref{fig_Re2300_iter1}: when $N=60, 70, 80$, the measures grow exponentially in time and then stay at high values, but when $N=90, 100$, the measures are small (less than 0.1). This implies that the violation of the divergence free condition causes the failure of the conventional projection method. 
\begin{figure}[!htbp]
\centering
\includegraphics[width=1.4in]
{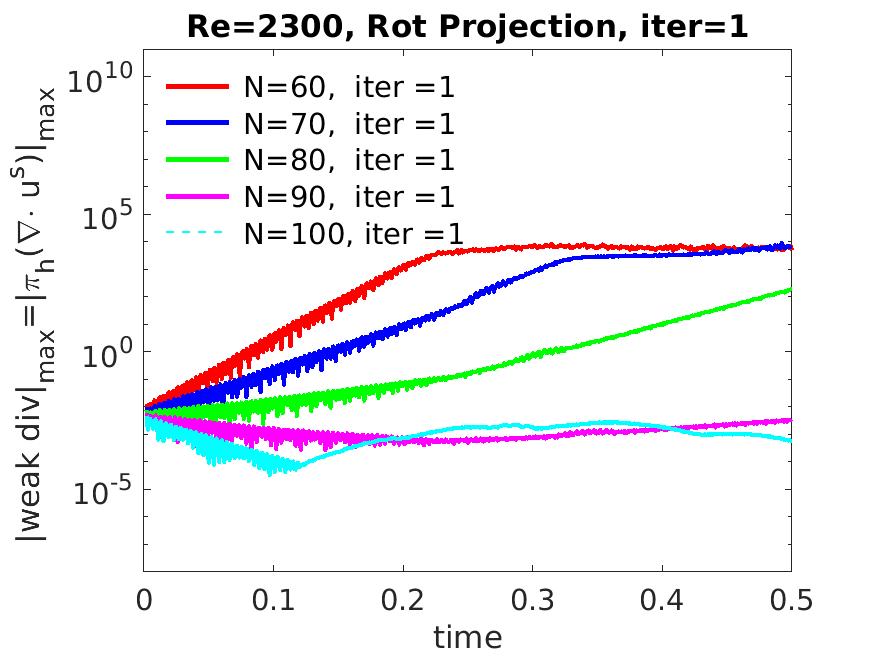}[a]
\includegraphics[width=1.4in]{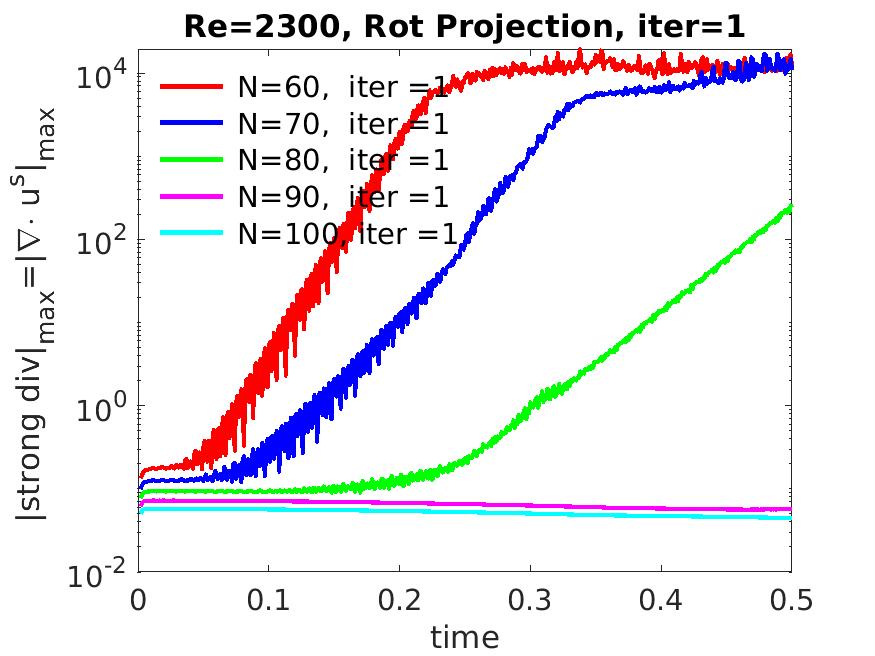}[b]
\caption{\scriptsize
Weak [a] and strong [b] incompressibility measures when $Re=2300$, iteration=1.
%{\color{red} The data are saved in HPCC: DU/Example3/YesconvectionRe2300/, [a,b]: G1, G2, G3, G4, Yw; [c,d]:Y1, Y2, Y4, Y6, Y8, Yw.  The Matlab command for [a] is `cmp_diverrmesh.m', for [b]: 'cmp_diverrmesh_errmax.m'.  }
\label{fig_Re2300_iter1}
}
\end{figure}

To see the effects of iterative projections, we  focus on the mesh $N=60$ with a fixed number of iterations per time step:  iteration=1, 2, 4, or 6.  Based on Table\,\ref{table_pig9} (Left), the errors decrease when more iterations are used and when iteration=6, the errors are optimal in the given time and spatial discretizations. 
The plottings of the numerical solutions in Figure\,\ref{Y2Y4Y6w} reveal that there are strongly unphysical oscillations  when only one iteration is used, but the oscillations reduce when more iterations are added and completely vanish when iteration=6. 
\begin{table}[!htbp]
\begin{center}
\scriptsize
\caption{
\scriptsize
Errors at $t=0.5$ when  $Re=2300$. 
Left: a fixed number of iterations per time step.
Right: stopping criterion is \eqref{stoppingcriterion} with $\epsilon=10^{-2}$
\label{table_pig9}
%{\color{red}these data are in HPCC: DU/Example3/YesconvectionRe2300/Y1, Y2, Y4, Y6, Y8, Y100, }
}
\begin{tabular}{|c|c|c|c|c|c|c|c|c|c|c|c|}
\hline
$Re$ & N & iter  & $\|{ u}_h - { u}\|_{H^1}$  & 	$|p_h-p|_{L^2}$ & CPU time	  \\
\hline
2300 & 60 & 1 & $0.12\times 10^{+3}$   &  $0.31\times 10^{+1}$  & 25 hrs \\
2300 & 60 & 2 & $0.43\times 10^{+2}$   &  $0.13\times 10^{+1}$  & 26 hrs\\ 
2300 & 60 & 4 & $0.29\times 10^{+1}$   &  $0.26\times 10^{-1}$  & 57 hrs\\
2300 & 60 & 6 & $0.23\times 10^{-1}$   &  $0.15\times 10^{-4}$  & 58 hrs\\
\hline
\end{tabular}
\quad
\begin{tabular}{|c|c|c|c|c|c|c|c|c|c|c|c|}
\hline
N & $\|{ u}_h - { u}\|_{H^1}$  & 	$|p_h-p|_{L^2}$ & average iter & CPU time\\
\hline
60  & $0.23\times 10^{-1}$ & $0.15\times 10^{-4}$ &  5.2 & 36 hrs\\
70  & $0.16\times 10^{-1}$ & $0.11\times 10^{-4}$ &  3.5 & 60 hrs\\
80  & $0.12\times 10^{-1}$ & $0.83\times 10^{-5}$ &  1.7 & 66 hrs\\
90  & $0.94\times 10^{-2}$ & $0.66\times 10^{-5}$ &
1   & 78 hrs\\
100 & $0.75\times 10^{-2}$ & $0.53\times 10^{-5}$ &  
1   & 102 hrs\\
\hline
\end{tabular}
\end{center}
\end{table} 
\begin{figure}[htbp]
\centering
\includegraphics[width=1.4in]{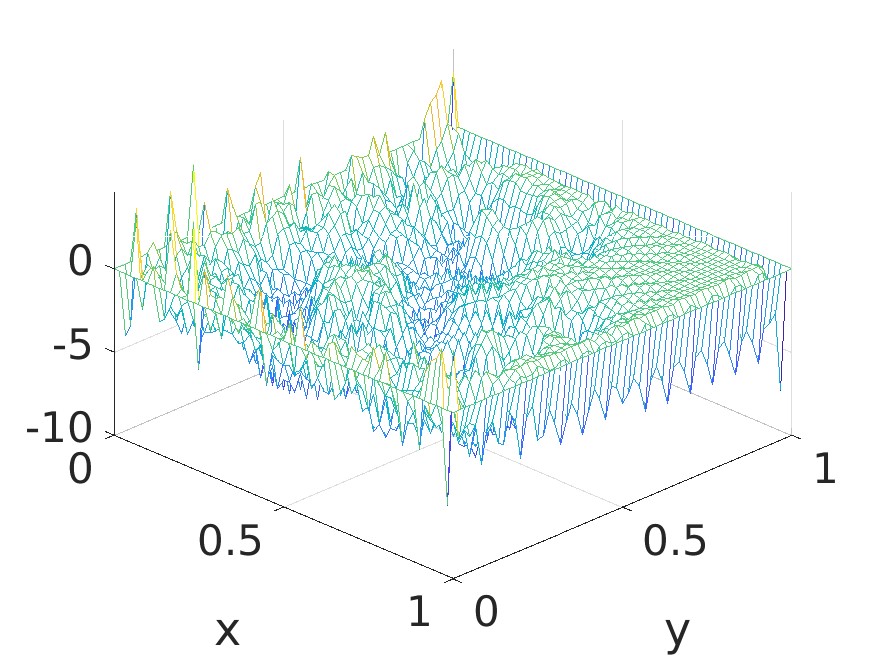}[a]
\includegraphics[width=1.4in]{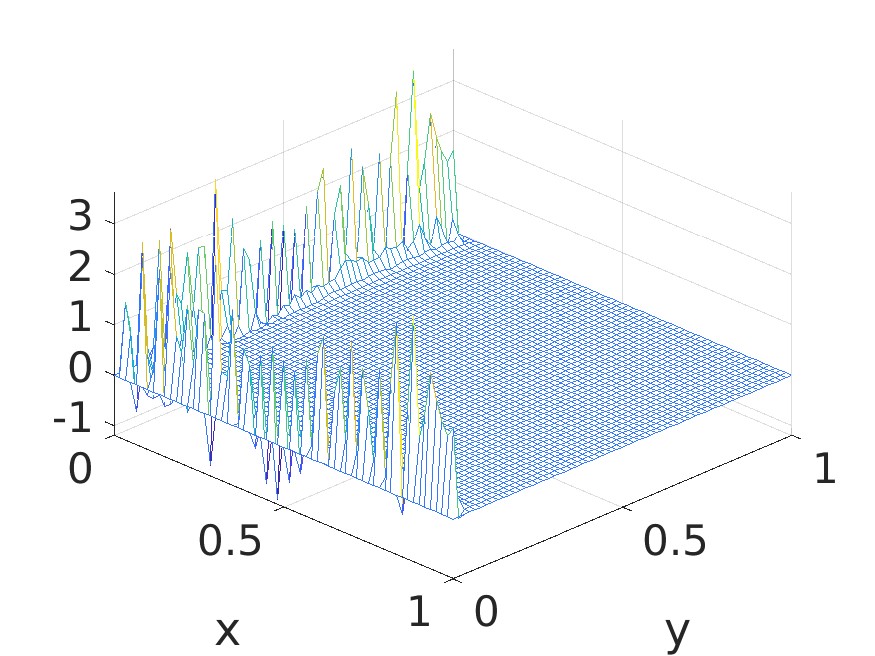}[b]
\includegraphics[width=1.4in]{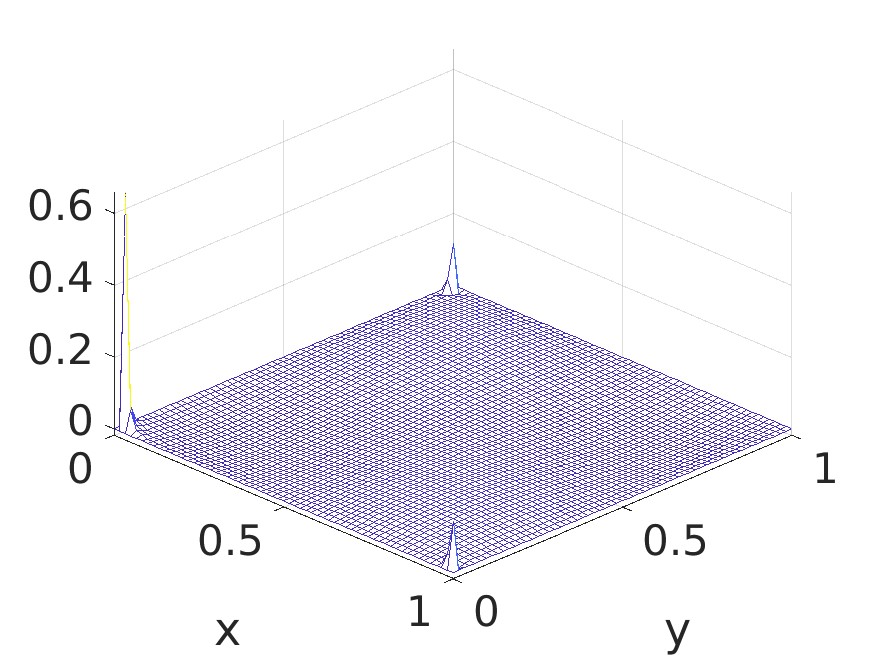}[c]
\includegraphics[width=1.4in]{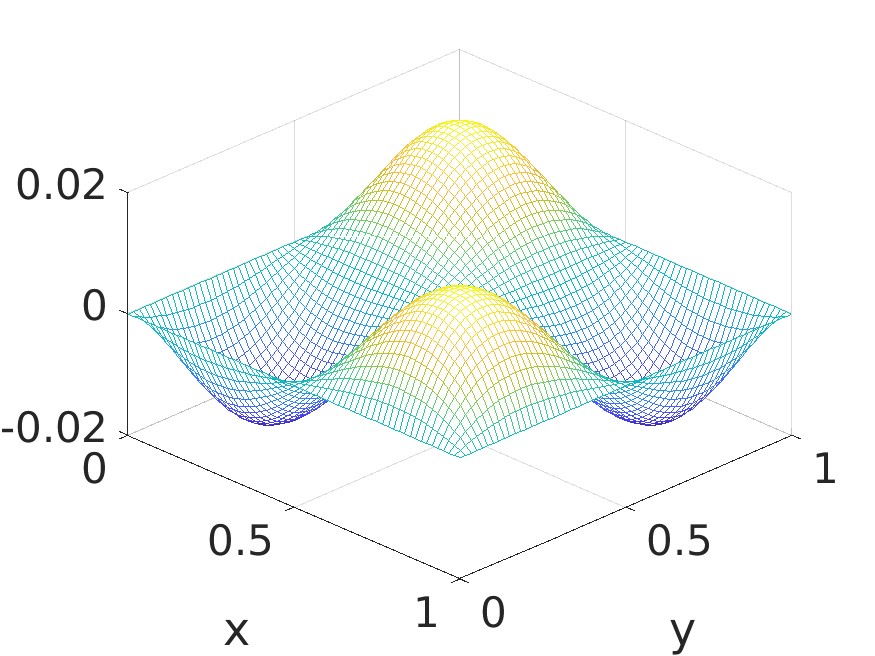}[d]
\caption{\scriptsize
Numerical solution $u_3$ on the plane $z=\frac{59}{60}$ at t=0.5 when $Re=2300$ and $N=60$. Iteration=1 in [a], 2 in [b], 4 in [c], 6 in [d].  
%{\color{red} Y2\_w.png in HPCC: DU/Example3/YesConvectionRe2300/Y1, Y2, Y4, Y6, using matlab code pp('duvp0010', 3, 60, 6), view(45, 40).}
\label{Y2Y4Y6w}
}
\end{figure}

We then use the stopping criterion \eqref{stoppingcriterion} with $\epsilon=10^{-2}$  and re-run the simulations with $Re=2300$. The optimal convergence is recovered as shown in Table\,\ref{table_pig9} [Right]. 
Note the iterative scheme spends 36 hours in the mesh with $N=60$ to obtain a solution with $\|u_h-u\|_{H^1}=0.23\times 10^{-1}$. In contrast, the convectional projection method needs to use a mesh with at least $N=90$ and 78 hours to achieve the same or better accuracy, through the comparison between Table\,\ref{table_pig1}[Right] and Table\,\ref{table_pig9}[Right].
This highlights the advantage of the iterative scheme in efficiently achieving reasonably accurate solutions while significantly reducing CPU time.

\item[(3)] {\bf Some higher Reynolds numbers}

This set of simulations applies the iterative projection method to handle some high Reynolds numbers $Re=5\times 10^3$, $1\times 10^4$, $2\times 10^4$, $4\times 10^4$,  and the stopping criterion is \eqref{stoppingcriterion} with $\epsilon=10^{-2}$. The computational results are shown in Figure\,\ref{highRe}. 
Both the $L^2$ and $H^1$ errors of velocity are optimal and increase when the Reynolds number is larger, which is consistent with the non-robust error estimates of this non-div-free FEM. The average iterations per time step also increase with the Reynolds numbers and are below 9 in these simulations.
\begin{figure}[!htbp]
\centering
\includegraphics[width=1.4in]{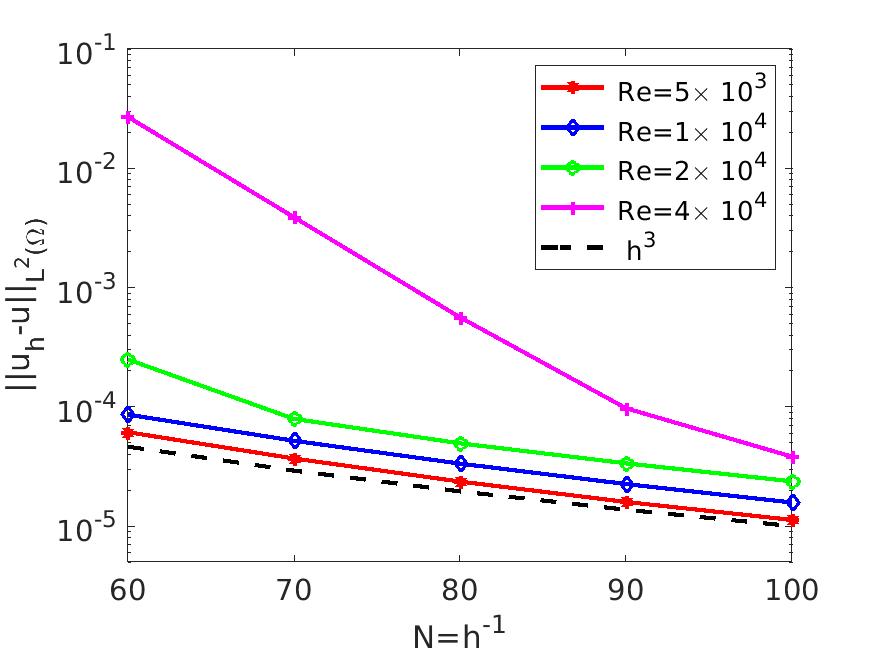}[a]
\includegraphics[width=1.4in]{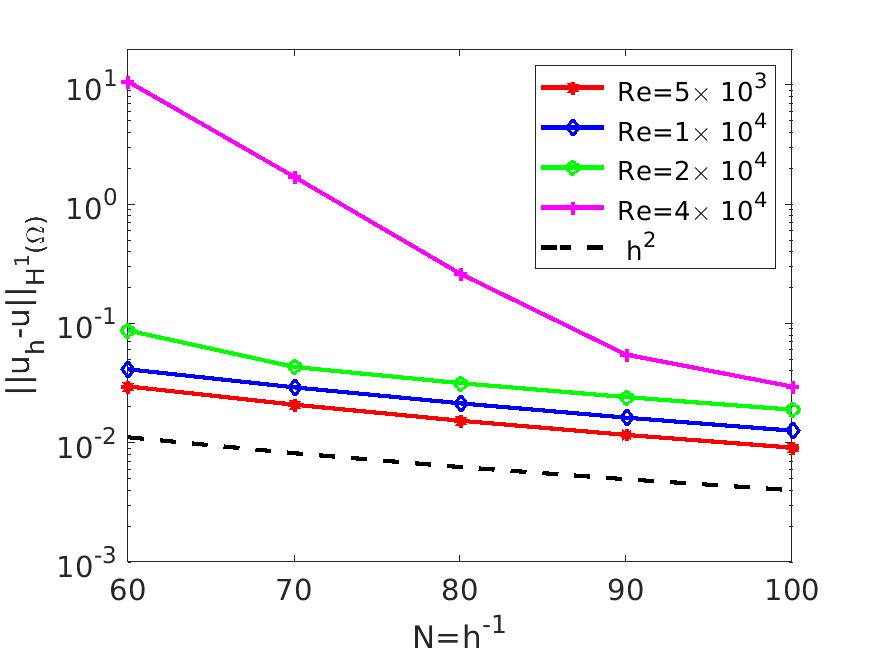}[b]
\includegraphics[width=1.4in]{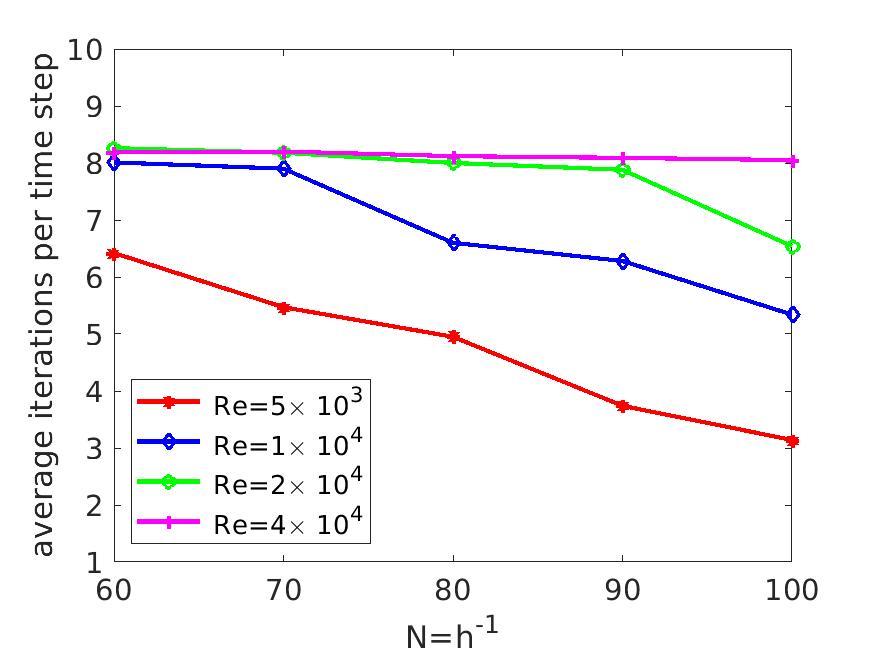}[c]
\caption{
\scriptsize
Simulation results for some high Reynolds numbers. [ab]: $L^2$ and $H^1$ errors of velocity. [c]: average iteration numbers. 
%{\color{red} these data are in HPCC: DU/Retest/ The Matlab command is `cmp(3).m', 'cmp(4_.m', cmpiter.m
\label{highRe}
}
\end{figure}

\item[(4)] {\bf Comparison with iterative projections with explicit convection}.

\label{sec_IMEX}
Following \cite{Aoussou2018},  we use a BDF2-IMEX scheme where the convection is fully explicit. The only difference with the method proposed in this work (\eqref{eqn_u}, \eqref{eqn_phi_FEM}, and \eqref{p_update})  is that the equation \eqref{eqn_u} is replaced with 
\begin{align}
1.5 ( u_h^{n+1,s}, v_h) 
&+ k\nu (\nabla u_h^{n+1,s}, \nabla  v_h)
- k (\nabla\cdot v_h, p^{n+1,s}_h) \notag\\
&=(2u^n_h -0.5u^{n-1}_h, v_h)
+ k ( f^{n+1}, v_h) -k ((w^{n+1}_h\cdot\nabla) w^{n+1}_h, v_h),
\label{eqn_u_IMEX}
\end{align}
where $w^{n+1}_h=2u^{n}_h - u^{n-1}_h$.
The numerical simulations  suggest this iterative IMEX projection scheme is not so stable as the proposed scheme 
(\eqref{eqn_u}, \eqref{eqn_phi_FEM}, and \eqref{p_update}). When $Re=920$ and $k=0.001$, this scheme blows up around $t=0.1$ (see Figure\,\ref{fig_IMEXRe920}).
In contrast, the proposed  scheme obtains the optimal convergence with just one iteration in each time step according to Section\,\ref{sec_Re920}.
\begin{figure}[!htbp]
\centering
\includegraphics[width=1.4in]{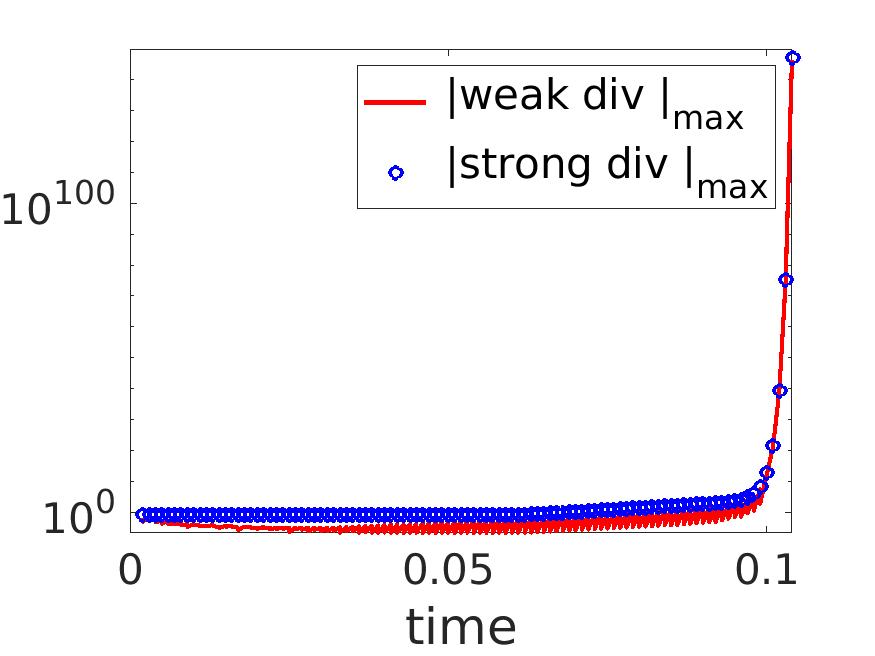}
\caption{
\scriptsize
Blowup of the iterative IMEX method with $Re=920$ and $k=0.001$. The stopping criterion is \eqref{stoppingcriterion} with $\epsilon=10^{-7}$.
%{\color{red} these data are in HPCC: DU/Example3/IMEX//IMEX60c/out/ 
%The Matlab command is `comp.m'}
\label{fig_IMEXRe920}
}
\end{figure}

\end{enumerate}

%%%%%%%%%%%%%
\subsection{Problem 2:  lid driven cavity flow}
\label{sec_liddriven}
Problem 2 is a classic lid driven cavity problem.  The domain is $\Omega=(-0.5,0.5)^3$  and the boundary velocity is ${u} = (0,1,0)$ on the sliding wall $x=-0.5$ and zero on the other walls for all time. The initial velocity is set as zero everywhere in the domain.  The meshes use Gauss-Lobatto points (e.g. \cite{ALBENSOEDER2005536}), where 
$x_i = 0.5 \cos(\frac{i\pi}{N})$, 
$y_j = 0.5\cos(\frac{j\pi}{N})$,
$z_l = \frac{l\pi}{N}-0.5$, 
$i,j,l=0,\cdots, N$,
where $N$ is the number of subdivisions in each direction.
The Reynolds number of this problem is $Re=\frac{1}{\nu}$.
In this section, the iterative projection scheme uses $\alpha=2$, $\rho=\nu$, and $k=0.001$. 

%%%%%%%%%%%%%%%%%%%%%%%%%%%
\subsubsection{Convergence of iterations at one single time step}
\label{sec_prob2_onestep}
\begin{figure}[!htbp]
\centering
\includegraphics[width=1.5in]{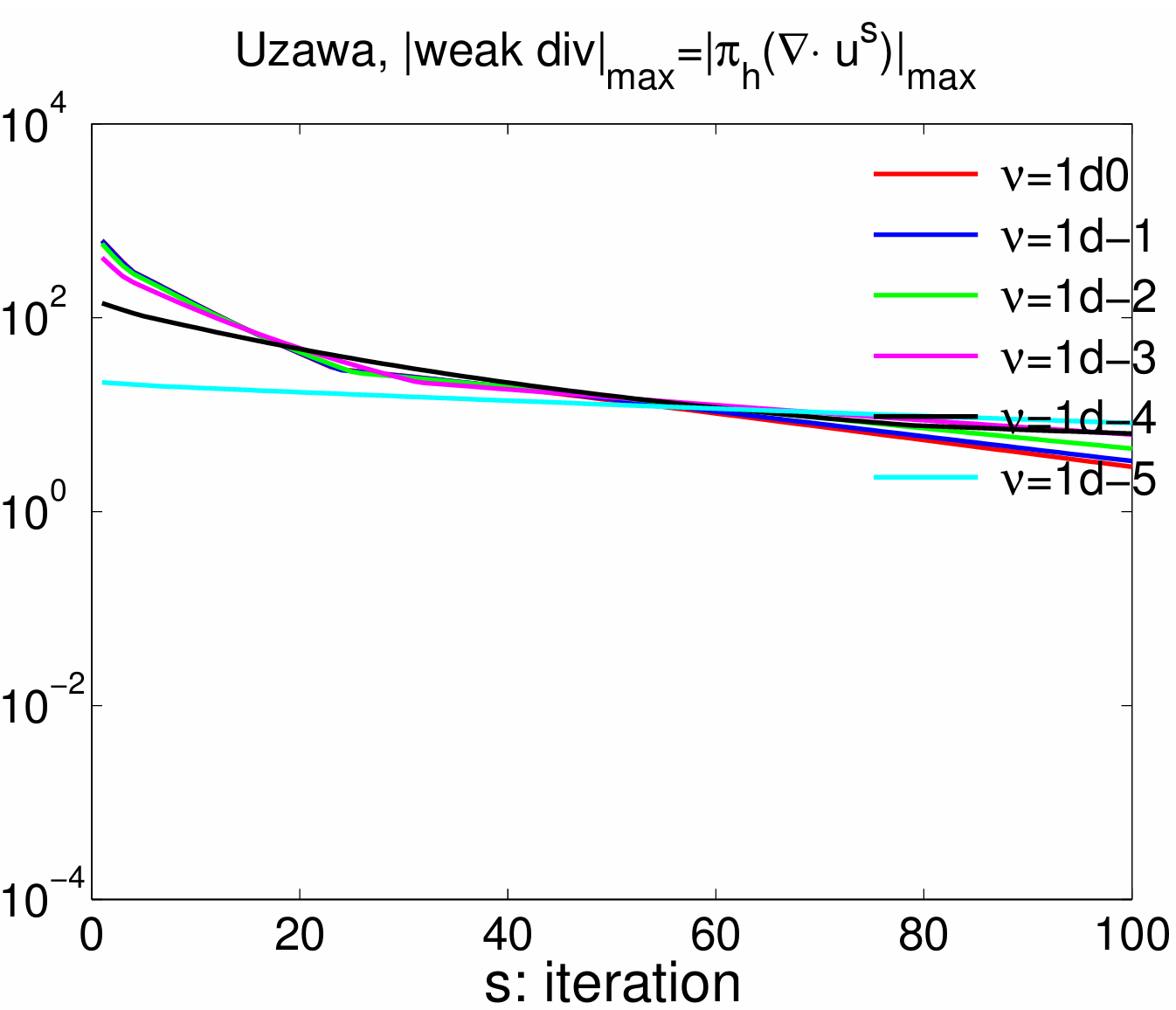}
\includegraphics[width=1.5in]{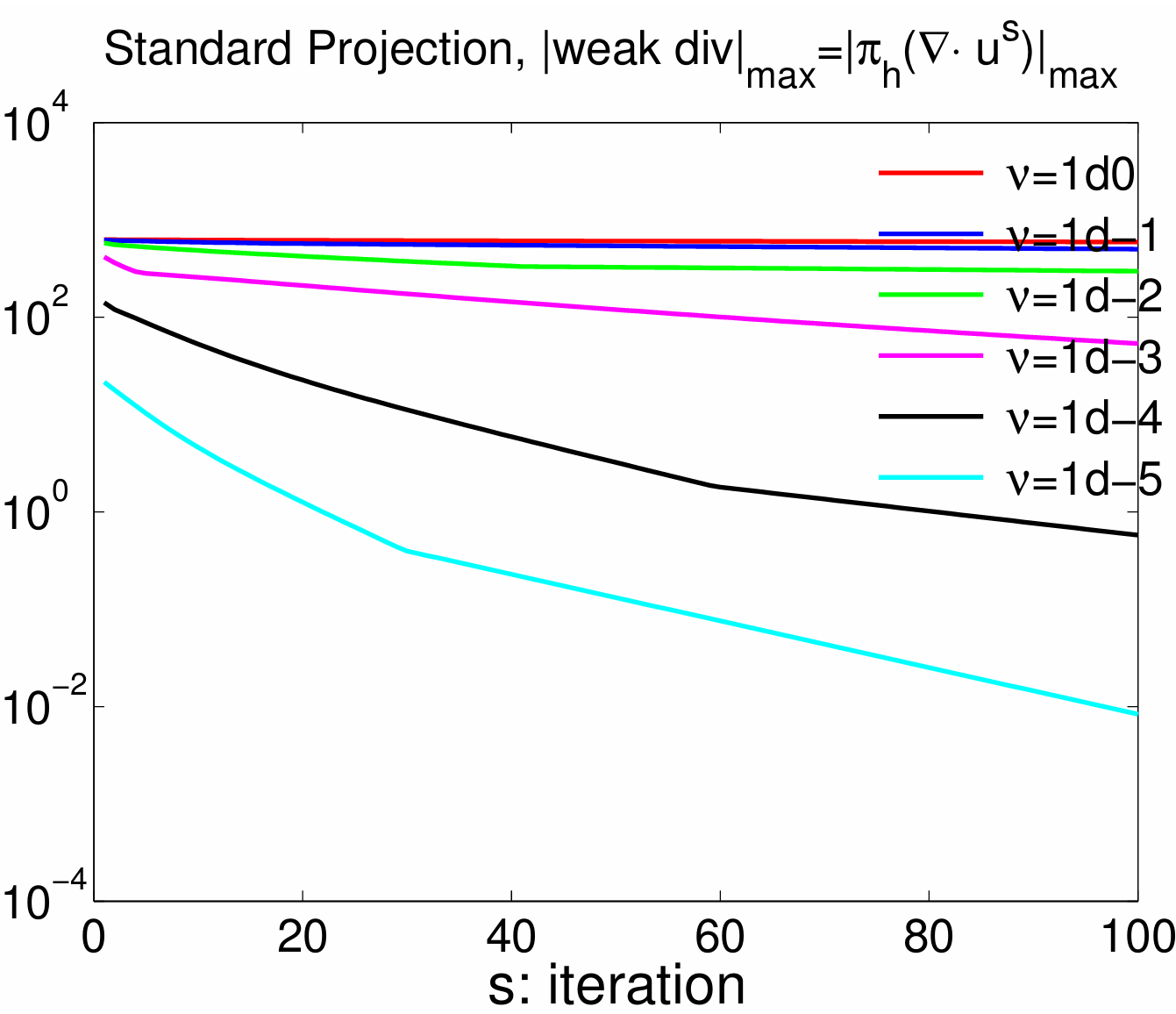}
\includegraphics[width=1.5in]{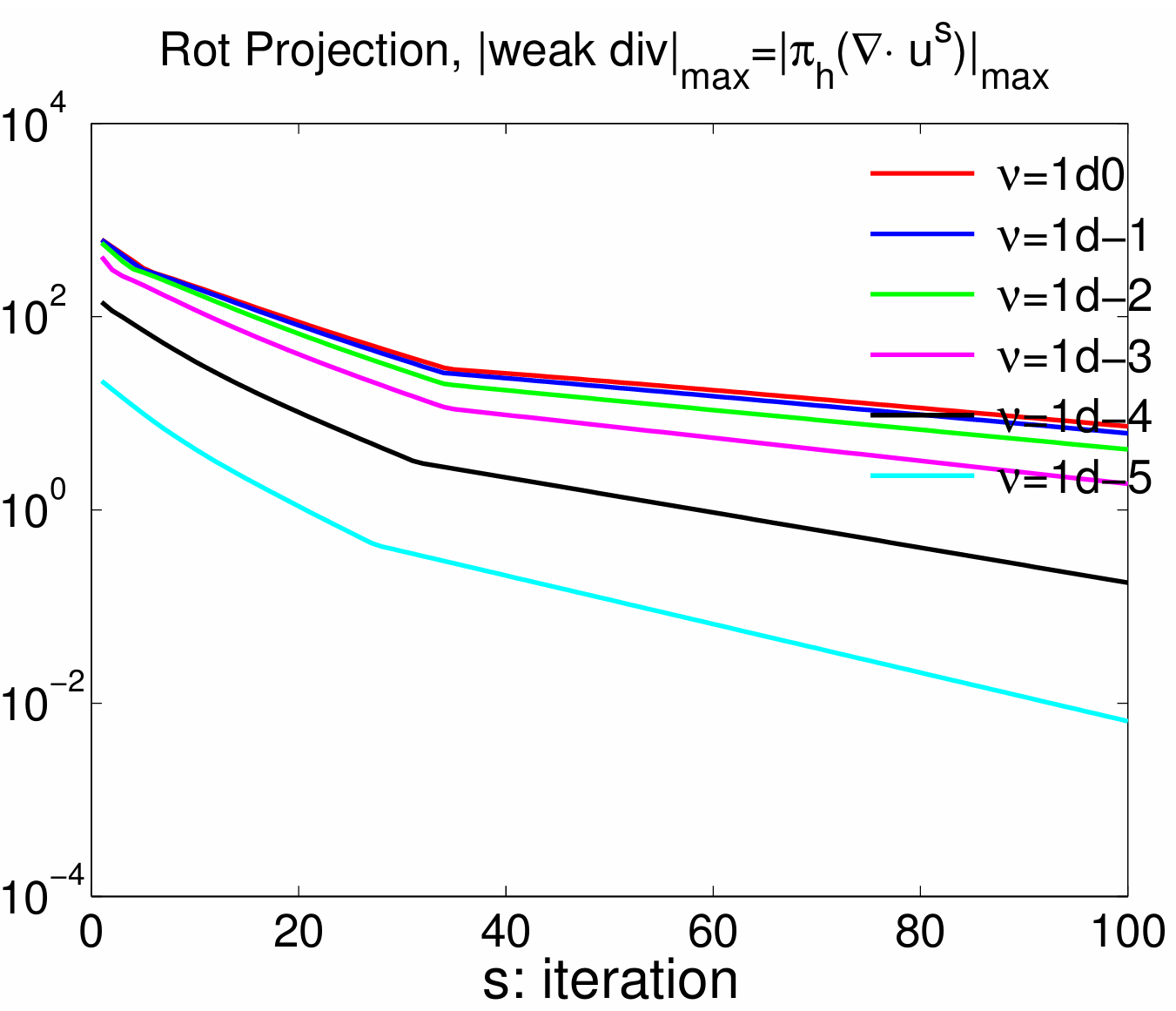}\\
\includegraphics[width=1.5in]{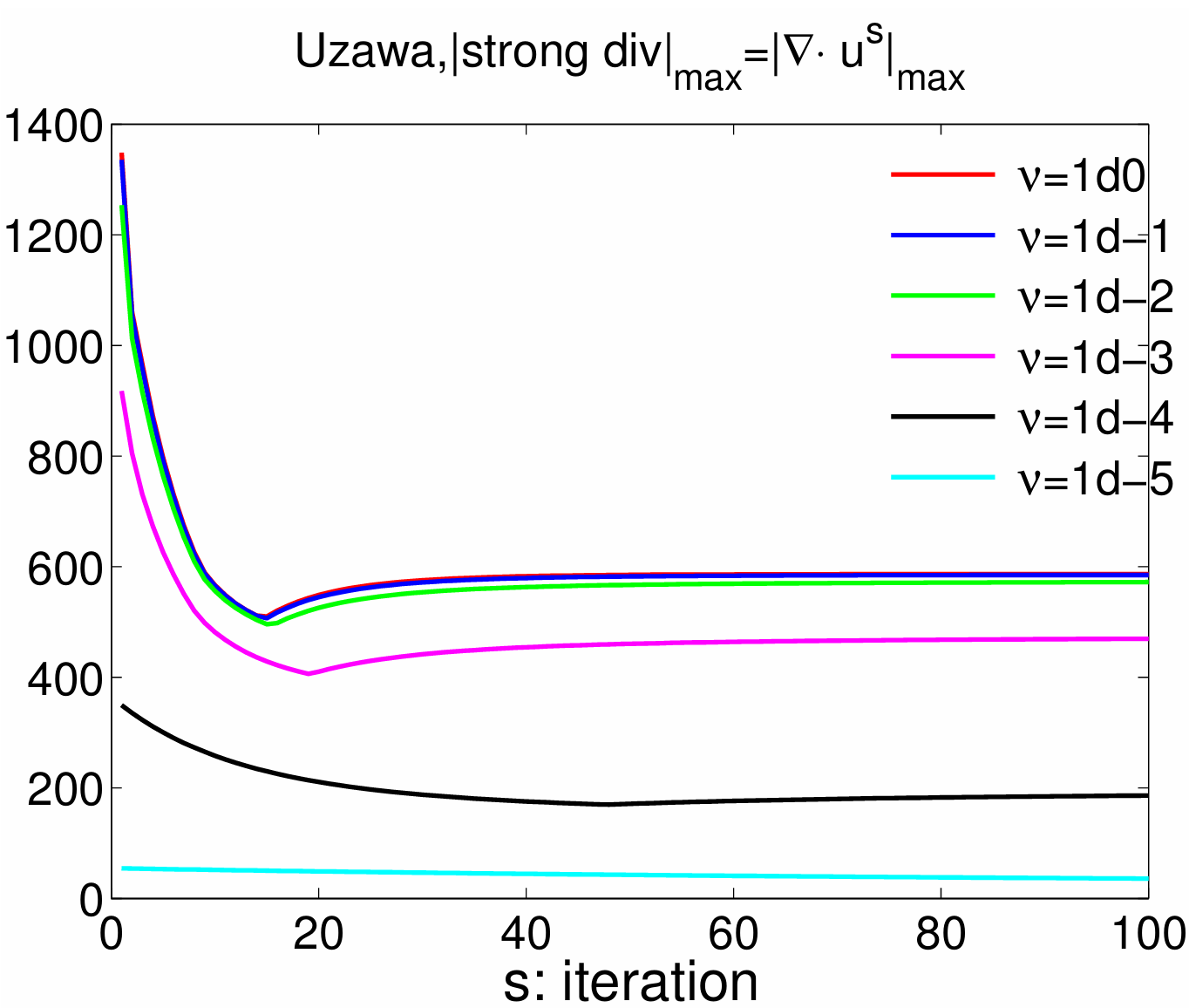}
\includegraphics[width=1.5in]{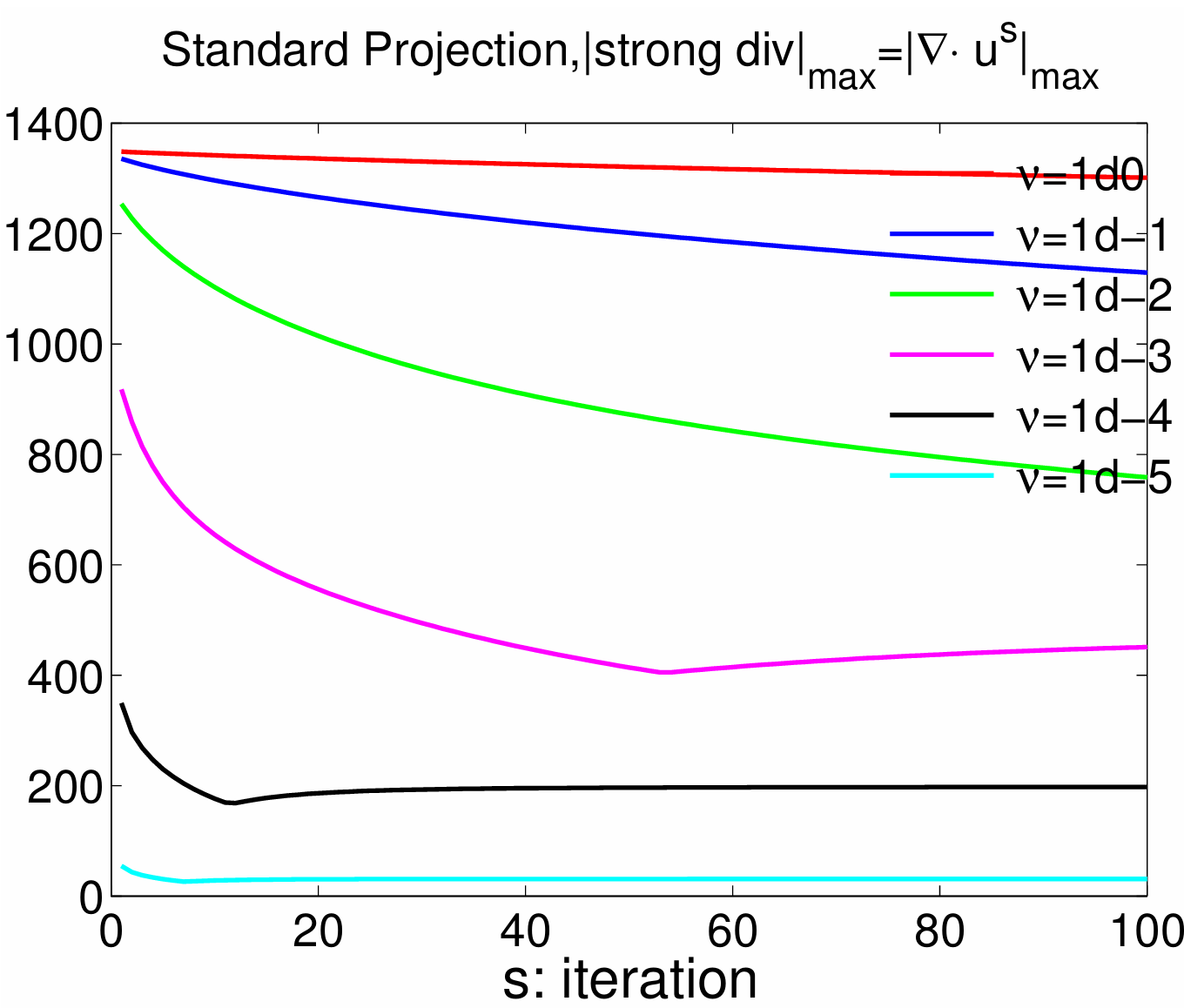}
\includegraphics[width=1.5in]{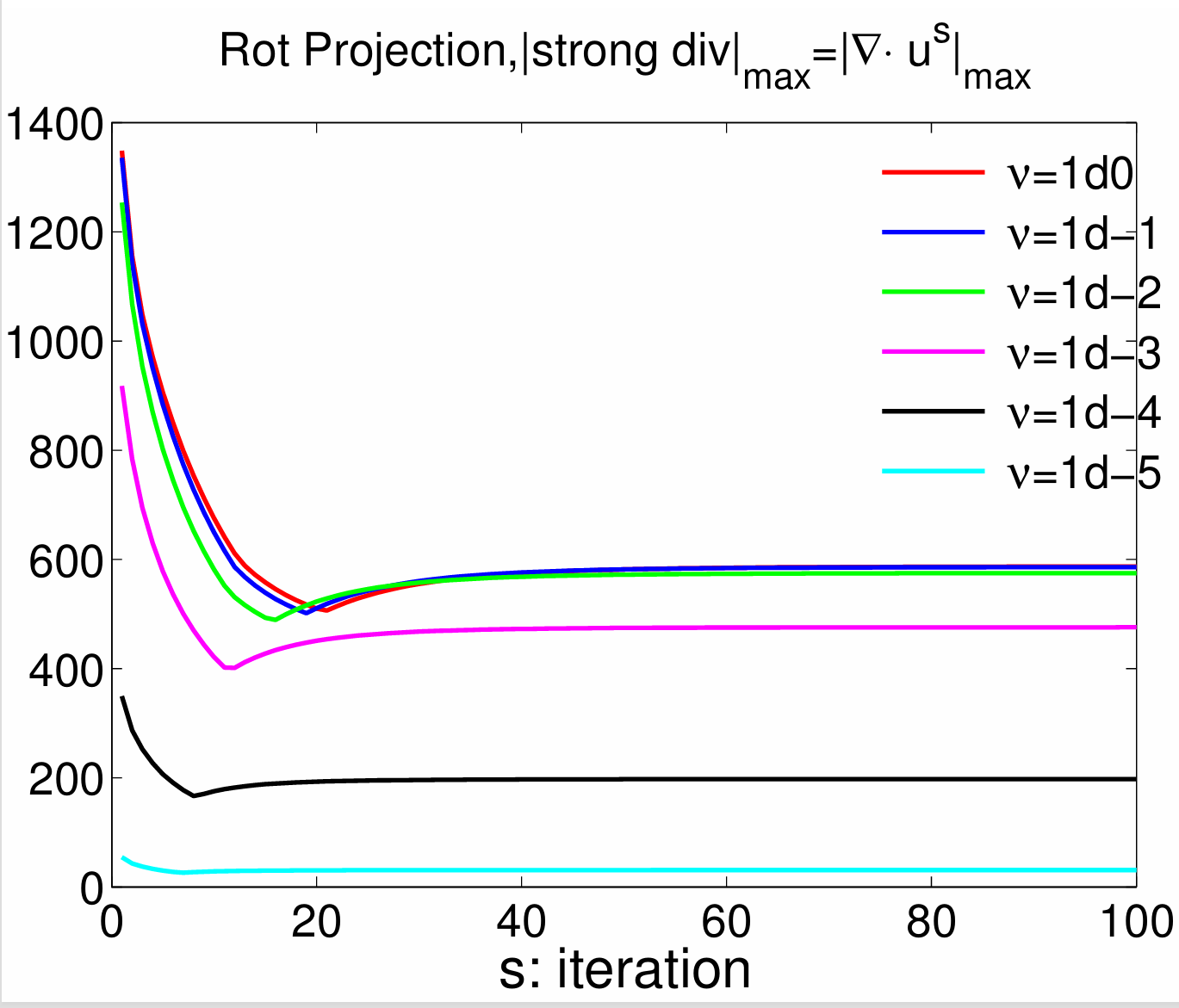}
\caption{\scriptsize
\label{April032019}
Problem 2:  Weak and strong measures of velocity divergence over iterations of problem 2. First row: weak measure. Second row: strong measure. 
%{\color{blue} These simulations have the convection term, but without Temam stabilizing term.  The data is saved in HPCC: DU/NewTest2/}
}
\end{figure}

Similar to Section,\ref{sec_test_iterations}, we first test the convergence of iterations at the time step $t_2$ for various values of the viscosity. Since the exact solution is not available, the solution at $t_1$ is taken the same as $t_0=0$. The solution at $t_2$ is obtained by running the scheme until convergence. 
We choose the number of subdivisions $N=40$ in this test. 
As shown in Figure\,\ref{April032019}, the Uzawa iterations show roughly the same behavior of weak divergence for all the viscosity values, which decrease only moderately.  The standard projection iterations perform excellently for small viscosity values but poorly for large viscosity values. The rotational projection iterations perform best. As for the strong measure of velocity divergence, both the Uzawa and rotational projection iterations exhibit roughly the same decay rates to the equilibrium states, but the standard projection iterations give very slow convergence for large viscosity values. Because of the singularity of this problem, both the weak and strong measures of the divergence are much larger than those for Problem 1 (compare Figure\,\ref{Jan182019} and Figure\,\ref{April032019}). 

%%%%%%%%%%%%
\subsubsection{Steady state solution when $Re=1000$}
\begin{figure}[!htb]
\centering
\includegraphics[height=1.5in]{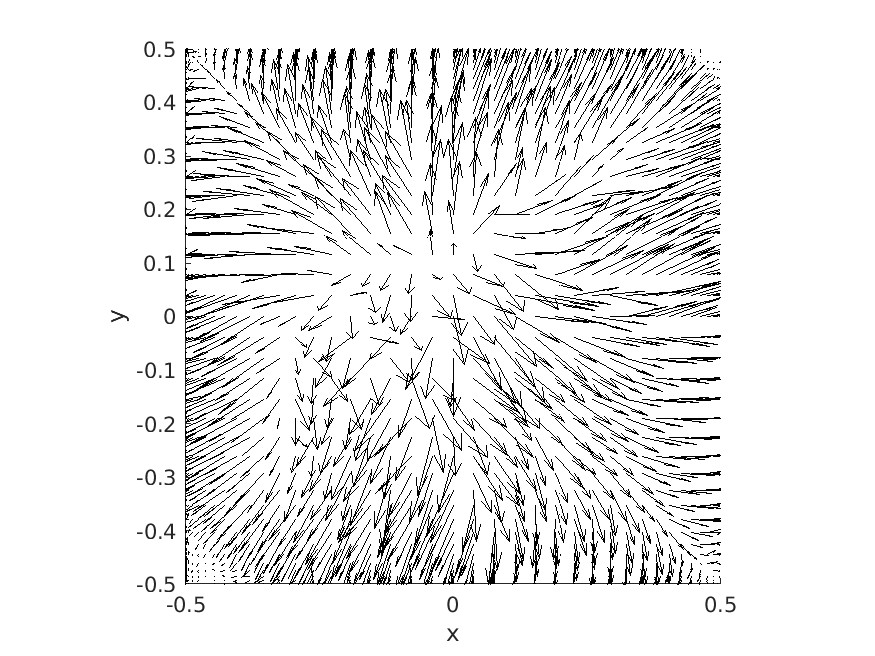}[a]
\includegraphics[height=1.5in]{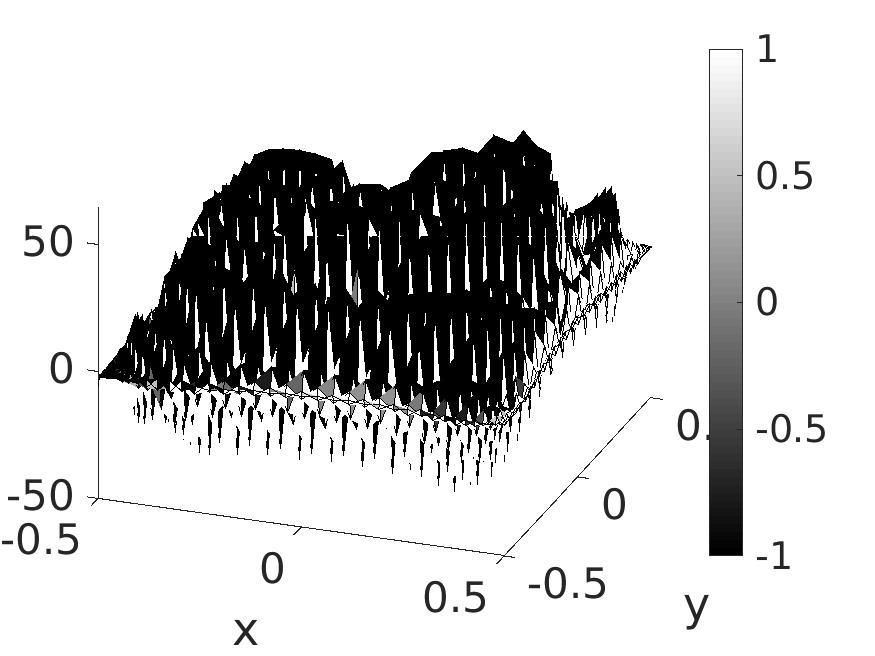}[b]
\caption{
\scriptsize
Velocity[a] and pressure[b] on the x-y plane of $z=0$ at time $t=20$. $Re=1000$, max speed=22.79. Produced by the conventional projection method.
\label{cavity_re1000_iter1}
%{\color{green}Command: ``velo4('duvp0001',2,21)'' and ``pp('duvp0001',2,21,7)''. The data are saved in folders  HPCC: DU/CubicLid_New/redi40_Re1000_1}
}
\end{figure}

\begin{figure}[!htb]
\centering
\includegraphics[height=1.5in]{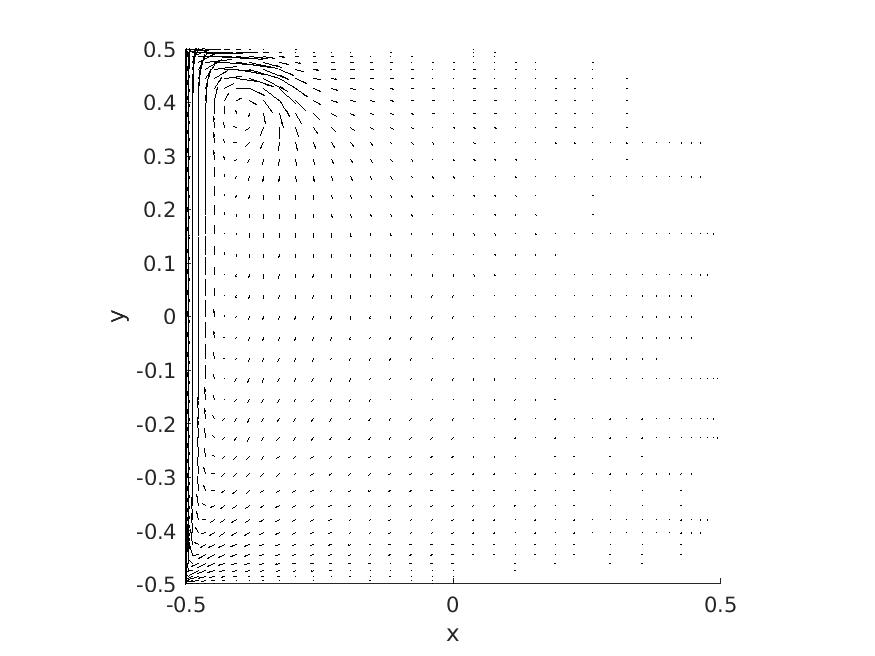}[b]
\includegraphics[height=1.5in]{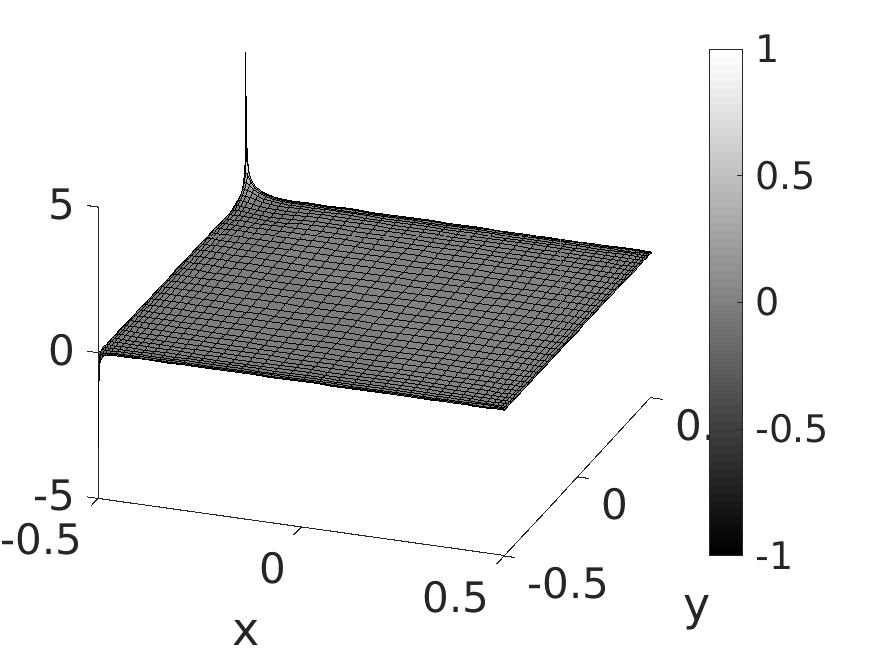}[b]
\caption{\scriptsize
Velocity[a] and pressure[b] on the x-y plane of  $z=0$  at time $t=20$. $Re=1000$, max speed=1. Produced by the iterative scheme.
\label{cavity_re1000_iter2}
%{\color{green} Command: ``velo4('duvp0001',2,21)'' and ``pp('duvp0001',2,21,7)''. Mesh size is n=40.
% The data are saved in folders HPCC: HPCC: DU/CubicLid_New/redi40_Re1000_1}
}
\end{figure}

When  $\nu=0.001$ or $Re=1000$, we run two sets of simulations, one with the conventional rotational projection method and the other with the iterative scheme by using the stopping criterion \eqref{stoppingcriterion} with $\epsilon=10^{-2}$. With this stopping criterion and $N=20$ or $40$, the average iteration number is 2 for the iterative scheme. 
When $N=40$, the conventional  projection method shows very chaotic velocity and pressure fields at $t=20$, as shown in Figure\,\ref{cavity_re1000_iter1}. But with the iterative scheme, the velocity field is very smooth (Figure\,\ref{cavity_re1000_iter2}).

\begin{figure}[!htb]
\centering
\includegraphics[height=1.5in]{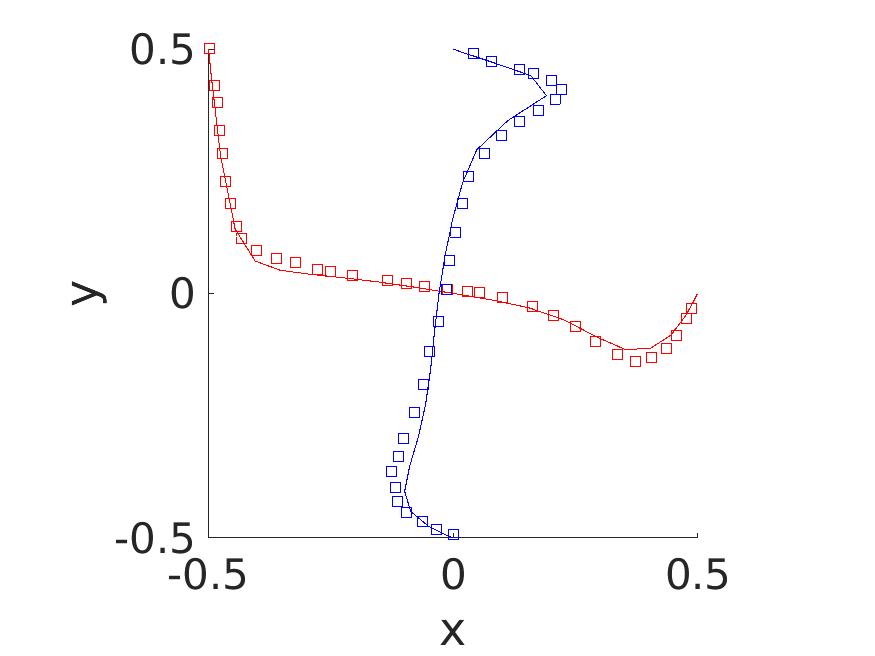}[a]
\includegraphics[height=1.5in]{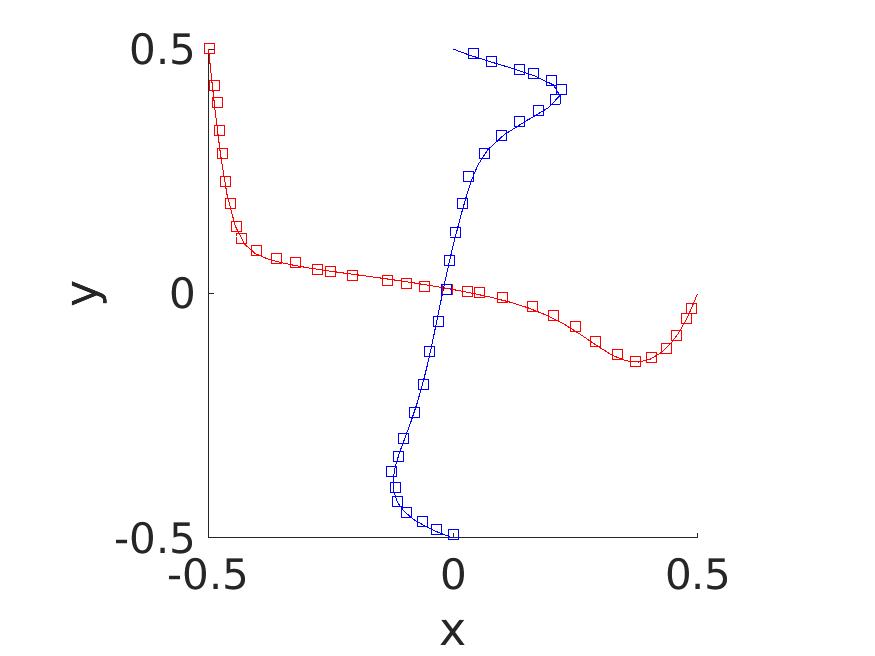}[b]
\caption{\scriptsize
Comparison of normalized velocities, $v/2$ (red) and $u/2$ (blue) on the center lines $(x, 0, 0)$ and $(0, y, 0)$ on the x-y plane $z=0$. The whole velocity vector is $(u,v,w)$. 
The square symbols are extracted from figures of 
\cite{ALBENSOEDER2005536}.  The solid lines are from the iterative projection method of this work.
[a]: $N=20$. [b]: $N=40$. 
$Re = 1000$, $t=20$. 
\label{comparewithAlbensoeder}
}
\end{figure}

The accuracy of the numerical solution is confirmed by comparing the steady state solution, achieved roughly around $t=20$,  with the result of \cite{ALBENSOEDER2005536}. The comparisons are shown in Figure\,\ref{comparewithAlbensoeder} when the mesh is refined from $N=20$ to $N=40$.  In \cite{ALBENSOEDER2005536}, the standard projection scheme with the BDF2 scheme in time and a spectral method in space is used. More importantly, they construct a special solution $u_c$ addressing the edge singularity to accompany the numerical solution. Therefore, without introducing the prescribed singular solution, the proposed iterative  projection method obtains the correct solution with averagely two projections in each time step.

\section{Concluding remarks}
\label{sec_conclusions}
This work studies a novel iterative projection method for the  Navier-Stokes equations, whose methodology is to iterate the projections to the divergence free space in every time step. To guarantee stability and convergence, a  semi-implicit skew-symmetric convection form and two tuning parameters are utilized. The benefits are significant. First, it produces the weakly divergence free velocity (pointwise divergence free velocity on div-free FEMs) when it is fully convergent. In contrast, the traditional projection method cannot produce weakly divergence free velocity.  Second, the stability and error estimates for the limit scheme, which is the fully implicit BDF2 scheme in time discretization, are theoretically established.
Third, this iterative projection method can handle high Reynolds numbers while the traditional projection method would fall short. Numerical experiments demonstrate that the desired stability and accuracy can be attained with just a few iterations per time step, without requiring full convergence.

There are many ways this work can be improved and extended. First, the scheme proposed in this work does not use any stabilization techniques mentioned in \cite{John2021}. A simple and appealing stabilization is the grad-div method as mentioned in \cite{John2018}, which has been shown to achieve robust error estimates for non-div-free FEMs. Second, there are many other forms of the convection term as shown in \cite{CHARNYI2017289}, whose theoretical and computational properties can be studied.
Third, the convergence rate of the projection iterations for \eqref{eqn01} and \eqref{eqn02} is first order, and its comparison in efficiency with Newton's method with GMRES approach of solving the nonlinear system  as in \cite{CHARNYI2017289} will be an interesting future work.
Fourth, a general BDF2 type formula proposed in \cite{HuangShen2023} may lift the temporal accuracy to second order theoretically.

\section{Appendix}

\subsection{One-step projection does not produce divergence free velocity}
\label{sec_non_divfree}
Using the matrix and vector notations in Section\,\ref{sec_iter_no_convection},
the equation \eqref{pp002} and \eqref{pp004} can be written as 
\begin{align}
G \vec{\phi}^{n+1} &= \frac{1}{k} B\vec{u}^*,\\
A_0 \vec{u}^{n+1}  &= A_0 \vec{u}^* - k B^T \vec{\phi}^{n+1}.
\end{align}
where $A_0$ is the mass matrix of $V_h$. That is, $(A_0)_{ij}=(w_i,w_j)$, $i,j=1,\cdots, l$. These two equations yield
\begin{align}
B\vec{u}^{n+1} = (I - BA_0^{-1} B^T G^{-1}) (B\vec{u}^*).
\end{align}
In general, $BA_0^{-1}B^T=G$ does not hold for mixed finite element spaces. Thus, $I-BA_0^{-1} B^T G^{-1}\ne 0$ and $B\vec{u}^{n+1}\ne 0$. This implies $u^{n+1}$ is not weakly divergence free in the one-step projection method.

\subsection{Proof of Theorem\,\ref{thm_stability}}
\label{sec-proof-stability}
\begin{pf}\\
It is easy to derive that for any $w,v\in V$, 
\begin{align}
(\text{NL}(w,v),v) = 0.
\label{kun036}
\end{align}
Thus, letting $v_h=u_h^{n+1}$ in \eqref{eqn_u_skew_conv} and 
$q_h=p^{n+1}_h$ in \eqref{eqn_p_skew_conv} eliminates the convection term  and $k(\nabla\cdot u^{n+1}_h,p^{n+1}_h)$, 
and  we obtain
\begin{equation}
(u^{n+1}_h, 1.5u^{n+1}_h-2u^n_h+0.5 u^{n-1}_h) + k\nu \|\nabla u^{n+1}_h\|^2_{L^2_x}
= k (f^{n+1}, u^{n+1}_h).
\label{nan001}
\end{equation}
Applying the following algebraic identity for real numbers $a^{n+1}, a^n, a^{n-1}$ (e.g., (4.7) in \cite{GuermondShen2004}), 
\begin{align}
a^{n+1} \cdot (1.5a^{n+1}- 2a^n+ 0.5 a^{n-1}) 
=\frac{1}{4} \big( &|a^{n+1}|^2 + |2a^{n+1}-a^n|^2 + |a^{n+1}-2a^n+a^{n-1}|^2  \nonumber \\
& - |a^n|^2 - |2a^n-a^{n-1}|^2 \big)
\label{Shen_identity}
\end{align}
on the first term in \eqref{nan001}, we obtain
\begin{align}
&\,\|u^{n+1}_h\|^2_{L^2_x} + \|2u^{n+1}_h - u^n_h\|^2_{L^2_x} + \|u^{n+1}_h -2u^n_h + u^{n-1}_h\|^2_{L^2_x}
+ 4k\nu \|\nabla u^{n+1}_h\|^2_{L^2_x} \nonumber \\
= &\,\|u^{n}_h\|^2_{L^2_x} + \|2u^{n}_h-u^{n-1}_h\|^2_{L^2_x} + 4k (f^{n+1}, u^{n+1}_h).
\end{align}
Dropping the third and fourth terms on the left and replacing the last term on the right by the upper bound in Young's inequality $|4k (f^{n+1}, u^{n+1}_h)| \le 4k \|f^{n+1}\|^2_{L^2_x} + k \|u^{n+1}\|^2_{L^2_x}$, we obtain
\begin{align}
(1-k) \|u^{n+1}_h\|^2_{L^2_x} 
+ \|2u^{n+1}_h-u^n_h\|^2_{L^2_x}
\le \|u^{n}_h\|^2_{L^2_x} 
+ \|2u^{n}_h - u^{n-1}_h\|^2_{L^2_x} 
+ 4k \|f^{n+1}\|^2_{L^2_x}.
\label{cat8921}
\end{align}
Denoting $A^{n}=\|u^{n}_h\|^2_{L^2_x} 
+ \|2u^{n}_h-u^{n-1}_h\|^2_{L^2_x}$, then \eqref{cat8921} yields
$
A^{n+1} \le \frac{1}{1-k} A^n+ \frac{4k}{1-k} \|f^{n+1}\|^2_{L^2_x}.
$
Repeating this inequality backward leads to 
$
A^{n+1} \le \frac{1}{(1-k)^n} A^1 +  4k \sum_{j=2}^{n+1} \frac{1}{(1-k)^{n+2-j}} \|f^{j}\|^2_{L^2_x}.
$
Using a basic inequality $(1+x)^j\le \exp(jx)$ for $x>-1$ and $j\in\mathbb{N}$, we get
$
\frac{1}{(1-k)^j} = \Big( 1 + \frac{k}{1-k} \Big)^j \le \exp\Big(\frac{jk}{1-k}\Big)
\le \exp\Big(\frac{nk}{1-k}\Big), j=1,\cdots, n.
$
Therefore,
$
A^{n+1} \le  \exp\Big(\frac{nk}{1-k}\Big) \cdot \Big( A^1 +  4k \sum_{j=2}^{n+1} \|f^j\|^2_{L^2_x}
\Big).
$
At any time step $1<n\le \lfloor \frac{T}{k} \rfloor$,  one has $nk\le T$ and thus the estimate \eqref{cat8910} follows. 

\end{pf}

\subsection{Proof of Theorem\,\ref{thm_errFI}}
In the error analysis, we adopt the approach from  
\cite{John2021} to treat the spatial terms. 
Let $\pi_s: V\to V^{div}_h$ be the Stokes projection with 
\begin{equation}
V^{div}_h = \{v_h\in V_h: (\nabla\cdot v_h, q_h)=0, \ \forall\, q_h\in Q_h  \}
\end{equation}
satisfying 
\begin{equation}
(\nabla \pi_s (u), \nabla v_h) = (\nabla u, \nabla v_h), \quad
\forall\, v_h\in V^{div}_h.
\end{equation}
In addition, let the $L^2$-projection $\pi_Q: Q\to Q_h$ be 
\begin{equation}
(\pi_Q(p), q) = (p,q), \quad \forall\, q\in Q_h.
\end{equation}

Suppose the degrees of local finite element polynomials in $V_h$ and $Q_h$ are $r$ and $r-1$, respectively, where $r\ge 2$. There exist the following estimates  according to \cite{John2021},
\begin{align}
\|u - \pi_s(u)\|_{H^m_x} &\le C h^{r+1-m} \|u\|_{H^{r+1}_x}, 
\quad m=0,1, 
\label{John_001}\\
\|\pi_s(u)\|_{L^\infty_x} &\le C \|u\|_{H^2_x},\label{John_002}\\
\|\nabla \pi_s(u)\|_{L^\infty_x} &\le 
  C \|\nabla u\|_{L^\infty_x}, \label{John_003}\\
\|\nabla \pi_s(u)\|_{L^p_x} &\le 
  C \|\nabla u\|_{L^p_x},
  \quad \quad\ \ \forall\, p\in[2,\infty), \label{John_004}\\
 \|p-\pi_Q(p)\|_{H^m_x} &\le C h^{r-m} \|p\|_{H^r_x}, \quad m=0,1. 
 \label{John_005}
\end{align}
Furthermore, the followings hold, where the first one is Agmon's inequality, 
\begin{equation}
\|u\|_{L^\infty_x} \le ( \|u\|_{H^1_x} 
\|u\|_{H^2_x} )^{1/2} \le \|u\|_{H^2_x}.
\label{agmon}
\end{equation}

Denote $s^n_h=\pi_s(u^n)$ and the error as 
\begin{equation}
e^n=u^n - u^n_h = (u^n-s^n_h) - (u^n_h-s^n_h) \triangleq 
\eta^n - \phi^n_h.
\label{err_defs}
\end{equation}

The following lemma is used in the proof of Theorem\,\ref{thm_errFI}.
\begin{lem}
\label{lem_eta}
When $u_t\in L^\infty_tH^{r+1}_x$, then 
$\|\eta^{n+1}-\eta^n\|_{L^2_x} \le Ck h^{r+1} \|u_t\|_{L^\infty_tH^{r+1}_x}$.
\end{lem}
\begin{pf}\\
Denote $u(t)= u(t,x)$ as  a map from $[0,T]$ to $V$. Note
\begin{align}
\eta^{n+1} - \eta^n=& (u^{n+1}-\pi_s(u^{n+1})) - (u^n-\pi_s(u^n))
=(I-\pi_s)(u(t^{n+1}) - u(t^n))
\nonumber \\
=& \int_{t^n}^{t^{n+1}} \partial_t (I-\pi_s) u(t) dt
=\int_{t^n}^{t^{n+1}} (I-\pi_s) u_t(t) dt.
\nonumber
\end{align}
Thus,
\begin{align}
\|\eta^{n+1}-\eta^n\|^2_{L^2_x} = &
\int_\Omega \bigg| \int_{t^n}^{t^{n+1}} (I-\pi_s) u_t(t) dt \bigg|^2 dx 
\le \int_\Omega k  \int_{t^n}^{t^{n+1}} |(I-\pi_s) u_t(t)|^2 dt dx \label{cat1102} \\
=&\, k  \int_{t^n}^{t^{n+1}}  \int_\Omega  |(I-\pi_s) u_t(t)|^2 dx dt 
\le\, k \int_{t^n}^{t^{n+1}}  \|(I-\pi_s) u_t(t)\|^2_{L^2_x}dt \label{cat1103}\\
\le&\, k^2 C h^{2r+2} \| u_t(t)\|^2_{L^\infty_tH^{r+1}_x}
\end{align}
where the Cauchy-Schwarz inequality is used in the last step of \eqref{cat1102}, Fubini's theorem in the first step of \eqref{cat1103} to switch the integration order, and the estimate \eqref{John_001} is applied in the last step.

\end{pf}

\begin{pf} \\
In the proof, we use $C$ to denote a generic constant which is independent of $u, p, k, h$, but may depend on the domain and the terminal time $T$. The value of the constant may vary line by line according to the context.

{\bf Step 1. Derivation of the error equation \eqref{eqn_NSu_06}}.
We rewrite \eqref{eqn_u_skew_conv} and \eqref{eqn_p_skew_conv} as 
\begin{align}
\Big( \frac{1.5  u_h^{n+1} - 2u^n_h + 0.5 u^{n-1}_h}{k}, v_h\Big) 
&+ \nu (\nabla u_h^{n+1}, \nabla  v_h)
+ (\text{NL}(u^{n+1}_h, u^{n+1}_h),v_h)
- (\nabla\cdot v_h, p^{n+1}_h) 
\nonumber \\
&= ( f^{n+1}, v_h), \qquad \forall\, v_h\in V_h,
\label{eqn_u_skew_02}\\
(\nabla \cdot u^{n+1}_h, q_h) &= 0,
\qquad \forall\, q_h\in Q_h.
\label{eqn_p_skew_02}
\end{align}
Consider the Navier-Stokes equations in the following weak form
\begin{align}
\Big( \frac{1.5  u^{n+1} - 2u^n + 0.5 u^{n-1}}{k}, v\Big) 
&+ \nu (\nabla u^{n+1}, \nabla  v)
+ (\text{NL}(u^{n+1}, u^{n+1}),v)
- (\nabla\cdot v, p^{n+1}) 
\nonumber \\
&= ( f^{n+1}, v) + (R^{n+1}(k),v), 
\qquad \forall\, v\in V,
\label{eqn_NSu}\\
(\nabla \cdot u^{n+1}, q) &= 0,
\qquad \forall\, q\in Q, 
\label{eqn_NSp}
\end{align}
where 
\begin{equation}
R^{n+1}(k) = 
\frac{1.5 u^{n+1}-2u^n+0.5 u^{n-1}}{k} -u^{n+1}_t
\label{def_Rk}
\end{equation}
is the truncation term from the time discretization. 
Subtracting \eqref{eqn_u_skew_02} from \eqref{eqn_NSu} with $v=v_h$ in \eqref{eqn_NSu}, we obtain
\begin{align}
&\, \Big( \frac{1.5  e^{n+1} - 2e^n + 0.5 e^{n-1}}{k}, v_h\Big) 
+ \nu (\nabla e^{n+1}, \nabla  v_h)
- (\nabla\cdot v_h, p^{n+1}- p^{n+1}_h) 
\nonumber \\
=&\, B(v_h)
+ ( R^{n+1}(k), v_h), 
\qquad\forall\, v_h\in V_h,
\label{eqn_NSu_02}
\end{align}
where 
\begin{align}
B(v_h)=(\text{NL}(u^{n+1}_h, u^{n+1}_h),v_h)
- (\text{NL}(u^{n+1}, u^{n+1}),v_h).
\label{def_B}
\end{align}
Using the notations in \eqref{err_defs}, the above equation becomes
\begin{align}
&\,\Big( \frac{1.5  \phi^{n+1}_h - 2\phi_h^n + 0.5 \phi_h^{n-1}}{k}, v_h\Big) 
+ \nu (\nabla \phi_h^{n+1}, \nabla  v_h)
+ (\nabla\cdot v_h, p^{n+1}- p^{n+1}_h)
+ B(v_h)
\nonumber \\
=&\,
\Big( \frac{1.5  \eta^{n+1} - 2\eta^n + 0.5 \eta^{n-1}}{k}, v_h\Big) 
+ \nu (\nabla \eta^{n+1},\nabla v_h)
- ( R^{n+1}(k), v_h), 
\qquad\forall\, v_h\in V_h.
\label{eqn_NSu_03}
\end{align}
Letting $v_h=\phi^{n+1}_h$ in \eqref{eqn_NSu_03}, 
we get
\begin{align}
&\,\frac{1}{k} (1.5\phi^{n+1}_h-2\phi^n_h+0.5\phi^{n-1}_h,\phi^{n+1}_h) 
+ \nu \|\nabla \phi^{n+1}_h\|^2_{L^2_x}
 + B(\phi^{n+1}_h) \nonumber\\
&\quad + (\nabla\cdot \phi^{n+1}_h, p^{n+1}-\pi_Q (p^{n+1}) ) 
+ (\nabla\cdot \phi^{n+1}_h, \pi_Q(p^{n+1})-p^{n+1}_h) 
\nonumber \\
=&\, \frac{1}{k} (1.5\eta^{n+1}-2\eta^n +0.5\eta^{n-1},\phi^{n+1}_h)
+ \nu (\nabla\eta^{n+1},\nabla\phi^{n+1}_h)
- (R^{n+1}(k),\phi^{n+1}_h).
\label{eqn_NSu_04}
\end{align}
Note that $(\nabla\cdot \phi^{n+1}_h, \pi_Q (p^{n+1})-p^{n+1}_h )=0$ because $\phi^{n+1}_h\in V^{div}_h$ and $\pi_Q (p^{n+1})-p^{n+1}_h\in Q_h$.
In addition, $ (\nabla\eta^{n+1},\nabla\phi^{n+1}_h)=0$
because $\eta^{n+1}=u^{n+1} - s^{n+1}_h$ 
and $(\nabla u^{n+1},\nabla\phi_h^{n+1}) = (\nabla s^{n+1}_h, \nabla\phi_h^{n+1})$ by Stokes projection with $\phi_h^{n+1}\in V^{div}_h$. Thus, \eqref{eqn_NSu_04} becomes
\begin{align}
&\,\frac{1}{k} (1.5\phi^{n+1}_h-2\phi^n_h+0.5\phi^{n-1}_h,\phi^{n+1}_h) 
+ \nu \|\nabla \phi^{n+1}_h\|^2_{L^2_x}
 \nonumber\\
=&\, 
\frac{1}{k} (1.5\eta^{n+1}-2\eta^n +0.5\eta^{n-1},\phi^{n+1}_h)
-(\nabla\cdot \phi^{n+1}_h, p^{n+1}-\pi_Q (p^{n+1}) ) 
 -B(\phi^{n+1}_h) - (R^{n+1}(k),\phi^{n+1}_h).
\label{eqn_NSu_05}
\end{align}

Applying the identity \eqref{Shen_identity} on the first term  on the left of \eqref{eqn_NSu_05} leads to 
\begin{align}
&\,\frac{1}{k} \big( 
\|\phi^{n+1}_h\|^2_{L^2_x}
 + \|2\phi^{n+1}_h-\phi^n_h\|^2_{L^2_x}
 + \|\phi^{n+1}_h - 2\phi^n_h + \phi^{n-1}_h\|^2_{L^2_x} 
 \big)
+ \nu \|\nabla \phi^{n+1}_h\|^2_{L^2_x}
 \nonumber\\
=&\, 
\frac{1}{k} 
\big(
\|\phi^{n}_h\|^2_{L^2_x}
 + \|2\phi^{n}_h-\phi^{n-1}_h\|^2_{L^2_x}
\big)
+ \frac{1}{k} (1.5\eta^{n+1}-2\eta^n +0.5\eta^{n-1},\phi^{n+1}_h)
\nonumber \\
&\, 
- (\nabla\cdot \phi^{n+1}_h, p^{n+1}-\pi_Q (p^{n+1}) ) 
- (R^{n+1}(k),\phi^{n+1}_h)
-  B(\phi^{n+1}_h).
\label{eqn_NSu_06}
\end{align}

%%%%%%%%%%%%%%%%%%%%%%%%%%%%%%%%%%%%%%%%
{\bf Step 2. Estimation of some terms  in  \eqref{eqn_NSu_06}}.
This refers to the last four terms on the right of \eqref{eqn_NSu_06}.

\begin{enumerate}
\item [(1)] 
 As for the second term on the right of \eqref{eqn_NSu_06}, we have
\begin{align}
&\,\frac{1}{k} \big|(1.5\eta^{n+1}-2\eta^n +0.5\eta^{n-1},\phi^{n+1}_h) \big| 
= \frac{1}{k} 
|( 1.5(\eta^{n+1}-\eta^n)- 0.5(\eta^n-\eta^{n-1}),\phi^{n+1}) \big|
\nonumber \\
\le&\, \frac{1}{k} ( 1.5 \|\eta^{n+1}-\eta^n\|_{L^2_x} 
+ 0.5 \|\eta^{n}-\eta^{n-1}\|_{L^2_x} ) 
\|\phi^{n+1}_h\|_{L^2_x}\label{cat2011} \\
\le &\, Ch^{r+1} \|u_t\|_{L^\infty_tH^{r+1}_x}
\|\phi^{n+1}_h\|_{L^2_x} 
\le 
\frac{1}{10}  \|\phi^{n+1}_h\|^2_{L^2_x}
+ Ch^{2r+2} \|u_t\|^2_{L^\infty_t H^{r+1}_x},
\label{mp001}
\end{align}
where Lemma\,\ref{lem_eta} is used in  \eqref{cat2011} and Young's inequality is used to earn \eqref{mp001}.

\item[(2)]
As for the third term on the right of \eqref{eqn_NSu_06}, 
we can get
\begin{align}
\big|(\nabla\cdot \phi^{n+1}_h,  p^{n+1}-\pi_Q (p^{n+1}) ) 
\big| &\le \|\nabla\cdot \phi^{n+1}_h \|_{L^2_x} \| p^{n+1} - \pi_Q(p^{n+1}) \|_{L^2_x}
\\
&\le \|\nabla \phi^{n+1}_h \|_{L^2_x} Ch^r   
\|p \|_{L^\infty_t H^r_x}
\label{cat1015} \\
&\le \frac{\nu}{4} \|\nabla\phi^{n+1}_h\|^2_{L^2_x} + 
\frac{Ch^{2r}\|p\|^2_{L^\infty_t H^{r}_x}}{\nu},
\label{mp002}
\end{align}
where the inequality \eqref{kun204} and the estimate \eqref{John_005} with $m=0$ are applied to get \eqref{cat1015} and Young's inequality is used in the last step. One way to avoid negative powers of $\nu$ is throwing derivatives to the pressure, thus we can get
\begin{align}
&\big|(\nabla\cdot \phi^{n+1}_h,  p^{n+1}-\pi_Q (p^{n+1}) )
\big| 
= \big|(\phi^{n+1}_h,  \nabla(p^{n+1}-\pi_Q (p^{n+1})) )
\big|  \nonumber \\
\le& \|\phi^{n+1}_h \|_{L^2_x} \| \nabla(p^{n+1} - \pi_Q(p^{n+1})) \|_{L^2_x}
\le \|\phi^{n+1}_h \|_{L^2_x} Ch^{r-1}   
\|p \|_{L^\infty_t H^r_x}
\label{cat1015b} \\
\le& \frac{1}{10} \|\phi^{n+1}_h\|^2_{L^2_x} 
+ Ch^{2r-2}\|p\|^2_{L^\infty_t H^{r}_x},
\label{mp002b}
\end{align}
where the estimate \eqref{John_005} is used in \eqref{cat1015b}.

Note if a divergence-free FE pair is used, then $\nabla\cdot \phi^{n+1}_h=0$ in $L^2$ sense and  
$(\nabla\cdot \phi^{n+1}_h,  p^{n+1}-\pi_Q (p^{n+1}) ) =0$.

\item[(3)]
As for the fourth term on the right of \eqref{eqn_NSu_06}, due to the general formula
\begin{equation}
\int_a^b (a-s)^2u_{ttt}(s)ds = (a-b)^2 u_{tt}(b) 
+ 2(a-b) u_t(b) + 2u(b) -2u(a),
\end{equation}
we obtain from \eqref{def_Rk} that 
\begin{equation}
R^{n+1}(k)=\frac{1}{k}\int_{t^n}^{t^{n+1}} (t^n-s)^2  u_{ttt}(s)ds
-\frac{1}{4k}\int_{t^{n-1}}^{t^{n+1}} (t^{n-1}-s)^2  u_{ttt}(s)ds.
\end{equation}
Thus, we can derive, by using the Cauchy-Schwarz inequality on the integrals on $s$,  the values of $\int_{t^{n}}^{t^{n+1}} (t^{n}-s)^4ds$ and $\int_{t^{n-1}}^{t^{n+1}} (t^{n-1}-s)^4ds$, and Fubini's Theorem to switch the integrals on $s$ and $x$, 
\begin{eqnarray}
\|R^{n+1}(k)\|^2_{L^2_x} 
%&\le& 
%\frac{C}{k^2} \left[
%\int_\Omega \left| \int_{t^n}^{t^{n+1}} (t^n-s)^2 u_{ttt}(s,x)ds \right|^2 dx
%+ 
%\int_\Omega \left| \int_{t^{n-1}}^{t^{n+1}} (t^{n-1}-s)^2 u_{ttt}(s,x)ds \right|^2 dx
%\right] \label{cat202501}\\
&\le& 
Ck^3 \left[
\int_{t^n}^{t^{n+1}} \int_\Omega |u_{ttt}(s,x)|^2dxds
+ 
\int_{t^{n-1}}^{t^{n+1}}\int_\Omega |u_{ttt}(s,x)|^2dx ds
\right]\label{202502}\\
&\le&
Ck^4 \|u_{ttt}\|_{L^\infty_t L^2_x}.
\label{202503}
\end{eqnarray}
%%%
This leads to the following estimate by applying Young's inequality,
\begin{align}
\big| (R^{n+1}(k),\phi^{n+1}_h) \big|
\le 
\frac{1}{10}\|\phi^{n+1}_h\|^2_{L^2_x}
+ Ck^4  \|u_{ttt}\|^2_{L^\infty_t L^2_x}.
 \label{mp003}
\end{align}

\end{enumerate}
%%%%%%%%%%%%%%%%%%%%%%%%%%%%%%%%%%%%%%%%%%%%%%
{\bf Step 3.} In $B(\phi^{n+1}_h)$, 
we replace $u^j_h$ with $u^j-\eta^j + \phi^j_h$ for $j=n+1, n, n-1$ and obtain
\begin{align}
B(\phi^{n+1}_h) 
=&\, -(\text{NL}(u^{n+1},\eta^{n+1}), \phi^{n+1}_h) 
-(\text{NL}(\eta^{n+1},s^{n+1}_h), \phi^{n+1}_h) 
+(\text{NL}(s^{n+1}_h,\phi^{n+1}_h), \phi^{n+1}_h) 
\nonumber \\
&\,
+(\text{NL}(\phi^{n+1}_h,s^{n+1}_h),\phi^{n+1}_h) 
+(\text{NL}(\phi^{n+1}_h,\phi^{n+1}_h),\phi^{n+1}_h).
\label{hai001}
\end{align}
Below is the estimate of the 5 terms on the right of \eqref{hai001}.
\begin{enumerate}
\item[(3A)] 
Estimate of 
$(\text{NL}(u^{n+1},\eta^{n+1}), \phi^{n+1}_h)$. 
Because $u^{n+1}$ is exact solution, 
$\nabla\cdot(u^{n+1})=0$ almost everywhere in $\Omega$. 
Thus $\text{NL}(u^{n+1},\eta^{n+1}) 
=(u^{n+1}\cdot \nabla) \eta^{n+1}$.
Therefore, 
\begin{align}
& | (\text{NL}(-u^{n+1},\eta^{n+1}), \phi^{n+1}_h) | 
\le  \|u^{n+1}\|_{L^\infty_x} 
 \|\nabla \eta^{n+1}\|_{L^2_x}
 \|\phi^{n+1}_h\|_{L^2_x} \nonumber \\
\le & \|u\|_{L^\infty_tH^2_x} 
Ch^r  \|u\|_{L^\infty_tH^{r+1}_x}
 \|\phi^{n+1}_h\|_{L^2_x}
\le 
\frac{1}{10}   \|\phi^{n+1}_h\|^2_{L^2_x}
+ Ch^{2r}  \|u\|^4_{L_t^\infty H^{r+1}_x},
 \label{mp003a}
\end{align}
where \eqref{agmon} is used to obtain 
$\|u^{n+1}\|_{L^\infty_x} \le \|u\|_{L^\infty_t H^2_x}$, 
the estimate \eqref{John_001} to get  $ \|\nabla \eta^{n+1}\|_{L^2_x}\le C h^r \|u\|_{L^\infty_t H^{r+1}_x}$,  and $r\ge 2$ and  Young's inequality are used in the last step.

\item[(3B)] 
Estimate of $(\text{NL}(\eta^{n+1},s^{n+1}_h), \phi^{n+1}_h)$. 
\begin{align}
& \big| (\text{NL}(\eta^{n+1}, s_h^{n+1}), \phi^{n+1}_h)
\big|  
= \big| ((\eta^{n+1}\cdot\nabla)s_h^{n+1},\phi^{n+1}_h)
 + \frac{1}{2} ((\nabla\cdot \eta^{n+1}) s_h^{n+1},\phi^{n+1}_h) \big| \label{cat1104} \\
\le &\, \|\nabla s_h^{n+1}\|_{L^\infty_x}
\|\eta^{n+1}\|_{L^2_x}
\|\phi^{n+1}_h\|_{L^2_x} 
 + \frac{1}{2} \|s^{n+1}_h\|_{L^\infty_x}
 \|\nabla\cdot \eta^{n+1} \|_{L^2_x}
 \|\phi^{n+1}_h\|_{L^2_x} \nonumber \\
 \le&\, \big(C h^{r+1} \|u\|_{L^\infty_t H^{3}_x}
 \|u\|_{L^\infty_t H^{r+1}_x}
+ Ch^{r}  \|u\|_{L^\infty_t H^{2}_x}
\|u\|_{L^\infty_t H^{r+1}_x}
\big)
\|\phi^{n+1}_h\|_{L^2_x}
\label{cat0311}\\
\le&\, Ch^r \|u\|^2_{L^\infty_t H^{r+1}_x} 
\|\phi^{n+1}_h\|_{L^2_x}
\label{cat1011}\\
\le&\, \frac{1}{10} \|\phi^{n+1}_h\|^2_{L^2_x}
+  Ch^{2r}  \|u\|^{4}_{L^\infty_t H^{r+1}_x}.
\label{mp003b}
\end{align}
To achieve \eqref{cat0311}, we have used $ \|\nabla s_h^{n+1}\|_{L^\infty_x}
\le C\|\nabla u\|_{L^\infty_t L^\infty_x}
\le C\|u\|_{L^\infty_t H^3_x}$ 
due to \eqref{John_003} and \eqref{agmon},  
$ \|\eta^{n+1}\|_{L^2_x}
\le Ch^{r+1} \|u\|_{L^\infty_t H^{r+1}_x}$ 
and $\|\nabla\cdot \eta^{n+1}\|_{L^2_x}\le Ch^{r} \|u\|_{L^\infty_t H^{r+1}_x}$
due to \eqref{John_001}, 
and $\|s^{n+1}_h\|_{L^\infty_x}\le 
C \|u\|_{L^\infty_t H^2_x}$ 
due to \eqref{John_002}.
To get \eqref{cat1011}, we have used 
$\|u\|_{L^\infty_t H^2_x}
\le \|u\|_{L^\infty_t H^3_x}
\le \|u\|_{L^\infty_t H^{r+1}_x}$ when $r\ge2$ and $h<1$.
In the step \eqref{mp003b}, Young's inequality is applied.

In the div-free FEMs, the second term in \eqref{cat1104} vanishes because $\nabla\cdot \eta^{n+1}=0$
and the final upper bound 
in\eqref{mp003b} becomes 
$\frac{1}{10} \|\phi^{n+1}_h\|^2_{L^2_x}
+  Ch^{2r+2}  \|u\|^{4}_{L^\infty_t H^{r+1}_x}$.

\item[(3C)] $(\text{NL}(s^{n+1}_h,\phi^{n+1}_h), \phi^{n+1}_h) =0$. 

\item[(3D)] 
Estimate of $(\text{NL}(\phi^{n+1}_h,s^{n+1}_h), \phi^{n+1}_h)$. 
\begin{align}
&\,\big| (\text{NL}(\phi^{n+1}_h,s^{n+1}_h), \phi^{n+1}_h) \big|
= \big| 
((\phi^{n+1}_h\cdot\nabla) s^{n+1}_h,\phi^{n+1}_h)
+\frac12 ((\nabla\cdot\phi^{n+1}_h) s^{n+1}, \phi^{n+1}_h) \big| \nonumber\\
\le&\, \|\nabla s^{n+1}_h\|_{L^\infty_x} 
\|\phi^{n+1}_h\|^2_{L^2_x} 
 + \frac{1}{2} 
\|s^{n+1}_h\|_{L^\infty_x} 
\|\nabla\cdot \phi^{n+1}_h\|_{L^2_x}
\|\phi^{n+1}_h\|_{L^2_x}
\label{cat1013}
\\
\le&\, C \|u\|_{L^\infty_t H^{r+1}_x}
\|\phi^{n+1}_h\|^2_{L^2_x} 
 + C \|u\|_{L^\infty_t H^{r+1}_x} 
 \|\nabla \phi^{n+1}_h\|_{L^2_x}
\|\phi^{n+1}_h\|_{L^2_x}
\label{cat2015}\\
\le&\, C \bigg(  \|u\|_{L^\infty_tH^{r+1}_x} + \frac{\|u\|^2_{L^\infty_t H^{r+1}_x}}{\nu} \bigg) \|\phi^{n+1}_h\|^2_{L^2_x}
+ \frac{\nu}{4} \|\nabla\phi^{n+1}_h\|^2_{L^2_x}. 
\label{mp003d}
\end{align}
From \eqref{cat1013} to \eqref{cat2015} then to \eqref{mp003d},  the estimates \eqref{John_003} and \eqref{agmon} are employed to obtain $\|\nabla s^{n+1}_h\|_{L^\infty_x}\le C\|\nabla u\|_{L^\infty_t L^\infty_x} \le C \|u\|_{L^\infty_t H^3_x} \le C \|u\|_{L^\infty_t H^{r+1}_x}$ when $r\ge 2$, 
and \eqref{John_002} is used to obtain $\|s^{n+1}_h\|_{L^\infty_x}\le \|u\|_{L^\infty_t H^2_x} \le 
\|u\|_{L^\infty_t H^{r+1}_x}$. Young's inequality is applied from \eqref{cat2015} to \eqref{mp003d}.

Note for  the div-free FEMs, the viscosity terms on the right of \eqref{mp003d} vanish because $\nabla\cdot \phi^{n+1}_h=0$.

\item[(3E)] $(\text{NL}(\phi^{n+1}_h, \phi^{n+1}_h), \phi^{n+1}_h) =0$.
\end{enumerate}

{\bf Step 4}. 
Dropping $\|\phi^{n+1}_h - 2\phi^n_h + \phi^{n-1}_h\|^2_{L^2_x}$ on the left of \eqref{eqn_NSu_06} and replacing the last four terms by their upper bounds in \eqref{mp001}, \eqref{mp002}, \eqref{mp003}, \eqref{mp003a}, \eqref{mp003b}, \eqref{mp003d}, 
we obtain for the non-divergence free FEMs
\begin{align}
&\,\frac{1}{k} \big( 
\|\phi^{n+1}_h\|^2_{L^2_x}
 + \|2\phi^{n+1}_h-\phi^n_h\|^2_{L^2_x}
 \big)
+ \nu \|\nabla \phi^{n+1}_h\|^2_{L^2_x}
\le
\frac{1}{k} 
\big(
\|\phi^{n}_h\|^2_{L^2_x}
+ \|2\phi^{n}_h-\phi^{n-1}_h\|^2_{L^2_x}
\big) \nonumber\\
&+ \frac{\nu}{2} \|\nabla \phi^{n+1}_h\|^2_{L^2_x}
+
\Big(1+ C \Big(\|u\|_{L^\infty_t H^{r+1}_x} 
+ \frac{1}{\nu} \|u\|^2_{L^\infty_t H^{r+1}_x} \Big) \Big)
\|\phi^{n+1}_h\|^2_{L^2_x} 
\nonumber\\
&\,+
C \Big( h^{2r+2} \|u_t\|^2_{L^\infty_t H^{r+1}_x} 
+ \frac{h^{2r} \|p\|^2_{L^\infty_t H^{r}_x}}{\nu}
+ h^{2r}   \|u\|^{4}_{L^\infty_t H^{r+1}_x} 
+ k^4  \|u_{ttt}\|^2_{L^\infty_t L^2_x } 
\Big).
\label{hai009}
\end{align}

According to the inequality $\sum_{i=1}^m x_i^2 \le \left( \sum_{i=1}^m x_i \right)^2$ for any positive integer $m$ and positive real numbers $x_i, i=1, \cdots, m$,  the sum of the four terms in the last pair of parentheses on the right of \eqref{hai009} is bounded above by $M_2^2$,  where $M_2$ is defined in \eqref{def_M2}.
Then, by multiplying $k$ on both sides and rearranging terms,  we can rewrite \eqref{hai009} as 
\begin{align}
(1-kM_1) \|\phi^{n+1}_h\|^2_{L^2_x}
 + \|2\phi^{n+1}_h-\phi^n_h\|^2_{L^2_x}
 + \frac{k \nu}{2} \|\nabla\phi^{n+1}_h\|^2_{L^2_x}
\le
\|\phi^{n}_h\|^2_{L^2_x}
+  \|2\phi^{n}_h-\phi^{n-1}_h\|^2_{L^2_x} 
+ kC M_2^2,
\label{eqn_NSu_09}
\end{align}
where $M_1$ is defined in \eqref{def_M1}.

We further denote 
\begin{equation}
a^n\triangleq \|\phi^{n}_h\|^2_{L^2_x}
+ \|2\phi^{n}_h-\phi^{n-1}_h\|^2_{L^2_x}
+ \frac{k\nu}{2}
\|\nabla\phi^{n}_h\|^2_{L^2_x}.
\label{def_an}
\end{equation}
It is easy to check that when $k<1/M_1$,  the left side of \eqref{eqn_NSu_09} is bounded below by $(1-k M_1) a^{n+1}$, and the right side of \eqref{eqn_NSu_09}
is bounded above by 
$a^n + kCM_2^2 $. 
Then \eqref{eqn_NSu_09}  can be thrown into a recursive relation $a^{n+1}\le \beta a^n + \alpha$, 
where 
$\alpha = \frac{kC M^2_2}{1-kM_1}$ and 
$\beta = \frac{ 1 }{1-k M_1}.
$
This relation ultimately generates 
$a^{n} \le e^{n (\beta-1)} \big( a^1 + \frac{\alpha}{\beta-1} \big)$, which corresponds to
$
a^{n} \le 
\exp( \frac{nk M_1} { 1-k }) 
\cdot
(a^1 + \frac{C M^2_2}{M_1}).
$
By using $nk\le T$ for $n\le \lfloor \frac{T}{k} \rfloor$
and  dropping 
$\|2\phi^{n}_h-\phi^{n-1}_h\|^2_{L^2_x}$ in $a^n$, we obtain
\begin{align}
\|\phi^n_h\|^2_{L^2_x}
 + \frac{k\nu}{2} \|\nabla \phi^n_h\|^2_{L^2_x} 
\le  
\exp{ \Big( \frac{T M_1}{1-kM_1}  \Big) }
\cdot 
\Big(a^1 + \frac{ CM^2_2}{M_1}
  \Big).
\label{error_anal_01}
\end{align}
As for the full error $u^n - u^n_h=(u^n-s^n_h) + \phi^n_h$, we have
\begin{align}
&\, \|u^n-u^n_h\|^2_{L^2_x}
 + \frac{k\nu}{2} \| \nabla(u^n-u^n_h)\|^2_{L^2_x}
\nonumber \\
\le&\,
 \|\phi^n_h\|^2_{L^2_x}
 + \frac{k\nu}{2}  \|\nabla \phi^n_h\|^2_{L^2_x} 
+ \|u^n- s^n_h\|^2_{L^2_x}
+  \frac{k\nu}{2} \| \nabla(u^n-s^n_h)\|^2_{L^2_x}
\\
\le&\, 
  \|\phi^n_h\|^2_{L^2_x}
 + \frac{k\nu}{2} \|\nabla \phi^n_h\|^2_{L^2_x} 
 + Ch^{2r+2} \|u\|^2_{L^\infty_t H^{r+1}_x}
 + C\nu k h^{2r} \|u\|^2_{L^\infty_t H^{r+1}_x} 
\label{error_anal_02}
\end{align}
where the inequality \eqref{John_001} is applied in the last step.
Finally, the error estimate \eqref{error_anal_FI} is obtained by combining \eqref{error_anal_01} and \eqref{error_anal_02}.

In the case of the div-free FEMs, some terms in \eqref{hai009} are changed or deleted based on the  analysis in Step 3, but the same estimate  \eqref{eqn_NSu_09} holds with  $M_1$  and $M_2$  replaced with  $M^{div0}_1$ defined in \eqref{def_M1} and  $M^{div0}_2$ in \eqref{def_M2div0}, respectively. The analysis after  \eqref{eqn_NSu_09} is the same as the case with the non-div-free FEMs. 
\end{pf}

\vskip 0.4cm

\section*{Declarations}

\noindent {\bf Funding} \quad X. Zheng was partially supported by NSF grant DMS-2309747.
K. Zhao was partially supported by the Simons Foundation Collaboration Grant for Mathematicians No.\,413028. 
J. Wu was partially supported by  NSF grants DMS-2104682 and DMS-2309748.
W. Hu was partially supported by the  NSF grants DMS-2005696 (previously DMS-1813570), DMS-2111486 and DMS-2205117.

\noindent {\bf Conflict of interest} \quad
The authors declare no competing interests.

\noindent {\bf Author contributions} \quad
Conceptualization and Methodology: X. Zheng, J. Wu, and D. Du. Software and Visualization:  X. Zheng.
Formal Analysis: all authors.
Writing: all authors.

\vskip 0.4cm

%%%%%%%%%%%%%%%%%%%%%%%%%%%%%%%%%%%
%\bibliographystyle{spmpsci}      % mathematics and physical sciences
%\bibliographystyle{spphys}       % APS-like style for physics
\bibliographystyle{plain}
\bibliography{refer}   % name your BibTeX data base

\end{document}